\documentclass[11pt,leqno]{amsart}
\usepackage{amsthm,amsfonts,amssymb,amsmath,oldgerm}
\numberwithin{equation}{section}
\usepackage{fullpage}
\usepackage{amsmath}
\usepackage{amsfonts}
\usepackage{amssymb,bm}
\usepackage{setspace}
\usepackage{stmaryrd}
\usepackage{mathrsfs}
\usepackage{fancyhdr}
\usepackage{graphicx}
\usepackage{subcaption}
\usepackage[linkcolor=blue,colorlinks=true]{hyperref}
\usepackage{psfrag}
\usepackage{listings} 
\usepackage{pgfplots}
\usepackage[top=20mm,bottom=20mm,left=25mm,right=25mm,a4paper]{geometry}
\usetikzlibrary{shapes.geometric}

\renewcommand\d{\partial}
\DeclareMathOperator\dD{d}
\DeclareMathOperator\eD{e}

\def\eps{\varepsilon }
\DeclareMathOperator{\Id}{Id}

\DeclareMathOperator{\Tr}{Tr}

\DeclareMathOperator{\Div}{div}
\DeclareMathOperator{\Curl}{curl}

\newcommand\br{\begin{remark}}
\newcommand\er{\end{remark}}
\newcommand\bp{\begin{pmatrix}}
\newcommand\ep{\end{pmatrix}}
\newcommand{\be}{\begin{equation}}
\newcommand{\ee}{\end{equation}}
\newcommand{\ba}[1]{\begin{array}{#1}}
\newcommand{\ea}{\end{array}}
\newcommand\ds{\displaystyle}
\newcommand\nn{\nonumber}
\newcommand{\beg}{\begin{example}}
\newcommand{\eeg}{\end{exaplem}}
\newcommand{\bpr}{\begin{proposition}}
\newcommand{\epr}{\end{proposition}}
\newcommand{\bt}{\begin{theorem}}
\newcommand{\et}{\end{theorem}}
\newcommand{\bc}{\begin{corollary}}
\newcommand{\ec}{\end{corollary}}
\newcommand{\bl}{\begin{lemma}}
\newcommand{\el}{\end{lemma}}
\newcommand{\bd}{\begin{definition}}
\newcommand{\ed}{\end{definition}}
\newcommand{\brs}{\begin{remarks}}
\newcommand{\ers}{\end{remarks}}

\newtheorem{theorem}{Theorem}[section]
\newtheorem{proposition}[theorem]{Proposition}
\newtheorem{corollary}[theorem]{Corollary}
\newtheorem{lemma}[theorem]{Lemma}
\newtheorem{remark}[theorem]{Remark}
\newtheorem{definition}[theorem]{Definition}

\newtheorem{example}[theorem]{Example}

\newcommand\R{\mathbf R}

\newcommand{\N}{\mathbf N}
\newcommand{\Z}{\mathbf Z}

\newcommand\bA{{\mathbf A}}
\newcommand\bB{{\mathbf B}}

\newcommand\bE{{\mathbf E}}
\newcommand\bF{{\mathbf F}}

\newcommand\bJ{{\mathbb J}}
\newcommand\bK{{\mathbf K}}
\newcommand\bL{{\mathbf L}}

\newcommand\bS{{\mathbf S}}

\newcommand\bU{{\mathbf U}}

\newcommand\bY{{\mathbf Y}}
\newcommand\bZ{{\mathbf Z}}

\newcommand\bfF{\bF}
\newcommand\bfU{\bU}

\newcommand\bfa{{\mathbf a}}
\newcommand\bfb{{\mathbf b}}

\newcommand\bfj{{\mathbf j}}

\newcommand\bfu{{\mathbf u}}
\newcommand\bfv{{\mathbf v}}

\newcommand\bfx{{\mathbf x}}
\newcommand\bfy{{\mathbf y}}
\newcommand\bfz{{\mathbf z}}
\newcommand\bfZ{{\mathbf Z}}

\newcommand\bfeta{{\boldsymbol \eta}}

\newcommand\bfchi{{\boldsymbol \chi}}
\newcommand\Sig{{\boldsymbol \Sigma}}
\newcommand\bfsigma{\bm{\sigma}}

\newcommand\urho{{\underline \rho}}

\newcommand\cA{{\mathcal A}}

\newcommand\cC{{\mathcal C}}

\newcommand\cE{{\mathcal E}}

\newcommand\cL{{\mathcal L}}
\newcommand\cM{{\mathcal M}}

\newcommand\cO{{\mathcal O}}

\newcommand\cU{{\mathcal U}}
\newcommand\cV{{\mathcal V}}
\newcommand\cW{{\mathcal W}}
\newcommand\cX{{\mathcal X}}
\newcommand\cY{{\mathcal Y}}
\newcommand\cZ{{\mathcal Z}}

\newcommand\teta{{\widetilde \eta}}

\newcommand\tby{\widetilde{\bfy}}

\newcommand\tbA{\widetilde{\bA}}

\DeclareMathOperator{\beD}{\mathbf{e}}
\newcommand\eDx{\beD_x}
\newcommand\eDy{\beD_y}
\newcommand\eDr{\beD_r}
\newcommand\eDz{\beD_z}
\newcommand\eDpar{\beD_\mypar}
\newcommand\eDperp{\beD_\perp}
\newcommand\eDa{\beD_a}
\newcommand\eDb{\beD_b}
\newcommand\vpar{v_\mypar}
\newcommand\bvperp{\bfv_\perp}
\newcommand\wperp{w_\perp}
\newcommand\eperp{w_\perp}
\newcommand\e{w}
\newcommand\muperp{\mu_\perp}
\newcommand\Epar{E_\mypar}
\newcommand\Eperp{\bE_\perp}
\newcommand\Er{E_r}
\newcommand\Ez{E_z}
\newcommand{\Trperp}[1]{\Tr_\perp^{#1}}

\newcommand{\ZGC}{\bfZ_{\rm gc}}
\newcommand{\GC}{\bfx_{\rm gc}}
\newcommand{\eGC}{w_{\rm gc}}
\newcommand{\vGC}{v_{\rm gc}}
\newcommand{\uGC}{\bfu_{\rm gc}}
\newcommand{\muGC}{\mu_{\rm gc}}
\newcommand{\EcB}{\bU_{\bE\times\bB}}
\newcommand{\gradB}{\bU_{\nabla B\times\bB}}
\newcommand{\curvB}{\bU_{\rm curv}}
\newcommand{\rotB}{\bU_{\Curl\eDpar}}
\newcommand{\dtB}{\bU_{\d_t}}
\newcommand{\Ud}{\bU_{\rm drift}}

\newcommand{\tobs}{T_{\rm obs}}

\newcommand{\mypar}{{\mkern3mu\vphantom{\perp}\vrule depth 0pt\mkern2mu\vrule depth 0pt\mkern3mu}}

\title[The Vlasov equation with strong magnetic field]{Asymptotics of the three dimensional Vlasov equation in the large magnetic field limit}

\author{Francis Filbet}
\address{
Universit\'e de Toulouse III \& IUF,
UMR5219, Institut de Math\'ematiques de Toulouse,
118, route de Narbonne;
F-31062 Toulouse Cedex, FRANCE}
\email{{\tt francis.filbet@math.univ-toulouse.fr}}
\thanks{FF was supported by the EUROfusion Consortium and has received funding
from the Euratom research and training programme 2014-2018 under grant
agreement No 633053. The views and opinions expressed herein do not
necessarily reflect those of the European Commission.\\
}

\author{L.~Miguel Rodrigues}
\address{
Univ Rennes \& IUF, CNRS, IRMAR - UMR 6625, F-35000 Rennes, FRANCE}
\email{{\tt luis-miguel.rodrigues@univ-rennes1.fr}}
\thanks{Research of LMR has received funding from the city of Rennes.}

\begin{document}

\begin{abstract}
We study the asymptotic behavior of solutions to the  Vlasov equation in the presence of a strong external magnetic field. In particular we provide a mathematically rigorous derivation of the guiding-center approximation in the general three dimensional setting under the action of large inhomogeneous magnetic fields. First order corrections are computed and justified as well, including electric cross field, magnetic gradient and magnetic curvature drifts. We also treat long time behaviors on two specific examples, the two dimensional case in cartesian coordinates and a toroidal axi-symmetric geometry, the former for expository purposes. Algebraic manipulations that underlie concrete computations make the most of the linearity of the stiffest part of the system of characteristics instead of relying on any particular variational structure. At last, we analyze a smoothed Vlasov-Poisson system thus show how our arguments may be extended to deal with the nonlinearity arising from self-consistent fields.
\end{abstract}

\date{\today}
\maketitle




\tableofcontents


\section{Introduction}\label{s:introduction}

Since fusion configurations involve very hot plasmas, they typically
require a careful design to maintain fast moving particles inside the
core of the device on sufficiently long times. In the magnetic
confinement approach
\cite{bellan_2006_fundamentals,chen_introduction,freidberg2008plasma,haz_mei_03,miyamoto_2006_plasma,
  piel2010plasma}, in particular in tokamak plasmas, a strong external
field is applied to confine the plasma by enforcing the oscillatory nature of the fast motions.

Various models are in use to describe such phenomena. In the kinetic modeling, the unknowns are the number densities of particles, $f\equiv f(t,\bfx,\bfv)$ depending on time $t\geq 0$, position $\bfx\in\Omega\subset \R^3$ and velocity $\bfv\in\R^3$. Such kinetic models provide an appropriate description of turbulent transport in a fairly general context, but in fusion configurations their numerical simulations require to solve a stiff six-dimensional problem, leading to a huge computational cost. To bypass this obstacle, it is classical --- see for instance \cite{Garbet-et-al_2010} --- to use reduced asymptotic models that describe only the slowest part of the plasma dynamics hence effectively reducing both the stiffness of the problem and the number of variables (since fastest variables are omitted). Over the years, due to its rich and fundamental nature, the physically-based derivation of such models has grown as a --- still very active --- field of its own, often referred to as gyrokinetics. Besides the already mentioned general monographs \cite{bellan_2006_fundamentals,chen_introduction,freidberg2008plasma,haz_mei_03,miyamoto_2006_plasma, piel2010plasma}, the reader may consult \cite{Krommes,bri_hahm_07,Matteo-PhD,Scott_gyrokinetic,PDFF} and references therein as more specialized entering gates to the field.

Despite considerable efforts in recent years, concerning mathematically rigorous derivations from collisionless\footnote{See for instance \cite{herda_2016_massless,herda_2016_anisotropic} and references therein for an introduction to the corresponding collisional issues.} kinetic equations,  the state of art is such that one must choose between linear models that neglect couplings due to self-consistent fields or nonlinear ones set in a deceptively simple geometry. See for instance the introductions and bibliographies of \cite{HanKwan_PhD,Lutz_PhD,Herda_PhD} for relatively recent panoramas on the question. For instance, for the kind of problem considered here, on the nonlinear side of the literature the most significant mathematical result --- which requires a careful analysis --- is restricted to a two-dimensional setting with a constant magnetic field and interactions described through the Poisson equation, and yet validates only half\footnote{The nontrivial half, however. This is possible there only because a very specific geometric cancellation uncouples part of the slow dynamics from the remaining one, which is expected to be slaved to it. See however the recent \cite{Bostan_2D-VP} for a more complete model, derived under more stringent assumptions.} of the slow dynamics; see \cite{laure0}, building on \cite{gol_lsr_99} and recently revisited in \cite{Miot-2D-gyrokinetic}. 

We consider here a plasma confined by a strong unsteady inhomogeneous magnetic field without any a priori geometric constraint but, in order to allow for such a generality, in most of the present paper\footnote{See however Section~\ref{s:nl} where we analyze a smoothed Vlasov-Poisson system.} we do neglect effects of self-consistent fields. The plasma is thus entirely modeled with a scalar linear kinetic equation, where the unknown is one of the number densities of particles. The approach that we follow focuses on the characteristic equations associated with the kinetic conservation law. By itself the study of those equations may follow the classical roadmap of the averaging of ordinary differential equations, as expounded in \cite{Bogoliubov-Mitropolsky_oscillations,Sanders-Verhulst-Murdock_averaging}. Yet, here, beyond the body of work already required to follow this road in usual ODE problems, a careful track of the dependence of averaging estimates on initial data, living here in an unbounded phase space, is necessary so as to derive asymptotics for the solutions of the original partial differential equations problem.

To be more specific, the Lorentz force term in our original nondimensionalized kinetic equation is scaled by a large parameter, $1/\eps$, where $\eps$ stands for the typical cyclotron period, {\it i.e.} the typical rotation period of particles about a  magnetic field line (or Larmor rotation). The dynamical time scales we focus on are in any case much larger than the cyclotron period and we establish asymptotic descriptions in the limit $\eps\to0$. As is classical in the field, we distinguish between short-time scales that are $\cO(1)$ with respect to $\eps$, and long time scales that are $\sim1/\eps$ in the limit $\eps\to0$. Correspondingly, slow dynamics refer to dynamics where typical time derivatives are at most of order $\cO(1)$ on short-time scales, and at most of order $\cO(\eps)$ on long-time scales so that on long time scales two kinds of fast dynamics may co-exist, principal ones at typical speed of order $1/\eps$ and subprincipal ones at typical speed of order $1$; see for instance \cite{cheve2} for a description of those various oscillations in a specific class of axi-symmetric geometries, without electric field and with a magnetic field nowhere toroidal and whose angle to the toroidal direction is also independent of the poloidal angle. With this terminology in hands, our results may be roughly stated as the identification and mathematical proofs of
\begin{enumerate}
\item a second-order --- that is, up to $\cO(\eps^2)$ --- description of the slow dynamics on short time scales but in arbitrary geometry;
\item a first-order description of the slow dynamics on long time scales but in an axi-symmetric geometry with a magnetic field everywhere poloidal and an electric field everywhere orthogonal to the magnetic field.
\end{enumerate}
The geometry of the latter is very specific and the proof of such a description is mostly carried out here to illustrate that the short-time second-order description contains all the ingredients to analyze long-time dynamics at first-order. Note that in any case, on long-time scales some restrictions are indeed necessary to ensure that sub-principally fast dynamics do not prevent long-time confinement and are of oscillatory type so that the issue of the identification of a long-time slow dynamics becomes meaningful. In Section~\ref{s:nl} we also prove a second-order description of the dynamics driven by a smoothed Vlasov-Poisson system, hence allowing for both nonlinear self-consistent effects and arbitrary geometry, but we restrict there to initial data that are well-prepared in the sense that their initial dependence on fast angles is weak.

A key feature of our analysis that underpins a treatment of
essentially arbitrary fields is that we make no explicit use of any
geometric structure, neither Hamiltonian (see for instance the pioneer
work of R. G. Littlejohn \cite{littleJ1, littleJ2, littleJ3} and later
\cite{Benettin-Sempio,FrenodLutz_geometrical_gyro-kinetic}) nor Lagrangian (see \cite{Possanner}). The main role of these structures in the averaging process is to ease the identification of terms that are asymptotically irrevelant as time-derivatives of small terms. Instead, in the present contribution this explicit identification hinges heavily on the linearity of principal oscillations. As an upset, besides generality, we gain the freedom to use change of variables that are also arbitrary and to focus on slow variables instead of carrying geometric constraints all along. 

A key motivation for our methodology is that in the design of well-adapted numerical schemes, that capture the slow part of the dynamics even with discretization meshes too rough to compute stiff scales, one might correspondingly aim at large classes of schemes of arbitrary order; see for instance \cite{Lee,FR1,FR2}. Likewise our choice of studying first characteristics instead of using directly partial differential equations techniques and our will to prove error estimates echoes the particle-in-cell methodology and its numerical analysis. Alternative PDE-based methods include most notably two-scale convergence analyses \cite{fre_son_97,fre_son_98} and filtering techniques hinging on ergodic von Neumann's theorem \cite{Bostan_transport,bostan_10}. Two main advantages of going through characteristics are that the limiting partial differential equation is by construction a conservation law for a density distribution and that increasing the order of description may be carried out merely by continuing the argument used to identify the leading order. We benefit from the latter to \emph{prove} for the first time a second-order description in full generality.


\section{Definitions and main results}
\label{s:results}

We consider the kinetic PDE 
\be\label{eq:vlasov}
\d_t f^\eps\,+\,\Div_\bfx(f^\eps\,\bfv)
\,+\,\Div_\bfv\left(f^\eps\,\left(\frac{\bfv\wedge \bB(t,\bfx)}{\eps}\,+\,\bE(t,\bfx)\right)\right)\,=\,0
\ee
and its characteristic flow encoded by the following ODEs
\be
\label{eq:xv}
\left\{
\ba{l}
\ds\frac{\dD\bfx}{\dD t}\,=\,\bfv\,,
\\[0.9em]
\ds\frac{\dD\bfv}{\dD t}\,=\,\frac{\bfv\wedge \bB(t,\bfx)}{\eps} \,+\,\bE(t,\bfx)\,,
\ea
\right.
\ee
where $\wedge$ denotes the standard vector product on $\R^3$, $\bB$ stands for the external magnetic field, $\bE$ for the external electric field.

As a preliminary we begin by recalling the classical link between \eqref{eq:vlasov} and \eqref{eq:xv} and making explicit how it can be used to analyze the slow part of the dynamics hidden in the stiff \eqref{eq:vlasov}.

\subsection{From ODEs to PDEs}\label{s:PDE}

Throughout the present contribution we shall use the following notational
conventions. We denote $\Psi_*(\mu)$ the push-forward of $\mu$ by $\Psi$,
which can be defined for instance when $\mu$ is a distribution and $\Psi$
is a smooth proper map by, for any test-function $\varphi$, 
$$
\int_{B} \varphi  \dD \Psi_*(\mu)\,=\, \int_{\Psi^{-1}(B)} \varphi\circ \Psi \dD \mu.
$$
When considering flows associated with ODEs, $\Phi(t,s,\bfy)$ denotes the value at time $t$ of the solution starting from $\bfy$ at time $s$ and the associated map is $\bfy\mapsto \Phi(t,s,\bfy)$. In particular the solution to \eqref{eq:vlasov} starting from $f_0$ at time $0$ is given at time $t$ by $\Phi(t,0,\cdot)_*\,(f_0)$ where $\Phi$ is the flow associated with \eqref{eq:xv}.

For general purpose we shall state an abstract proposition, almost tautological, converting estimates on characteristics into estimates on densities. First, to enlighten the meaning of the following statement, identifying measures with their densities, we recall that the ``value'' at $a$ of $\cA_*(\mu)$ the push-forward of $\mu$ by $\cA$ is essentially the average of $\mu$ on the level set $\cA^{-1}(\{a\})$. Indeed for any function $f$ at any regular value $a$ of $\cA$
$$
\cA_*\,(f)(a)\,=\,\int_{\cA^{-1}(\{a\})}\,f(\bfy)\,\frac{\dD \sigma_a(\bfy)}{\sqrt{\det(\dD\cA(\bfy)(\dD \cA(\bfy))^*)}}\,,
$$
where $\sigma_a$ denotes the surface measure on $\cA^{-1}(\{a\})$, $\dD$ denotes the differential operator and ${}^*$ the adjoint operator. For instance if $\bfy=(y_1,y_2)\in\R^2$, then with $\cA(\bfy)=\|\bfy\|=\sqrt{y_1^2+y_2^2}$, 
$$
\cA_*\,(f)(r)\,=\,\int_0^{2\pi}\,f(r\,\beD(\theta))\,r\,\dD \theta\,,
$$
where $\beD(\theta)=(\cos(\theta),\sin(\theta))$, whereas with $\cA(\bfy)=\tfrac12\|\bfy\|^2$,
$$
\cA_*(f)(e)\,=\,\int_0^{2\pi}\,f(\sqrt{2\,e}\,\beD(\theta))\dD \theta\,.
$$
It turns out that the correct way to ``average" the stiff equation \eqref{eq:vlasov} is precisely to push $f$ by a map $\cA$ defining a complete\footnote{So that an uncoupled system is obtained in closed form (at the required order).} set of slow variables.

\bpr
\label{p:ODEtoPDE}
Let $\Phi$ and $\Phi_{\rm slow}$ be flows associated with respective ODEs 
$$
\frac{\dD\bfy}{\dD t}=\cX(t,\bfy)\qquad\textrm{and}\qquad 
\frac{\dD\bfa}{\dD t}=\cX_{\rm slow}(t,\bfa)
$$
and assume that there exist time-dependent slow maps $\cA(t,\cdot)$ and weights $\cM(t,\cdot)$ such that for a.e. $t\geq0$,
$$
\|\cA(t,\Phi(t,0,\cdot))-\Phi_{\rm slow}(t,0,\cA(0,\cdot))\|\,\leq\,\cM(t,\cdot)\,.
$$
Then if $f$ solves
$$
\d_t f\,+\,\Div_\bfy (\cX\,f)\,=\,0,
$$
with initial data a measure $f_0$ and $F(t,\cdot)=\cA(t,\cdot)_*\,f(t,\cdot)$ is the push-forward of $f$ by the slow map $\cA$ then for a.e. $t\geq0$
$$
\|F(t,\cdot)-G(t,\cdot)\|_{\dot{W}^{-1,1}}\,\leq\,\int
\cM(t,\cdot)  \dD|f_0|\,,
$$
where $G$ solves 
$$
\d_t G \,+\,\Div_\bfa (\cX_{\rm as} \,G)\,=\,0,
$$
with initial data $F_0:=\cA(0,\cdot)_* \,f_0$.
\epr

In the former we have denoted $\dot{W}^{-1,1}$ the dual of $\dot{W}^{1,\infty}$. Incidentally we observe that the distance on $\dot{W}^{-1,1}$ coincides with the $1$-Wasserstein distance from optimal transportation. Explicitly 
$$
\|\mu\|_{\dot{W}^{-1,1}}\,=\,\sup_{\|\nabla\varphi\|_{L^\infty}\leq
  1}\int \varphi\dD \mu\,.
$$
Here and throughout $L^p$ denotes the classical Lebesgue space of index $p$, $W^{s,p}$ and $\dot{W}^{s,p}$ their corresponding Sobolev spaces at regularity $s$, respectively in inhomogeneous and homogeneous versions. Associated (semi-)norms are denoted $\|\cdot\|_{L^p}$, $\|\cdot\|_{W^{s,p}}$ and $\|\cdot\|_{\dot{W}^{s,p}}$.

\begin{proof}
This stems readily from 
$$
F(t,\cdot)\,=\,\cA(t,\Phi(t,0,\cdot))_*\,(f_0)\,,\qquad
\qquad
G(t,\cdot)=\Phi_{\rm slow}(t,0,\cA(0,\cdot))_*\,(f_0)\,,
$$
and
$$
\|\varphi\circ\cA(t,\Phi(t,0,\cdot))-\varphi\circ \Phi_{\rm slow}(t,0,\cA(0,\cdot))\|
\leq \|\nabla\varphi\|_{L^\infty}\,\cM(t,\cdot)\,.
$$
\end{proof}

Note that in the foregoing statement, for readability's sake, we have deliberately left domains in time, original variables and slow variables, unspecified. However this may be straightened by classical ODE considerations, notably when fields are continuous, and locally Lipschitz in respectively $\bfy$ and $\bfa$ and either the support of $f_0$ is compact or involved vector-fields grow at most linearly.

\subsection{Slow variables and first-order asymptotics}

Getting back to our concrete system we begin our identification of a slow dynamics.

First, as is classical, we split the magnetic field $\bB$ as 
$$
\bB(t,\bfx)\,=\,B(t,\bfx)\,\eDpar(t,\bfx),
$$
with $B(t,\bfx)=\|\bB(t,\bfx)\|$. Accordingly we define, for any $\bfx\in\R^3$ and any time $t$, the linear operator $\bJ(t,\bfx)$ as
\be
\label{def:J}
\bJ(t,\bfx) \,\bfa\,=\, \bfa\wedge \eDpar(t,\bfx)\,.
\ee

The direction of the magnetic field plays a very special role and it is expedient to introduce for velocities an associated decomposition into parallel and orthogonal components
$$
\left\{
\ba{l}
\ds\vpar(t,\bfx,\bfv) \,=\,\langle\bfv\,,\,\eDpar(t,\bfx)\rangle\,,
\\[0.9em]
\ds\bvperp(t,\bfx,\bfv)=\bfv-\vpar(t,\bfx,\bfv)\,\eDpar(t,\bfx)
\ea\right.
$$
and similarly for the electric field $\bE$,
$$
\left\{
\ba{l}
\ds\Epar(t,\bfx)=\left\langle\bE(t,\bfx)\,,\,\eDpar(t,\bfx)\right\rangle,
\\[0.9em]
\ds\Eperp(t,\bfx)=\bE(t,\bfx)-\Epar(t,\bfx,\bfv)\,\eDpar(t,\bfx)\,,
\ea\right.
$$
where hereabove $\langle\cdot,\cdot\rangle$ denotes the canonical Euclidean scalar product, and below $\|\cdot\|$ denotes the associated Euclidean norm.

From system~\eqref{eq:xv} it is clear that at least one component out of the six-dimensional $(\bfx,\bfv)$ must obey a dynamics forcing oscillations of amplitude of typical size $1$ and typical frequency $1/\eps$. However at typical size $1$ a five-dimensional slow dynamics survives. This is already suggested by the fact that one may derive from \eqref{eq:xv} for the slow variables $(\bfx,\vpar,\eperp)$, 
\be
\label{eq:1}
\left\{
\ba{l}\ds
\frac{\dD\bfx}{\dD t} \,=\,\ds\bfv,
\\[0.9em]
\ds\frac{\dD\vpar}{\dD t}\,=\,\Epar(t,\bfx)\,+\,\left\langle \bvperp ,
\d_t\eDpar(t,\bfx)+\dD_{\bfx}\eDpar(t,\bfx)\bfv \right\rangle,
\\[0.9em]
\ds\frac{\dD\eperp}{\dD t}\,=\, \left\langle \Eperp(t,\bfx) \,-\,\vpar\,\left(
\d_t\eDpar(t,\bfx)+\dD_{\bfx}\eDpar(t,\bfx)\,\bfv\right),  \bvperp\right\rangle,
\ea \right.
\ee
where $\eperp=\tfrac12\|\bvperp\|^2$ and we have used the shorthand $\vpar(t)$ for
$\vpar(t,\bfx(t),\bfv(t))$ and similarly for $\bvperp$. 

Our goal is to identify such a slow dynamics, uncoupled from fast oscillations. Roughly speaking, since $\bvperp$ is expected to weakly converge to zero when $\eps$ goes to zero, at leading order the only issue is to identify the asymptotic behavior of quadratic terms in $\bvperp$ in \eqref{eq:1}. It turns out that those are responsible for the apparition of terms $\e\,\Div_\bfx\eDpar$ in the asymptotic model, set on a reduced phase space, where slow variables $\bfZ=(\bfy,v,\e)$ live. Introducing the limiting vector field
\be
\label{V:0}
\cV_0(t,\bfZ)
\,=\,\bp\ds
v\,\eDpar(t,\bfy)\\[1em]\ds
\Epar(t,\bfy)+\e\,\Div_\bfx\eDpar(t,\bfy)\\[1em]\ds
-v\,\e\,\Div_\bfx\eDpar(t,\bfy)
\ep
\ee
we may state our first significant result.

\bt\label{th:1}
Let $\bE\in W^{1,\infty}$ and $\bB$ be such that $1/B\in W^{1,\infty}$ and $\eDpar\in W^{2,\infty}$. There exists a constant $C$ depending polynomially on $\|\bE\|_{W^{1,\infty}}$,  $\|B^{-1}\|_{W^{1,\infty}}$ and $\|\eDpar\|_{W^{2,\infty}}$ such that if $f^\eps$ solves \eqref{eq:vlasov} with initial data a nonnegative density\footnote{Results would equally well hold with measures of arbitrary sign, but we stick to densities to provide nicer integral formulations for push-forwards when available, and to nonnegative densities to remove absolute values in error bounds.} $f_0$, then $F^\eps$ defined by
$$
F^\eps(t,\bfx,v_\mypar,\eperp)=\int_{\bS_{t,\bfx}}\,f^\eps(t,\bfx,v_\mypar\,\eDpar(t,\bfx)\,+\,\sqrt{2\,\eperp}\ \widehat{\beD})\ \dD \sigma_{t,\bfx}(\widehat{\beD}),
$$
with $\bS_{t,\bfx}=\{\eDpar(t,\bfx)\}^\perp\cap\bS^2$ and $\sigma_{t,\bfx}$ its canonical line-measure, satisfies for a.e. $t\geq0$
$$
\|F^\eps(t,\cdot)-G(t,\cdot)\|_{\dot{W}^{-1,1}}\leq C\,\eps\,e^{C\,t^4}\,\int_{\R^3\times\R^3} e^{C\,t\,\|\bfv\|^3}\,\|\bfv\|\,(1+\|\bfv\|^2)\,f_0(\bfx,\bfv)\,\dD\bfx\,\dD\bfv,
$$
where $G$ solves 
\be
\label{eq:1st}
\d_t G\,+\,\Div_\bfZ \left(\cV_0\,G \right) \,=\,0,
\ee
with $\cV_0$ given by (\ref{V:0}) and the initial datum $G_0$ is
$$
G_0(\bfZ)=\int_{\bS_{0,\bfy}}\,f_0(\bfy,v\,\eDpar(0,\bfy)\,+\,\sqrt{2\,\e}\ \widehat{\beD})\ \dD \sigma_{0,\bfy}(\widehat{\beD})\,.
$$
\et
Theorem~\ref{th:1} is proved in Section~\ref{sec:proofTh1} ; a nonlinear counterpart for a system of Vlasov-Poisson type is both stated and proved in Section~\ref{s:1st-nl}. 

The underlying vector field $\cV_0$ of the asymptotic model being divergence-free, many conservation laws already come as consequences of the asymptotic model. Yet as we state below a few more may be obtained if one assumes classical extra structure on electromagnetic fields.

\bpr
\label{prop:cons1}
Assume that $\bE=-\nabla_\bfx\phi$ where the couple $(\phi,B)$ does not depend on time and suppose that the confining magnetic field satisfies the Gauss' law 
$$
\Div_\bfx\bB\,=\, 0.
$$
Then solutions to the asymptotic model \eqref{eq:1st} satisfy
\begin{itemize}
\item the conservation of energy
$$
\d_t \left(\left(\frac{v^2}{2}+\e+\phi\right)\,G\right)\,+\,\Div_{\bfZ} \left(\left(\frac{v^2}{2}+\e+\phi\right)\,\cV_0\, G\right)
\,=\,0;
$$
\item the conservation of the classical adiabatic invariant $\mu_\perp=\e/B$
$$
\d_t \left(\frac{\e}{B}\,G\right)\,+\,\Div_{\bfZ} \left(\frac{\e}{B}\,\cV_0\,G\right)
\,=\, 0.
$$
\end{itemize} 
\epr
\begin{proof}
 For the asymptotic model \eqref{eq:1st}, the balance law for the kinetic energy is 
$$
\d_t \left(\left(\frac{v^2}{2}+\e\right)\,G\right)\,+\,\Div_{\bfZ} \left(\left(\frac{v^2}{2}+\e\right)\,\cV_0\, G\right)
\,=\, G\,v\,\Epar,
$$
which is a conservation law only if $\Epar\equiv0$. Then if $\bE$ derives from a potential, $\bE=-\nabla\phi$, the corresponding balance law for the total energy of the asymptotic model is
$$
\d_t \left(\left(\frac{v^2}{2}+\e+\phi\right)\,G\right)\,+\,\Div_{\bfZ} \left(\left(\frac{v^2}{2}+\e+\phi\right) \,\cV_0\, G\right)
\,=\, G\,\d_t\phi,
$$
which reduces to the claimed conservation law when $\d_t\phi\equiv0$. 

Note moreover that from
$$
\frac{\Div_\bfx\bB(t,\bfx)}{(B(t,\bfx))^2}\,=\,
\frac{\Div_\bfx\eDpar(t,\bfx)}{B(t,\bfx)}
-\eDpar(t,\bfx)\cdot\nabla_\bfx\left(\frac1B\right)(t,\bfx)
$$
follows for the asymptotic model the balance law 
$$
\d_t \left(\frac{\e}{B}\,G\right)\,+\,\Div_{\bfZ} \left(\frac{\e}{B}\,\cV_0\, G\right)
\,=\,-\,G\,\frac{\e}{B^2}\,\left(\d_tB+v\Div_\bfx\bB\right),
$$
which is indeed a conservation law when $\d_tB\equiv0$ and $\bB$ is divergence-free.
\end{proof}

\subsection{Second-order asymptotics}

Though already instructive, equation~\eqref{eq:1st} fails to capture parts of the slow dynamics that are too slow, for instance it does not describe the evolution of $\e/B$ (when $\Div_\bfx\bB\equiv0$). One way to correct this is to derive a higher-order description of the slow dynamics.

It is at this next order that are found macroscopic velocities, including those classically known as the $\bE\times \bB$ drift, the curvature drift, the grad-$B$ drift and the magnetic rotational drift, that with notation below read respectively  $\EcB(t,\bfy)$, $v^2\,
\curvB(t,\bfy)$,  $e\,\gradB(t,\bfy)$, and $\e\,\rotB(t,\bfy)$. Those have simple expressions in terms of vectors fields depending only on time $t$ and space $\bfy$ variables, and defined themselves as
\be
\label{drift-0}
\left\{\ba{l}
\ds\EcB \ds\,:=\,\frac{\bJ\,\bE}{B}
\,=\,\frac{\bE\wedge\bB}{B^2}\,,
\\[0.95em]
\ds \curvB\,:=\,
-\frac{\bJ}{B}(\dD_{\bfx}\eDpar\,\eDpar)
\,=\,-\frac{1}{B^2}\,\left(\dD_{\bfx}\eDpar\,\eDpar\right)\wedge\bB\,,
\\[0.95em]
\ds\gradB \,:=\, 
\bJ\nabla_\bfx\left(\frac1B\right)
\,=\,-\frac{1}{B^3}
\nabla_\bfx B \wedge \bB,
\\[0.95em]
\ds\rotB\,:=\,\frac{1}{B}\langle\Curl_\bfx \eDpar,\eDpar\rangle\,\eDpar,
\ea\right.
\ee
where $\bJ$ is given in \eqref{def:J}. Since the direction of the magnetic field $\eDpar$ is allowed to depend on time, another drift is present, given by $v\, \dtB(t,\bfy)$ where
\be
\ds\dtB\,:=\,
-\frac{\bJ}{B}\,\d_t\eDpar
\,=\,-\frac{1}{B^2}\,\d_t\eDpar\,\wedge\,\bB\,.
\label{drift-1}
\ee
Since it appears repeatedly it is convenient to introduce a piece of notation for a special combination of $\dtB$ and $\curvB$,
\be
\Sig(t,\bfy,v) \,=\, \dtB(t,\bfy) \,+\,v\,\curvB(t,\bfy)\,.
\label{drift-2}
\ee

With the above definitions we may write the full drift vector field $\Ud(t,\bfZ)$ in the concise form
\begin{eqnarray}
\label{drift-3}
\Ud(t,\bZ) \,=\, \left(\EcB\,+\,v\,\Sig\right)(t,\bfy,v)\,+\, \e\,\left( \rotB+\gradB\right)(t,\bfy).
\end{eqnarray}
where $\bZ=(\bfy,v,\e)$ stands for our set of slow variables in the
asymptotic model. For the sake of comparison with the existing
literature we observe the equivalent reformulations that may be derived from $(\dD_\bfx\eDpar)^*\eDpar=0$, a consequence of $\eDpar$ being unitary valued,
$$
\left\{
\ba{l}
\ds\dD_{\bfx}\!\eDpar\,\eDpar\,=\,\Curl_\bfx \eDpar\wedge\eDpar\,,
\\[0.9em]
\ds\curvB\,=\,-\frac{1}{B}(\Curl_\bfx \eDpar\wedge\eDpar)\wedge\eDpar
\ea\right.
$$
and observing that
$$
\left(\Curl_\bfx \eDpar \wedge \eDpar\right)\wedge \eDpar \,=\, -\Curl_\bfx\eDpar
\,+\,   \langle \Curl_\bfx\eDpar,\eDpar\rangle\,\eDpar
$$
we get that
\be
\label{formule:utile}
\rotB+\curvB\,=\, \frac{\Curl_\bfx \eDpar}{B}\,=\,\frac{\Div_\bfx \bJ}{B}
\ee

The vector-field involved in our higher-order description of the complete slow dynamics is then given by
\begin{equation}
  \label{V:eps}
\cV^\eps \,=\, \cV_0 \,+\, \eps\, \cV_1, 
\end{equation}
where the first order contribution $\cV_1(t,\bfZ)$ is
\begin{align}
\label{V:1}
\cV_1(t,\bfZ)
\,=\, \bp\ds
\Ud(t,\bfZ) 
\\[1em]
\ds\left\langle\Sig(t,\bfy,v),\bE(t,\bfy)\right\rangle
+\e\,\Div_\bfx\Sig(t,\bfy,v)
\\[1em]
\ds - \e\,\left[\,\left\langle
\curvB(t,\bfy),\bE(t,\bfy)\right\rangle + \Div_\bfx\left(\EcB+v\Sig\right)(t,\bfy,v)\,\right]
\ep\,.
\end{align}

The foregoing vector-field describes the dynamics of variables that are $\eps$-corrections of $(\bfx,\vpar,\eperp)$ but that are slower than those. The corrected spatial position 
\be
\label{GC}
\GC^\eps (t,\bfx,\bfv) \ :=\ \bfx+\eps\frac{\bJ(t,\bfx)\,\bfv}{B(t,\bfx)}
\ =\ \bfx+\eps\frac{\bfv\wedge \bB(t,\bfx)}{(B(t,\bfx))^2}
\ee
is well-known as the guiding center position, whereas  the corrected
parallel velocity is given as 
\be
\label{vGC}
\vGC^\eps (t,\bfx,\bfv)\,:=\, \vpar 
\ds +\; 
\eps\,\left\langle\bvperp,\,\Sig(t,\bfx,\vpar)\right\rangle
+\; \frac{\eps}{2 B(t,\bfx)} \,\left\langle\bJ(t,\bfx)\,\bvperp,\,\Re(\dD_\bfx\eDpar(t,\bfx))\,\bvperp\right\rangle
\ee
and the corrected version of the part of the kinetic energy in the plane perpendicular to the magnetic field direction is
\be
\label{eGC}
\left\{
\ba{lll}
\eGC^\eps(t,\bfx,\bfv) \,:=\,\eperp  &\ds
-\;\eps\,\left\langle\bvperp,\,
\EcB(t,\bfx)+\vpar\,\Sig(t,\bfx,\vpar)\right\rangle
\\[1.1em]
\; &\ds-\; \frac{\eps\,\vpar}{2B(t,\bfx)} \,\left\langle\bJ(t,\bfx)\,\bvperp,\,\Re(\dD_\bfx\eDpar(t,\bfx))\,\bvperp\right\rangle,
\ea\right.
\ee
where $\eperp=\frac{1}{2}\|\bfv_\perp\|^2$, whereas  $\Re$  denotes the symmetric part 
\be
\label{ReA}
\Re(\bA)\,=\,\frac{1}{2}\,(\bA+\bA^*),
\ee
with $\bA^*$ denoting the adjoint of $\bA$.
Therefore, our global sets of slower components are derived at time $t$ from
$(\bfx,\bfv)$ by 
$$
\bfZ_{\rm gc}^\eps(t,\bfx,\bfv) \,=\,(\GC^\eps, \vGC^\eps,
\eGC^\eps)(t,\bfx,\bfv).
$$
We can now  state our main theorem.

\bt
\label{th:2}
Let $\bE\in W^{2,\infty}$ and $\bB$ be such that $1/B\in W^{2,\infty}$ and $\eDpar\in W^{3,\infty}$. There exists a constant $C$ depending polynomially on $\|\bE\|_{W^{2,\infty}}$,  $\|B^{-1}\|_{W^{2,\infty}}$ and $\|\eDpar\|_{W^{3,\infty}}$ such that if $f^\eps$ solves \eqref{eq:vlasov} with initial data a nonnegative density $f_0$, then $F^\eps$ defined by
$$
F^\eps(t,\cdot)\,=\, \bfZ_{\rm gc}^\eps(t,\cdot)_* \,(f^\eps(t,\cdot))
$$
satisfies for a.e. $t\geq0$
$$
\|F^\eps(t,\cdot)-G^\eps(t,\cdot)\|_{\dot{W}^{-1,1}}\leq C\,\eps^2\,e^{C\,t^4\,(1+\eps\,t)}\,\int_{\R^6} e^{C\,t\,\|\bfv\|^3\,(1+\eps\,\|\bfv\|)}\,(1+\|\bfv\|^4)\,f_0(\bfx,\bfv)\,\dD\bfx\,\dD\bfv,
$$
where $G^\eps$ solves 
\be
\label{eq:2nd}
\d_t G^\eps\,+\,\Div_{\bfZ} \left(\cV^\eps\,G^\eps\right)
\,=\,0,
\ee
with $\cV^\eps$ given in \eqref{V:eps} and the initial data $G_{0}^\eps$ is
$$
G_{0}^\eps\,=\, \bZ_{\rm gc}^\eps(0,\cdot)_* \,(f_0)\,.
$$
\et
The proof of this asymptotic result is given in Section~\ref{sec:proofTh2}. A nonlinear counterpart for a system of Vlasov-Poisson type is also stated and proved in Section~\ref{s:2nd-nl} but restricts to well-prepared data in a sense detailed there.

A few comments on the structure of the asymptotic model are now in order. To begin with we observe that $\cV^\eps(t,\cdot)$ is still divergence-free. This follows from Lemma~\ref{lem:00} below and the fact that $\Div_\bfx\Div_\bfx(\bJ/B)=0$ by the skew-symmetry of values of
$\bJ$. 

\bl
\label{lem:00}
Consider $\rotB$, $\curvB$ and $\gradB$ defined in \eqref{drift-0}. Then we have
$$
\Div_\bfx\left( \frac{\bJ}{B}\right) \,=\, \rotB+\curvB+\gradB,
$$
where $\bJ$ is given in \eqref{def:J} and $B=\|\bB\|$.
\el
\begin{proof}
Straightforward by chain rule and \eqref{formule:utile} since $\rotB+\curvB=\Div_\bfx\left(\bJ\right)/B$ and $\gradB=\bJ\nabla_{\bfx}(1/B)$.
\end{proof}

Then we observe that we also have
 
\bpr
Assume that $\bE=-\nabla_\bfx\phi$, where $\phi$ does not depend on time. Then solutions to the asymptotic model \eqref{eq:2nd} satisfy the conservation of energy
$$
\d_t \left(\left(\frac{v^2}{2}+\e+\phi\right)\,G^\eps\right)\,+\,\Div_{\bfZ} \left(\left(\frac{v^2}{2}+\e+\phi\right)\,\cV^\eps\, G^\eps\right)
\,=\,0;
$$
\epr
\begin{proof}
If $\bE$ derives from a potential, $\bE=-\nabla\phi$, then one obtains the following balance law for the total energy of the second-order asymptotic model
$$
\d_t \left(\left(\frac{v^2}{2}+\e+\phi\right)G^\eps\right)\,+\,\Div_{\bfZ} \left(\left(\frac{v^2}{2}+\e+\phi\right)\,\cV^\eps \,G^\eps\right)
\,=\,G^\eps\,\d_t\phi
$$
by using Lemma~\ref{lem:00} and observing that
$$
-\left\langle
\Div_\bfx\left(\frac{\bJ}{B}\right),\bE\right\rangle 
+\Div_\bfx\EcB
=
\Tr\left(\frac{\bJ}{B}\dD_\bfx\bE\right)\,=\,0
$$
since $\bJ$ is skew-symmetric and $\dD_\bfx\bE$ is symmetric. From this stems the claimed conservation of energy when $\d_t\phi\equiv0$.
\end{proof}

\subsection{Long-time asymptotics in a toroidal axi-symmetric geometry}\label{s:poloidal-results}

Another way to unravel the dynamics of slower components is to derive asymptotics that hold on time scales of typical size $1/\eps$. Yet this seems doable only if the dynamical geometry of the first asymptotic model captured by Theorem~\ref{th:1} is sufficiently confining to ensure that the motion at speed of typical size $1$ is purely oscillatory and thus may be uncoupled from a dynamics evolving with macroscopic velocities of typical size $\eps$.

Our claim is that when such conditions are satisfied the proof of Theorem~\ref{th:2}, and more specifically the normal form on which it hinges (see System~\eqref{e:1new3}), contains sufficient ingredients to identify this long-time dynamics. To support this claim we illustrate it with a consideration of one of the simplest non trivial confining geometries.

We fix now a unitary vector $\eDz$ and for any $\bfx\in\R^3$ define the coordinate of $\bfx$ along $\eDz$ and its distance to the axis $\R\eDz$
$$
z(\bfx)\,=\,\langle\eDz,\bfx\rangle\,,\qquad\qquad
r(\bfx)\,=\,\|\eDz\wedge\bfx\|\,.
$$
We assume that for some $r_0>0$, where $r(\bfx)\geq r_0$, $\bB$ and $\bE$ are axi-symmetric, $\bB$ is stationary and toroidal and $\bE$ is orthogonal to $\bB$, that is,
\begin{align*}
\eDpar(\bfx)&=\frac{\eDz\wedge\bfx}{r(\bfx)}\,,&
B(t,\bfx)&=b(r(\bfx),z(\bfx))\,,&\Epar(t,\bfx)&=0\,,
\end{align*}
and
$$
\Eperp(t,\bfx)\,=\,\Er(t,r(\bfx),z(\bfx))\,\eDr(\bfx)+\Ez(t,r(\bfx),z(\bfx))\,\eDz
$$
for some smooth $b$, $\Er$, $\Ez$, with $b$ non vanishing and
$$
\eDr(\bfx)=\eDpar(\bfx)\wedge\eDz\,.
$$ 

Under the foregoing geometric assumptions, we have both $\Epar\equiv0$ and $\Div(\eDpar)\equiv0$ so that the only motion at speed of typical size $1$ is the rotation of $\bfx$ around $\eDz$ at angular velocity $\vpar$. Since by axi-symmetric assumption the corresponding angle is easily factored out one may expect to capture a slow dynamics at typical speed $\eps$ in variables $(r,z,\vpar,\eperp)$. This is the content of the next theorem. See Remark~\ref{rk:axi} for some hints on the relaxation of the assumptions made here for simplicity.
\begin{figure}[ht!]
\begin{center}
\begin{tikzpicture}[scale=1.75]
\begin{axis}[
axis equal image,
axis lines=none,
xmax=26,
ymax=23,
zmax=5,
ticks=none,
clip bounding box=upper bound,
colormap/blackwhite]
\addplot3[shader=interp,domain=0:360,y domain=0:250, samples=20,surf,z buffer=sort]
({(11 + 5.25 * cos(x)) * cos(y)}, 
{( 11 + 5.25 * cos(x)) * sin(y)},
{          5.25 * sin(x)});
\draw [-latex,thick,blue] (axis cs: -2,0,0) -- node[yshift=0.5em,xshift=3.em,scale=0.65]{$r(\bfx)$} (axis cs: 24,0,0);
\draw [-latex,thick,blue]  (axis cs: 0,0,-2) --node[pos=0.95,xshift=1.2em,scale=0.65]{$z(\bfx)$}(axis cs: 0,0,16);

\draw [-latex,black,thick,domain=40:80] plot ({90+180*cos(\x)},{300+180*sin(\x)});
\node[scale=0.65] at (260,420) {$\eDpar(\bfx)$};
\end{axis}
\end{tikzpicture} 
\end{center}
\label{fig:1}
\caption{Representation  of the torus local frame $({\bf e}_r(\bfx),\, {\bf
    e}_\parallel(\bfx), \,{\bf e}_z(\bfx))$ where the magnetic
field is along the unit vector field  ${\bf
    e}_\parallel$  whereas the electric field $\bE$ is
orthogonal to the magnetic field $\bB$.}
\end{figure}
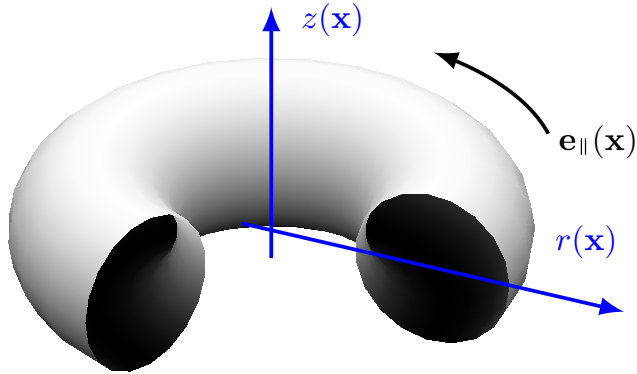

The involved asymptotic vector field is $\eps\cW_1$ with $\cW_1$ defined as
$$
\cW_1(t,\bfZ)
\,=\,\bp\ds
-\frac{\Ez}{b}(t,r,z)-\e\,\d_z\left(\frac1b\right)(r,z)\\[1em]\ds
\frac{\Er}{b}(t,r,z)+\frac{v^2}{r\,b(r,z)}
+\e\,\d_r\left(\frac1b\right)(r,z)
\\[1em]
\ds\frac{v}{r}\,\left(\frac{\Ez}{b}(t,r,z)
+\e\,\d_z\left(\frac1b\right)(r,z)\right)\\[1em]\ds
\e\,\left[\d_r\left(\frac{\Ez}{b}\right)(t,r,z)-\d_z\left(\frac{\Er}{b}\right)(t,r,z)
-\frac{v^2}{r}\,\d_z\left(\frac1b\right)(r,z)\right]
\ep,
$$
where the new slow variable is $\bfZ=(r,z,v,\e)$. 

\br
On the two first components of $\cW_1$ one readily identifies that in the present geometry along slower variables only survive as spatial drifts, the $\bE\times \bB$ and grad-$B$ drifts. This is due to the fact that here $\Curl_\bfx \eDpar$ vanishes identically in the zone of interest.
\er

\bt
\label{th:3}
Let $\bB$ be a stationary, axi-symmetric and toroidal magnetic field and $\bE$ be
an axi-symmetric electric field orthogonal to $\bB$, with $(\Er,\Ez,1/b)\in
W^{2,\infty}$ in the region where $r(\bfx)\geq r_0$ for some
$r_0$. For any $r_1>r_0$, there exist positive constants $\eps_0$ $\tau_0$ and $C_0$, $(1/\eps_0,1/\tau_0,C_0)$ depending polynomially on $1/r_0$, $1/(r_1-r_0)$ and $\|(\Er,\Ez,1/b)\|_{W^{2,\infty}([r_0,\infty[\times\R)}$, such that the following holds with
$$
\eps_{max}(R_0):=\frac{\eps_0}{1+R_0}\qquad\textrm{and}\qquad
T_{max}(R_0):=\frac{\tau_0}{1+R_0^2}\,.
$$
Consider $f^\eps$ a solution to \eqref{eq:vlasov} with initial datum a nonnegative density $f_0$ supported where 
$$
r(\bfx)\geq r_1\qquad\textrm{and}\qquad \|\bfv\|\leq R_0 
$$
for some $R_0>0$ and define $F^\eps$ as
$$
F^\eps(t,r,z,v,e)=\int_0^{2\pi}\!\!\int_0^{2\pi}
\,f^\eps(t,r\,\eDr^\theta+z\eDz,
v\,\eDpar(r\,\eDr^\theta+z\eDz)
\,+\,\sqrt{2\,\e}\ \eDperp^{\theta,\,\varphi})\ r\,\dD\varphi\,\dD\theta
$$
with  
$$
\eDr^\theta\,=\,\cos(\theta)\eDx+\sin(\theta)\eDy\,,\qquad\qquad
\eDperp^{\theta,\,\varphi}\,=\,\cos(\varphi)\eDr^\theta+\sin(\varphi)\eDz\,,
$$
where $(\eDx,\eDy,\eDz)$ is a fixed\footnote{$(\eDx,\eDy)$ are somehow arbitrary and the particular choice made here does not change $F^\eps$. In contrast we recall that $\eDz$ plays a special role as it directs the axis of symmetry.} orthonormal basis.  Then provided that
$$
0<\eps\leq \eps_{max}(R_0)\,,
$$
we have for a.e. $0\leq t \leq T_{max}(R_0)/\eps$
$$
\|F^\eps(t,\cdot)-G^\eps(t,\cdot)\|_{\dot{W}^{-1,1}}\leq C_0\,\eps\,
\int_{\R^3\times\R^3} e^{C\,\eps\,t\,\|\bfv\|^4}\,(1+\|\bfv\|^3)\,f_0(\bfx,\bfv)\,\dD\bfx\,\dD\bfv,
$$
where $G^\eps$ solves 
\be\label{eq:axi}
\d_t G^\eps\,+\,\eps\Div_{\bfZ} \left(\cW_1\, G^\eps\right)
\,=\,0,
\ee
with initial datum $G_0$ given by
$$
G_0(\bfZ)\,=\,
\int_0^{2\pi}\!\!\int_0^{2\pi}
\,f_0(\,r\,\eDr^\theta+z\eDz,
v\,\eDpar(r\,\eDr^\theta+z\eDz)
\,+\,\sqrt{2\,\e}\ \eDperp^{\theta,\,\varphi})\ r\,\dD \varphi\,\dD\theta\,.
$$
\et

Again note that averaging formulas coincide with push-forwards by the
slow map $(\bfx,\bfv)\mapsto(r,z,\vpar,\eperp)(t,\bfx,\bfv)$. Furthermore, we observe that  $r\,\cW_1$ is divergence-free and provide the following analogous to Proposition~\ref{prop:cons1}.
 
\bpr
\label{prop:cons3}
Suppose that $\bE$ derives from a stationary axi-symmetric potential 
$$
\phi(\bfx)=\phi(r(\bfx),z(\bfx))\,,\qquad\qquad
\bE=-\nabla_\bfx\phi.
$$
Then solutions to the asymptotic model \eqref{eq:axi} satisfy
\begin{itemize}
\item the conservation of energy;
\item the conservation of the classical adiabatic invariant.
\end{itemize}
\epr
\begin{proof}
When $\bE$ derives from an axi-symmetric potential as above, the corresponding balance law for the total energy of the asymptotic model \eqref{eq:axi} is
$$
\d_t \left(\left(\frac{v^2}{2}+\e+\phi\right)G^\eps\right)\,+\,\eps\,\Div_{\bfZ} \left(\left(\frac{v^2}{2}+\e+\phi\right)\,\cW_1\,G^\eps\right)
\,=\,G^\eps\,\d_t\phi
$$
which is a conservation law if furthermore $\d_t\phi\equiv0$. Moreover observe that $\bB$ is divergence-free in the present configuration and that the asymptotic model \eqref{eq:axi} comes with the balance law
$$
\d_t \left(\frac{\e}{b}\,G^\eps\right)\,+\,\eps\,\Div_{\bfZ} \left(\frac{\e}{b}\,\cW_1\,G^\eps\right)
\,=\,\eps\,G^\eps\,\frac{\e}{b^2}\,\left(\d_r\Ez-\d_z\Er\right),
$$
which is a conservation law if $\bE$ is curl-free, hence in particular if $\bE$ derives from a potential.
\end{proof}

\subsection{Further comments and numerical illustrations}

Though we have chosen not add this level of (mostly notational) complexity, the introduction of parameter dependencies in fields $\bB$ and $\bE$ would be immaterial to our analysis provided they satisfy upper bounds on $1/B$, $\bE$ and the needed number of derivatives of $(\bB,\bE)$. See the related Remarks~\ref{rk:toy-scaled-time} and~\ref{rk:toy-scaled-time}. In particular in this context one may expand and simplify further our asymptotic systems if one assumes an expansion of $\bB$ and $\bE$ with respect to $\epsilon$ or an $\epsilon$-ordering of gradients, or likewise one may perform a second asymptotic expansion with respect to another small parameter...

We believe that our leading-order slow variables $(\bfx,\vpar,\eperp)$ are both simple and natural. Yet many other choices have been used in the literature, and for comparison we provide in Appendix~\ref{s:mu} versions of our main results with another commonly-used choice, $(\bfx,\vpar,\eperp/B)$. Once a leading-order choice has been made, the higher-order corrections added to it to reach varying order of slowness are uniquely determined provided that a normalization is chosen. All through our analysis our implicit choice is to enforce that corrections have no slow component in the sense that they have zero mean with respect to the fast angle. See the related Remarks~\ref{trace-free} and~\ref{trace-free-bis}.

Though we have chosen to focus on the description of the slow
dynamics, the method would equally well provide a detailed description
of the oscillations as slaved to the evolution of the slow
variables. We stress that in most of methods relying on variational
principles one needs to provide both descriptions jointly even though
the oscillating part is subordinated to the slow part, as those
methods proceed by performing full changes of variables in the
original phase-space preserving the geometric structure under
consideration. Note that in principle to be fully justified from a
mathematical point of view this requires a careful tracking of how
small $\eps$ must be to guarantee that performed transformations are
indeed changes of variables. Here instead, with the exception of
results from Section~\ref{s:poloidal-results} where $\eps$ is
constrained to ensure sufficient confinement on large times  and of Section~\ref{s:nl} where nonlinear effects are analyzed, our
results are free of smallness constraint on the gradient of force
fields $(\bE,\bB)$ and on the initial data.

In Section~\ref{s:nl}, we have chosen to exemplify how our analysis extends to nonlinear models on a smoothed Vlasov-Poisson system. On this system, the nonlinear first-order analysis appears almost as a corollary of Theorem~\ref{th:1}. However we believe that the significantly more involved second-order analysis --- yielding a nonlinear counterpart to Theorem~\ref{th:2} --- illustrates well typical difficulties of related nonlinear analyses, and provides robust versatile solutions to those. In particular, it shows how the identification of slow variables yields suitable changes of variables effectively reducing the derivation of nontrivial uniform estimates --- of paramount importance in the nonlinear asymptotic analysis --- to standard arguments of PDE analysis. Moreover, since there we do perform a change of variables to get the latter uniform estimates, we also show how instead of modifying $\eperp$, a path that in a change of variables would bring polar coordinates singularities, one may adapt our slow-correction strategy and directly modify $\bvperp$. For more singular models, such arguments providing nontrivial uniform bounds may already be required to prove nonlinear counterparts to Theorem~\ref{th:1}.

A feature specific to nonlinear models is that oscillations of particle densities force oscillations in fields themselves that could be resonant in return with the particle densities. When this phenomenon do occur at a relevant order, it impacts possible asymptotic models. For this reason, to provide a result as close as possible to the linear analysis of Theorem~\ref{th:2}, we have restricted our nonlinear second-order analysis to initial data satisfying a form of preparation and a significant part of our analysis aims at proving that the well-prepared character is propagated by the stiff intricate dynamics. This leaves fully open the investigation of nonlinear resonances, a question that we regard as one of the main open questions in the field. However, in contrast with the upshots of the present contribution, we expect that the outcome of such an analysis would be extremely sensitive to fine details of the model under consideration.

\bigskip

To conclude the presentation of our results, we provide the reader with some numerical simulations illustrating and hopefully making more concrete respective error bounds. Since it is simpler to visualize we restrict numerical experiments to single-particle simulations.

In the present numerical experiments, we choose the electric field equal to zero and the
magnetic field $\bB$ as
$$
\bB(t,\bfx) = -\frac{20}{\|\bfx\|^5} \left(
\ba{c}
 3 \,x\,z \\ 
3 \,y\, z \\ 
2 \;z^2\,-\, x^2 \,-\, y^2
\ea
\right),
$$ 
where $\bfx=(x,y,z)$. The initial data is
$$
\bfx(0)=(10,0,0),\qquad \bfv(0)=\left(0,6\cos\left(\frac{\pi}{3}\right), 6\sin\left(\frac{\pi}{3}\right)\right),
$$
We approximate the solution of the  initial system \eqref{eq:xv} using a fourth order
Runge-Kutta scheme with a time step sufficiently small to resolve oscillations and compare it with the numerical solution obtained with the first order approximation corresponding to the characteristic curves of \eqref{V:0}-\eqref{eq:1st} (given in Proposition \ref{1st-xve})  and with  the second order approximation corresponding to the characteristic curves of \eqref{V:0}, \eqref{V:1} and \eqref{eq:2nd} (given in Proposition \ref{prop:4.9}).

\begin{figure}[h!]
\begin{subfigure}[b]{0.49\textwidth}
\centering
 \includegraphics[width=\textwidth]{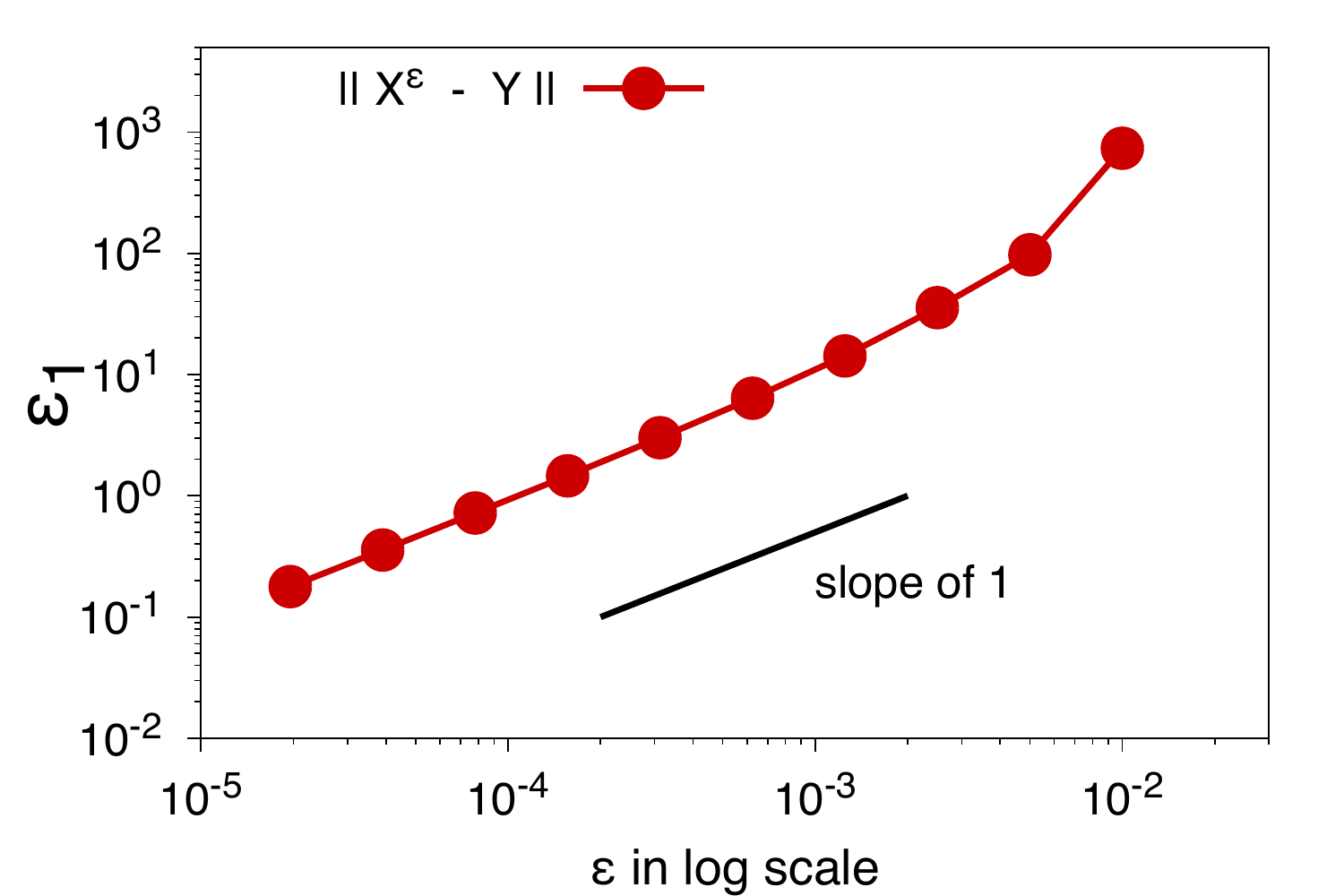}
\subcaption{}
\end{subfigure}
\begin{subfigure}[b]{0.49\textwidth}
\centering
 \includegraphics[width=\textwidth]{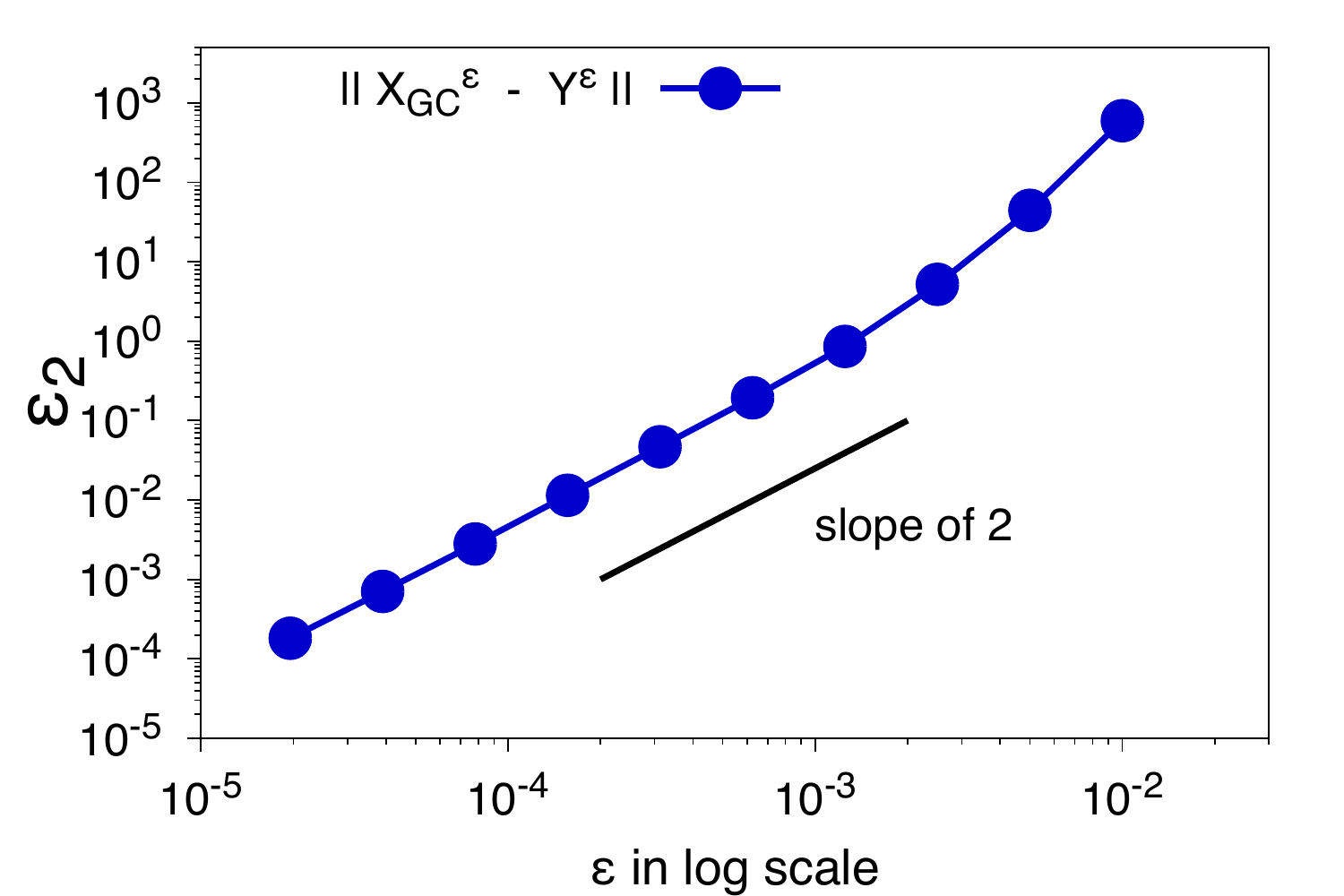}
\subcaption{}
\end{subfigure}
\begin{subfigure}[b]{0.49\textwidth}
\centering
 \includegraphics[width=\textwidth]{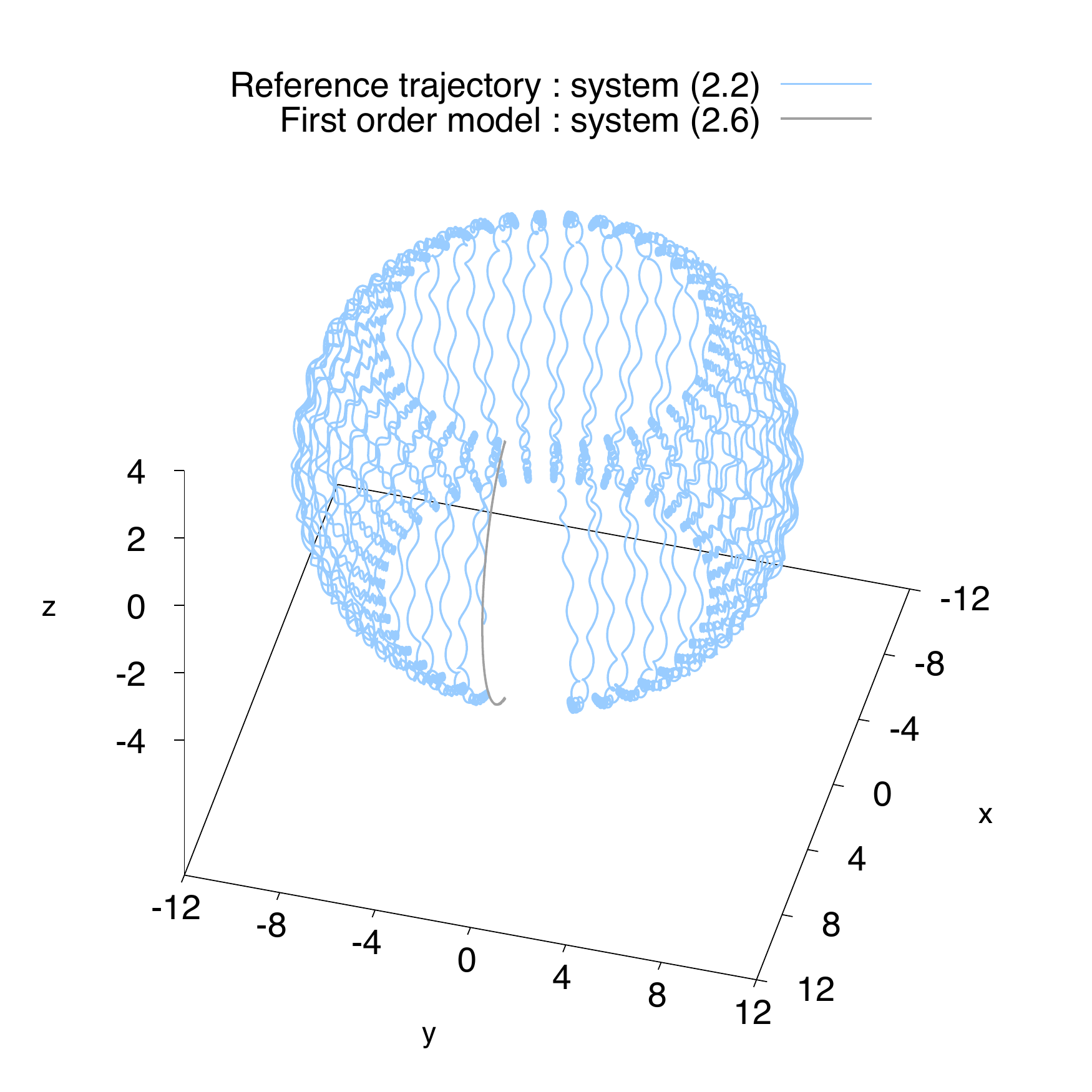}
\subcaption{}
\end{subfigure}
\begin{subfigure}[b]{0.49\textwidth}
\centering
 \includegraphics[width=\textwidth]{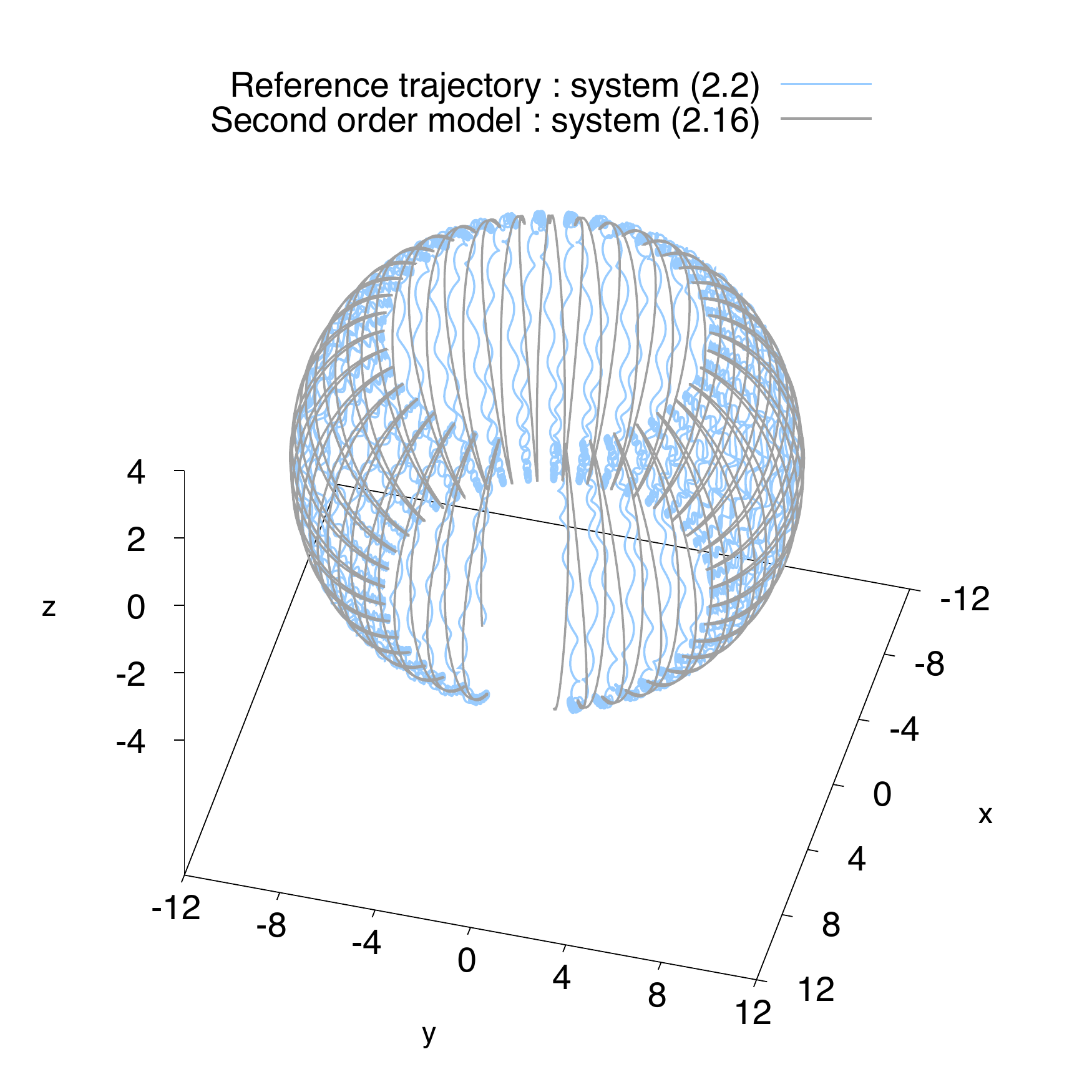}
\subcaption{}
\end{subfigure}
\caption{Numerical simulation: (A) error
  between the first order model and (\ref{eq:xv}) (B) error
  between the second order model and (\ref{eq:xv}) (C) particle
  trajectory obtained with the first order model and (\ref{eq:xv}) (D) particle
  trajectory obtained with the second order model and (\ref{eq:xv}).}
\label{fig:01}
\end{figure}

On the one hand, Figure~\ref{fig:01}-(A) represents the error 
$$
\cE_1(\eps) = \int_0^{10} \| \bfx^\eps(t) - \bfy(t) \| \dD t,
$$
where $\bfx^\eps$ is the spatial component of the solution $(\bfx^\eps,\bfv^\eps)$ to \eqref{eq:xv} and $\bfy$ is the spatial component of $(\bfy,v,\e)$ satisfying  the differential system with initial data as in Proposition~\ref{1st-xve} (and corresponding to~\eqref{eq:1st}), whereas Figure~\ref{fig:01}-(B) represents the error 
$$
\cE_2(\eps) = \int_0^{10}\| \GC^\eps(t) - \bfy^\eps(t) \| dt,
$$
where $\GC^\eps$ is obtained through \eqref{GC} from $(\bfx^\eps,\bfv^\eps)$
solving \eqref{eq:xv} and $\bfy^\eps$ is the spatial component of $(\bfy^\eps,v^\eps,\e^\eps)$ satisfying the Cauchy problem as in Proposition~\ref{prop:4.9} (and corresponding to~\eqref{eq:2nd}). These numerical results illustrate the order of accuracy stated in Theorems~\ref{th:1} (first order) and~\ref{th:2} (second order). On the other hand, we have also claimed that to capture long-times dynamics it is also crucial to include second-order terms in the asymptotic models. Theorem~\ref{th:3} provides some quantitative support to the claim, in a specific geometry. We now provide in a different configuration another, qualitative, illustration of the claim, by plotting in Figures~\ref{fig:01}-(C) and~\ref{fig:01}-(D) examples of spatial parts of particle trajectories obtained with original, first and second-order models. Here we take $\eps=10^{-3}$ and solve on $[0,250]$. Roughly speaking, the solution to the original problem exhibits a superposition of three kinds of spatial motions, namely, with decreasing velocity, the cyclotron oscillation about magnetic field lines, an oscillation along magnetic field lines, and a slower drift, responsible of a horizontally circular displacement. By construction, both asymptotic descriptions remove the cyclotron motion. However, whereas the second-order asymptotic model seems able to reproduce the slow part of the complicate multi-scale behavior, the first-order one only describes oscillations along the magnetic field lines. Indeed, since the first order model does not include classical drifts $\gradB$, $\curvB$ and $\rotB$, it is not adequate to follow accurately the correct trajectory on times sufficiently long to feel the effects of those. 

\section{Two-dimensional homogeneous case}\label{s:homogeneous}

As a warm up we begin our analysis by revisiting the two-dimensional homogeneous case. The goal is to expound the tenets of the method without being slowed down by computational complexity. For the sake of exposition, for this simple system we prove first results that are even weaker than what the method may prove but that correspond to the best that is expected from the general $3$-D system without assuming special symmetry.
\begin{figure}[ht!]
\begin{center}
\begin{tikzpicture}[plane/.style={trapezium,draw,fill=black!30,trapezium left angle=60,trapezium right angle=120,minimum height=2.5cm},scale=1.75]
\node (p)[plane] at (0,0){.};
\draw (p.center) edge ++(0,0.75cm) edge[densely dashed] (p.south) (p.south) edge ++(0,-0.25cm);
\draw[-latex,blue,thick] (0,0) --  (0.75,0) node[anchor=north east]{$x$};
\draw[-latex,blue,thick] (0,0) -- (0.27,0.5) node[anchor=north west]{$y$};
\draw[-latex,black,thick] (-1,-0.5) -- (-1,0.2) node[anchor=north west]{$\bB$};
\end{tikzpicture}
\end{center}
\label{fig:2}
\caption{Representation  of perpendicular plane to  the magnetic
field $\bB= B(\bfx)\,{\bf e}_z$.}
\end{figure}
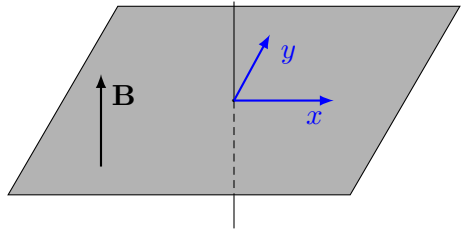

Since here the parallel direction is fixed and we follow only
perpendicular motions we drop temporarily ${}_\perp$ and ${}_\mypar$
indices.  Thus, we consider for any $(\bfx,\bfv)\in\R^2\times\R^2$ and
$t\geq 0$,
\be
\label{eq:vla2d}
\ds\d_t f^\eps\,+\,\Div_\bfx \left(f\,\bfv\right)
\,+\,\Div_{\bfv}\left(f\,\left(\frac1\eps\,B\,\bJ\bfv +\bE\right)\right)
\,=\,0\,.
\ee
Characteristics of the underlying PDE are obtained by solving
\be\label{eq:xv-hom}
\left\{
\ba{l}
\ds\frac{\dD\bfx^\eps}{\dD t}\,=\;\bfv^\eps,
\\[0.75em]
\ds\frac{\dD\bfv^\eps}{\dD t}\,=\,\frac{1}{\eps}\,B\,\bJ\bfv^\eps\,+\,\bE(t,\bfx^\eps),
\ea
\right.
\ee
with $B>0$ and
$$
\bJ\left(a_1,a_2\right)=(a_2,-a_1)\,,
$$
important properties of $\bJ$ being
$$
\bJ^2=-\Id\,,
\qquad\qquad
\bJ^*=-\bJ\,.
$$
  
For the sake of readability, from now on, when no confusion is possible, that is, when no asymptotic comparison is under consideration, we shall drop ${}^\eps$ exponents on solutions.

We shall perform a series of transformations so as to extract from system~\eqref{eq:xv-hom} a normal form where some slow variables satisfy a system of ODEs uncoupled from fast scales up to error terms. It is worth pointing out that under stringent assumptions on fields one may expect to perform at once an infinite number of transformations and uncouple at infinite order slow variables from fast variables. We shall not pursue this line of investigation here but as a consequence one should keep in mind that variables that we designate as slow are slow only up to a certain order and that depending on the objective at hand the level of slowness required may vary. As an example, anticipating a bit the analysis below, note that depending on the aimed conclusion one may be allowed to work directly with the spatial position $\bfx$ or need to manipulate the gyrocenter $\bfx+\eps B^{-1}\bJ\bfv$, or even be compelled to use a version of those corrected by higher-order powers of $\eps$.

\subsection{Uniform bounds}

Both to enforce that terms expected to be irrelevant are indeed
irrelevant and to ensure that solutions persist on a sufficiently long
time interval, uniform bounds on the solution are needed. Let us
obtain them by introducing a kinetic energy variable $\e(\bfv) =\tfrac12\|\bfv\|^2$ and noting that system~\eqref{eq:xv-hom} yields
$$
\left\{\ba{l}
\ds\frac{\dD\bfx}{\dD t}\,=\,\bfv,
\\[0.9em]
\ds\frac{\dD \e}{\dD t}\,=\,
\langle\bfv,\bE(t,\bfx)\rangle,
\\[0.9em]
\ds\frac{\dD\bfv}{\dD t}\,=\,\frac{B}{\eps}\bJ\bfv\,+\,\bE(t,\bfx)\,.
\ea\right.
$$

\br\label{rddnt-hom}
We warn the reader that though we write the latter system as if $\e$ and $\bfv$ were independent variables this is mostly an algebraic trick here. In particular one should keep in mind that the system does contain some redundancy that is kept for the sake of simplicity of algebraic manipulations. In contrast an \emph{augmented} formulation was in turn crucially used in \cite{FR2} jointly with suitable numerical schemes so as to allow the discretization to disconnect the weak convergence of $\bfv$ to zero from the strong convergence of $\e$ to a non trivial limit.
\er

Of course here one could obtain from a Lipschitz assumption on $\bE$ global-in-time existence and some bounds growing exponentially in time from a standard Gr\"onwall lemma. Yet for expository reasons we show how to perform simple better estimates. Note however that, as we derive below in the long-time analysis, those are still deceptively pessimistic.

\bl
\label{toy-bnd0}
Solutions to \eqref{eq:xv-hom} starting from $(\bfx_0,\bfv_0)$ are
defined globally in time and satisfy for any $t\geq0$
$$
\left\{
\ba{l}
\|\bfx(t)\|\,\leq\, \|\bfx_0\|
+\,t\,\|\bfv_0\|+\,t^2\,\|\bE\|_{L^\infty},
\\[0.9em]
\|\bfv(t)\|\,\leq\,
\|\bfv_0\|+2\,t\,\|\bE\|_{L^\infty}\,.
\ea\right.
$$
\el
\begin{proof}
From the equation on $\e$ stems, for any $t\geq0$, as long as the solution exists
$$
\max_{s\in[0,t]}\|\bfv(s)\|^2\,\leq\, 
\|v_0\|^2+2\max_{s\in[0,t]}\|\bfv(s)\|\,t\,\|\bE\|_{L^\infty},
$$
hence by solving the second-order inequality, for any $t\geq0$, as long as the solution exists
$$
\max_{s\in[0,t]}\|\bfv(s)\|\leq 
\sqrt{\|\bfv_0\|^2+t^2\|\bE\|_{L^\infty}^2}+\,t\,\|\bE\|_{L^\infty}\,.
$$
This yields the estimate on $\bfv$. In turn it implies the estimate on $\bfx$ by a mere integration, and jointly they prove global well-posedness by ruling out finite-time blow-up.
\end{proof}

\subsection{Elimination of linear terms}

We begin the uncoupling process. The thrust of the method is that the
equation that forces $\bfv$ --- or more exactly its argument --- to
evolve on fast scales also provides a way to eliminate at leading
order $\bfv$ --- or more exactly dependences on its argument --- in
slow equations. This general philosophy, that may be turned into
rigorous arguments, explain why slow evolutions may be uncoupled from
fast scales at any prescribed order. Explicitly elimination, at leading order, of linear terms in $\bfv$ is summarized as

\bl
\label{toy-1st}
Consider $\bL\in W^{1,\infty}\left(\R^+_t;\,\cL_1(\R^2,\R^p)\right)$, $p\in\N^*$ and $(\bfx,\bfv)$ a solution to \eqref{eq:xv-hom}. Then
for a.e. $t\geq 0$, we have
\be
\label{toy-res:0}
\bL(t)(\bfv(t))\,=\,-\eps\,\frac{\dD}{\dD t}\left(\bL(t)\left(\frac{\bJ\,\bfv}{B}\right)\right)
\,+\,\eps \,\bL'(t)\left(\frac{\bJ\bfv}{B}\right)
\,+\,\eps \,\bL(t)(\EcB(t,\bfx)),
\ee
with
$$
\EcB(t,\bfx)\,=\,\frac{\bJ \,\bE}{B}(t,\bfx)\,.
$$
\el

In the former we have used the following notational convention. For any $\alpha\in\N$, $\cL_\alpha(V,W)$ denotes the space of $\alpha$-linear operators from
$V^\alpha$ to $W$. In particular, $\cL_1(V,W)$ is the set of linear operators from $V$ to $W$. 

\begin{proof}
This follows directly from 
$$
\bfv\,=\,
-\frac{\eps}{B}\frac{\dD}{\dD t}(\bJ\bfv)\,+\,\frac{\eps}{B}\,\bJ\,\bE(t,\bfx)\,.
$$
\end{proof}

The latter term $\EcB$ identifies the classical $\bE\times\bB$ velocity drift from gyrokinetic theory. The foregoing lemma singles out the prominent role played by the $\EcB$ drift in two-dimensional gyrokinetics. 

\subsection{First elimination and partial asymptotics}

By using Lemma~\ref{toy-1st} first with $\bL(t)(\bfv)=\bfv$ then
with $\bL(t)(\bfv)=\langle\bE(t,\bfx(t)),\bfv\rangle$, one derives
that from system~\eqref{eq:xv-hom} follows for the guiding center
variable,
\be
\label{toy-x1}
\frac{\dD}{\dD t}\left[\bfx
+\eps\frac{\bJ\bfv}{B}\right]
\,=\, \eps\,\EcB(t,\bfx)
\ee
and for the corrected kinetic energy $\e(\bfv)=\frac{1}{2}\|\bfv\|^2$,
\be
\label{toy-e1}
\frac{\dD}{\dD t}\left[\e(\bfv)\,-\,\ds
\eps\,\left\langle\EcB(t,\bfx),\bfv\right\rangle\right]
\,=\,
-\eps\,\left\langle\d_t\EcB(t,\bfx)
+\dD_\bfx\EcB(t,\bfx)(\bfv),\bfv\right\rangle\,.
\ee

\br
In the present paper the stiffer part of the fast equation is always linear in $\bfv$. This leads to a quite simple elimination of terms linear in $\bfv$. In particular, since the slow equations on $(\bfx,\e)$ are linear in $\bfv$ the first elimination comes almost for free. However in general each simplification increases the level of nonlinearity in $\bfv$ of slow equations and subsequent simplifications get more and more algebraically cumbersome. 
\er

A specific feature of System~\eqref{eq:xv-hom} is that slow variables evolve with speeds of typical size $\cO(\eps)$ and not $\cO(1)$. Therefore on time intervals $[0,\tobs^\eps]$, one hopes to validate approximation of the slow part by the solution of an uncoupled system up to error terms of size $\cO(\eps^2\tobs^\eps)$ with $\tobs^\eps=\cO(\eps^{-1})$. We first prove this claim with $\tobs^\eps$ of size $1$ then refine the analysis to reach $\tobs^\eps$ of size $\eps^{-1}$. Note that when $\tobs^\eps$ is of size $\eps^{-1}$ we aim at an error of size $\cO(\eps)$ and thus we may use directly $(\bfx,w)$ as slow variables whereas when $\tobs^\eps$ is $\cO(1)$ we aim at precision $\cO(\eps^2)$ thus we should use 
$$
\left(\bfx+\eps \frac{\bJ\bfv}{B},\,\e(\bfv)\,+\,\ds
\eps\,\left\langle\bE(t,\bfx),\frac{\bJ\bfv}{B}\right\rangle\right),
$$
or a higher-order version of the latter.

Note that without further simplification the aforementioned asymptotics may not be readily derived since the equation on $e$ still contains $\bfv$-terms at leading order. However an aspect even more peculiar to System~\eqref{eq:xv-hom} is that at leading order the equation for $\bfx$ uncouples not only from the argument of $\bfv$ but also from $\e$. At this stage an asymptotic description of the slow part corresponding to $\bfx$ may be guessed without any further computation. 


\bpr
\label{toy-cheap}
Assume $\bE\in L^\infty\left(\R^+_t;\,W^{1,\infty}(\R^2)\right)$ and let
$(\bfx^\eps,\bfv^\eps)$ be the solution to~\eqref{eq:xv-hom} starting from
$(\bfx_0,\bfv_0)$. Then the guiding center variable (\ref{GC}) satisfies
for a.e. $t\geq 0$,
$$
\left\|\GC^\eps(t)-\bfy^\eps(t)\right\|
\,\leq\, 
\frac{\eps^2}{B^2}\|\dD_\bfx\bE\|_{L^\infty}\,e^{\frac{\eps\,t}{B}\|\dD_\bfx\bE\|_{L^\infty}}\,
(\,t\,\|\bfv_0\|\,+\,t^2\,\|\bE\|_{L^\infty}\,),
$$
where $\bfy^\eps$ solves
\be
\label{GC2D}
\left\{
\ba{l}
\ds\frac{\dD\bfy^\eps}{\dD t}\,=\,\eps\,\EcB(t,\bfy^\eps),
\\[0.95em]
\ds\bfy^\eps(0)=\GC(0).
\ea\right.
\ee
\epr

\begin{proof}
We consider $\GC^\eps$ given in (\ref{GC}), which satisfies  
$$
\frac{\dD\GC^\eps}{\dD t}\,=\,\eps\,\EcB(t,\GC^\eps)
\,+\,\eps\,\left[\EcB\left(t,\GC^\eps-\eps\,\frac{\bJ\bfv^\eps}{B}\right)-\EcB(t,\GC^\eps)\right]\,.
$$
This implies for a.e. $t\geq0$
$$
\|\GC^\eps(t)-\bfy^\eps(t)\|
\leq \frac{\eps}{B}\|\dD_\bfx\bE\|_{L^\infty}\int_0^t\|\GC^\eps(s)-\bfy^\eps(s)\|\dD s\\
+\frac{\eps^2}{B^2}\|\dD_\bfx\bE\|_{L^\infty}\int_0^t\|\bfv^\eps(s)\|\,\dD s\,.
$$
Thus by the Gr\"onwall lemma, for a.e. $t\geq0$,
$$
\|\GC^\eps(t)-\bfy^\eps(t)\|
\,\leq\, \frac{\eps^2}{B^2}\|\dD_\bfx\bE\|_{L^\infty}\,e^{\frac{\eps\,t}{B}\|\dD_\bfx\bE\|_{L^\infty}}\,\int_0^t\|\bfv^\eps(s)\|\,\dD s\,.
$$
Then the result follows from Lemma~\ref{toy-bnd0}.
\end{proof}

The foregoing bound is very simple but is not sharp with respect to $\eps$. Indeed the principal part of the error term of the equation is linear in $\bfv$ thus may also be eliminated. 

\br\label{antiderivative}
The special structure of equation~\eqref{toy-x1} is somewhat fortuitous. However the fact that the error introduced by replacing $\bfx^\eps$ with its $\eps$-correction $\GC^\eps$ may be cast away at leading order is not mere luck. It is due to our choice in \eqref{toy-res:0} of an antiderivative 
$$
\eps\,\bL(t)\left(\frac{\bJ\bfv}{B}\right),
$$
that at leading order contains no slow part. Henceforth in similar cases enforcing such properties will always streamline our particular choices.
\er

The announced further elimination yields the following refinement.

\bpr
\label{toy-x}
Assume that $\bE\in W^{2,\infty}$. There exists
a constant $C>0$, depending polynomially on $\|\bE\|_{W^{2,\infty}}$ and
$B^{-1}$, such that if $(\bfx^\eps,\bfv^\eps)$ is a solution to
\eqref{eq:xv-hom} starting from $(\bfx_0,\bfv_0)$, then it satisfies for a.e. $t\geq0$,
$$
\Big\|\GC^\eps(t)-\bfy^\eps(t)\Big\|
\,\leq\, C\,\eps^3\,e^{\frac{\eps\,t}{B}\|\dD_\bfx\bE\|_{L^\infty}}\,
(\,1\,+\,t\,(1+\|\bfv_0\|\,+\,t\,\|\bE\|_{L^\infty})\,)\,
(\|\bfv_0\|\,+\,t\,\|\bE\|_{L^\infty}),
$$
where $\GC^\eps$ is as in \eqref{GC} and $\bfy^\eps$ solves
(\ref{GC2D}).
\epr

\begin{proof}
The term to weed out is linear in $\bfv$ and by applying Lemma \ref{toy-1st}  with 
$$\bL(t)(\bfv)=-\dD_\bfx\EcB(t,\GC^\eps(t))\left(\frac{\bJ\bfv}{B}\right)$$ 
one obtains 
\begin{eqnarray*}
&&\frac{\dD}{\dD t}\left[\GC^\eps+\eps^3\dD_\bfx\EcB(t,\GC^\eps)\left(\frac{\bfv^\eps}{B^2}\right)
\right]
\\
&&=\eps\,\EcB(t,\GC^\eps)
-\eps^3\dD_\bfx\EcB(t,\GC^\eps)\left(\frac{\bJ\,\EcB(t,\bfx^\eps)}{B}\right)\\
&&\,+\,\eps^3\dD_\bfx\d_t\EcB(t,\GC^\eps)\left(\frac{\bfv^\eps}{B^2}\right)
\,+\,\eps^3\dD_\bfx^2\EcB(t,\GC^\eps)\left(\EcB(t,\bfx^\eps),\frac{\bfv^\eps}{B^2}\right)\\
&&\,+\,\eps\,\left[\EcB\left(t,\GC^\eps-\eps\,\frac{\bJ\bfv^\eps}{B}\right)-\EcB(t,\GC^\eps)
+\eps\dD_\bfx\EcB(t,\GC^\eps)\left(\frac{\bJ\bfv^\eps}{B}\right)\right]\,.
\end{eqnarray*}

Therefore, for a.e. $t\geq 0$, one has
\begin{align*}
\|\GC^\eps(t)-\bfy^\eps(t)\|
&\leq \frac{\eps}{B}\|\dD_\bfx\bE\|_{L^\infty}\int_0^t\|\GC^\eps(s)-\tby^\eps(s)\|\dD s
\,+\,\frac{\eps^3}{B^3}\,\|\dD_\bfx\bE\|_{L^\infty}\,(\|\bfv_0\|+\|\bfv^\eps(t)\|)\\
&+\,\frac{\eps^3}{B^3}
\,(\|\dD_\bfx\d_t\bE\|_{L^\infty}+B^{-1}\|\dD_\bfx^2\bE\|_{L^\infty}\|\bE\|_{L^\infty})
\,\int_0^t\|\bfv^\eps(s)\|\,\dD s\\
&+\,\frac{\eps^3}{2\,B^3}\|\dD_\bfx^2\bE\|_{L^\infty}
\,\int_0^t\|\bfv^\eps(s)\|^2\,\dD s\,.
\end{align*}
At this stage the result follows from Lemma~\ref{toy-bnd0} and the Gr\"onwall lemma.
\end{proof}

One may go on by correcting $\GC^\eps$ into a ``higher-order'' approximation
$$
\bfx_{\rm ho}^\eps
\,=\,\GC^\eps+\eps^3\dD_\bfx\EcB(t,\GC^\eps)\left(\frac{\bfv^\eps}{B^2}\right),
$$
then expanding from $\bfx_{\rm ho}^\eps$ and eliminating terms involving $\bfv^\eps$. But the expansion process would involve terms quadratic in $\bfv$ whose elimination brings a coupling with $\e^\eps$ as may be seen from Lemma~\ref{toy-2nd} below. Proposition~\ref{toy-x} is therefore expected to be optimal with respect to $\eps$-scaling on time intervals of length $\cO(1)$.

\subsection{Elimination of quadratic terms}

In order to obtain asymptotics for the full set of slow variables $(\bfx,\e)$ we study now the extraction of slow components from expressions that are quadratic in $\bfv$.

\bl
\label{toy-2nd}
Consider $\bA\in W^{1,\infty}\left(\R^+_t; \,\cL_2(\R^2,\R^p)\right)$, $p\in\N^*$ and $(\bfx,\bfv)$ a  solution to \eqref{eq:xv-hom}. Then, for
a.e. $t\geq$, we have
$$
\bA(t)(\bfv(t),\bfv(t))
\,=\,\e(\bfv)\,\Tr(\bA(t))-\eps\frac{\dD\chi_\bA}{\dD t}(t)\,+\,\eps\,\eta_\bA(t),
$$
where $\e(\bfv)=\|\bfv\|^2/2$, whereas $\Tr$ denotes the trace
operator, for the canonical basis $(\eDx,\eDy)$ of $\R^2$
$$
\Tr(\bA(t)) \,=\, \bA(t)(\eDx,\eDx) \,+\, \bA(t)(\eDy,\eDy)
$$
whereas $\chi_\bA$ is given by
$$
\chi_\bA \,=\,\frac12\,\Re(\bA)\left(\bfv,B^{-1}\bJ\,\bfv\right)
$$
and $\eta_\bA$ is
\begin{eqnarray*}
\eta_\bA& =&
\frac12\,\Re(\bA)(\bfv,\EcB(t,\bfx))
\,+\,\frac12\,\Re(\bA)'\left(\bfv,B^{-1}\bJ\,\bfv\right)
\\
\nonumber
&+&\frac12\,\Re(\bA)\left(\bE(t,\bfx),B^{-1}\bJ\,\bfv\right),
\end{eqnarray*}
with $\Re$ denoting the symmetric part defined in (\ref{ReA}).
\el

\br
\label{trace-free}
Consistently with Remark~\ref{antiderivative}, note that $\chi_\bA$ has itself no slow component at leading order since $\Re(\bA(t))\left(\cdot,\bJ(\cdot)\right)$ is trace-free. Indeed its trace is
\begin{eqnarray*}
\Re(\bA)\left(\bfa,\bJ\bfa\right)+\Re(\bA)\left(\bJ\bfa,\bJ\,\bJ\bfa\right)
&=&\Re(\bA)\left(\bfa,\bJ\bfa\right)-\Re(\bA)\left(\bJ\bfa,\bfa\right)
\\
&=&0,
\end{eqnarray*}
where $\bfa$ is any unitary vector. In the latter to express the trace we have used that $(\bfa,\bJ\bfa)$ form an orthonormal basis for any unitary $\bfa$.
\er

\begin{proof}
Note first that one may assume without loss of generality that $\bA$ is valued in symmetric bilinear forms. Thus we assume $\Re(\bA)=\bA$ for the sake of notational concision. By differentiation one derives
\begin{eqnarray*}
2\frac{\dD\chi_\bA}{\dD t}(t)&=&
\bA'(t)\left(\bfv(t),B^{-1}\bJ\bfv(t)\right)
\\
&+&\bA(t)\left[\frac{\dD\bfv}{\dD t}(t),B^{-1}\bJ\bfv(t)\right]
\,+\,\bA(t)\left[\bfv(t),B^{-1}\bJ\left(\frac{\dD\bfv}{\dD
                                 t}(t)\right)\right]
\\
&=&
\bA'(t)\left(\bfv(t),B^{-1}\bJ\bfv(t)\right)
\\
&+&\frac1\eps\bA(t)\left(\bJ\bfv(t),\bJ\bfv(t)\right)
+\bA(t)\left(\bE(t,\bfx(t)),B^{-1}\bJ\bfv(t)\right)
\\
&-&\frac1\eps\bA(t)\left(\bfv(t),\bfv(t)\right)
+\bA(t)\left(\bfv(t),\EcB(t,\bfx(t))\right)
\end{eqnarray*}
and the result follows by multiplying by $\eps/2$ then adding $\bA(t)(\bfv(t),\bfv(t))$ and using that
$$
\e\,\Tr(\bA)
\,=\,\frac12\,\left(\bA(\bfv,\bfv)
+\bA(\bJ\bfv,\bJ\bfv)\right)\,.
$$
\end{proof}

The last equality of the foregoing proof is essentially the definition of the trace operator. An elementary but fundamental point is that the right-hand side is invariant by rotation, thus the definition of $\Tr(\bA)$ does not depend on the vector $\bfv$ chosen to express it. 

\subsection{Second elimination and full asymptotics}\label{s:hom-full}
For the sake of concision and symmetry we introduce
\be
\label{toy-egc}
\eGC^\eps\,=\,
\e(\bfv^\eps)\,-\,\eps\,\left\langle\EcB(t,\bfx^\eps),\bfv^\eps\right\rangle\,,
\ee
which corresponds to the corrected kinetic energy, a two-dimensional version of \eqref{eGC}.

By applying Lemmas~\ref{toy-1st} with
$$
\bL(t)(\bfv)=
-\left\langle\d_t\EcB(t,\bfx(t)),\bfv\right\rangle
$$
and Lemma \ref{toy-2nd} with
$$
\bA(t)(\bfv,\bfu) \,=\,
-\left\langle\dD_\bfx\EcB(t,\bfx(t))(\bfv),\bfu\right\rangle,
$$
equation~\eqref{toy-e1} in system \eqref{toy-x1}-\eqref{toy-e1} may be turned into
\be
\label{toy-e2}
\frac{\dD}{\dD t}\left[\eGC^\eps+\frac{\eps^2}{B^2}\,\chi^\eps\right]
\,=\,-\eps\,\e(\bfv^\eps)\,\Div_\bfx(\EcB)(t,\bfx^\eps)
\,+\,\frac{\eps^2}{B^2}\,\eta^\eps,
\ee
where
$$
\chi^\eps(t,\bfx,\bfv)\,=\,\left\langle\d_t\bE(t,\bfx),\bfv\right\rangle
-\frac14\,\left[\left\langle\dD_\bfx\bE(t,\bfx)(\bfv),\bfv\right\rangle
\,-\,\left\langle\dD_\bfx\bE(t,\bfx)(\bJ\bfv),\bJ\bfv\right\rangle\right]
$$
and
\begin{eqnarray*}
\eta^\eps(t,\bfx,\bfv)
&=&-\,\left\langle\d_t\bE(t,\bfx),\bE(t,\bfx)\right\rangle
-\,\left\langle\d_t^2\bE(t,\bfx)
+\dD_\bfx\d_t\bE(t,\bfx)(\bfv),\bfv\right\rangle\\
&+&\frac14\,\left\langle\dD_\bfx\bE(t,\bfx)(\bJ\bE(t,\bfx)),\bJ\bfv\right\rangle
-\frac14\left\langle\dD_\bfx\bE(t,\bfx)(\bfv),\bE(t,\bfx)\right\rangle\\
&-&\frac14\left\langle\dD_\bfx\bE(t,\bfx)(\bE(t,\bfx)),\bfv\right\rangle
+\frac14\left\langle\dD_\bfx\bE(t,\bfx)(\bJ\bfv),\bJ\bE(t,\bfx)\right\rangle\\
&-&\frac14\left\langle\dD_\bfx\d_t\bE(t,\bfx)(\bfv)
+\dD_\bfx^2\bE(t,\bfx)(\bfv,\bfv),\bfv\right\rangle\\
&+&\frac14\left\langle\dD_\bfx\d_t\bE(t,\bfx)(\bJ\bfv)
+\dD_\bfx^2\bE(t,\bfx)(\bfv,\bJ\bfv),\bJ\bfv\right\rangle\,.
\end{eqnarray*}

Now we may complete Proposition~\ref{toy-x} to obtain leading-order asymptotics for $(\GC^\eps,\eGC^\eps)$.

\bpr
\label{toy-xe}
Assume  $\bE\in W^{2,\infty}$. There exists a constant $C>0$, depending
polynomially on $\|\bE\|_{W^{2,\infty}}$ and $B^{-1}$ such that the following holds. Let
$(\bfx^\eps,\bfv^\eps)$ be the solution to~\eqref{eq:xv-hom} starting from
$(\bfx_0,\bfv_0)$. Then,  for a.e. $t\geq0$,
$$
\left\{
\ba{l}
\left\|\GC^\eps(t)-\bfy^\eps(t)\right\|
\,\leq\, C\,\eps^3\,e^{C\,\eps\,t}
\,\Big(1+t\,(1+\|\bfv_0\|\,+\,t)\Big)\,(\|\bfv_0\|\,+\,t),
\\[0.95em]
\left\|\eGC^\eps(t)-\e^\eps(t)\right\|
\,\leq\, C\,\eps^2\,e^{C\,\eps\,t}\,
\Big(\,1+(1+t)\,(\|\bfv_0\|\,+\,t)+t\,(\|\bfv_0\|\,+\,t)^2\,\Big)\,(\|\bfv_0\|\,+\,t),
\ea\right.
$$
where $(\GC^\eps,\eGC^\eps)$ is as in \eqref{GC} and \eqref{toy-egc} and $(\bfy^\eps,\e^\eps)$ solves
\be
\label{eq:mhv}
\left\{
\begin{array}{l}
\ds\frac{\dD\bfy^\eps}{\dD t}\,=\, \eps\,\EcB(t,\bfy^\eps),
\\[1.1em]
\ds\frac{\dD \e^\eps}{\dD t}\,=\,-\eps\,\e^\eps\,\Div_\bfx(\EcB)(t,\bfy^\eps),
\end{array}\right.
\ee
with initial data $\bfy^\eps(0)= \GC^\eps(0)$ and  $\e^\eps(0)=\eGC^\eps(0)$.
\epr

\br
Note that if $\bE$ derives from a potential, that is, if $\bE$ is curl-free, then the equation on $\e^\eps$ is trivial since $\Div_\bfx (\EcB)=0$. Yet this cancellation does not improve any convergence rate. Incidentally we point out that in this case $\dD_\bfx\bE$ is symmetric so that the cancellation follows at a more abstract level from computations of Remark~\ref{trace-free}.
\er

\begin{proof}
We consider $(\bfx^\eps,\bfv^\eps)$ the  solution to \eqref{eq:xv-hom} starting from $(\bfx_0,\bfv_0)$ and the corresponding $(\bfy^\eps,\e^\eps)$, the solution to \eqref{eq:mhv} starting from $(\GC^\eps(0),\eGC^\eps(0))$, where $(\GC^\eps,\eGC^\eps)$ is as in \eqref{GC}-\eqref{toy-egc}. First, to ease comparisons, we recall that $\e(\bfv^\eps)=\|\bfv^\eps\|^2/2$ (and ban temporarily the confusing shorthands $\e^\eps$ and $\e$ for $\e(\bfv^\eps)$) and write \eqref{toy-e2} as
\begin{eqnarray*}
\frac{\dD}{\dD t}\left[\eGC^\eps+\frac{\eps^2}{B^2}\,\chi^\eps\right]
&=&\,-\eps\,\eGC^\eps\,\Div_\bfx\EcB(t,\bfy^\eps)\,+\,\frac{\eps^2}{B^2}\,\eta^\eps
\\
&& -\eps\,\left(\e(\bfv^\eps)-\eGC^\eps\right)\,\Div_\bfx\EcB(t,\bfy^\eps)
\\
&&-\eps\,\e(\bfv^\eps)\,\left(\Div_\bfx\EcB(t,\bfx^\eps)-\Div_\bfx\EcB(t,\GC^\eps)\right)
\\
&&-\eps\,\e(\bfv^\eps)\,\left(\Div_\bfx\EcB(t,\GC^\eps)-\Div_\bfx\EcB(t,\bfy^\eps)\right).
\end{eqnarray*}
Then, from subtracting the latter equation to the one for $\e^\eps$ in \eqref{eq:mhv} stems for a.e. $t\geq0$,
\begin{eqnarray*}
|\eGC^\eps(t)-\e^\eps(t)|\
&\leq&C\,\eps\,\int_0^t\,|\eGC^\eps(s)-\e^\eps(s)|\dD s\,
+\,C\,\eps^2\,\left(|\chi^\eps(0)|+|\chi^\eps(t)|\right)\,\\
&&+\,C\,\eps^2\,\int_0^t\,|\eta^\eps(s)|\dD s
\,+\,C\,\eps^2\,\int_0^t\,\|\bfv^\eps(s)\|\dD s\\
&&+\,C\,\eps^2\,\int_0^t\,\|\bfv^\eps(s)\|^3\dD s
+C\eps\,\int_0^t\,\|\bfv^\eps(s)\|^2\,\|\GC^\eps(s)-\bfy^\eps(s)\|\dD s,
\end{eqnarray*}
where $C$ depends polynomially on $\|\bE\|_{W^{2,\infty}}$ and $B^{-1}$. Finally the estimate on $\eGC^\eps-\e^\eps$ follows from Lemma~\ref{toy-bnd0}, Proposition~\ref{toy-cheap}  and the Gr\"onwall lemma.
\end{proof}

\subsection{Long-time asymptotics}
\label{s:toy-long}

As aforementioned, the fact that vector fields appearing in the leading-order asymptotics seem to be $\cO(\eps)$ suggests that it should also be possible to validate asymptotics for $(\bfx^\eps,\e(\bfv^\eps))$ on time intervals of length $\cO(\eps^{-1})$ with convergence rates $\cO(\eps)$. To carry out such achievement we need to refine bounds from Lemma~\ref{toy-bnd0}.

\bl
\label{toy-bnd-long} 
There exists a universal positive constant $C$ such that any solution to \eqref{eq:xv-hom} starting from $(\bfx_0,\bfv_0)$ satisfies for any $t\geq0$
$$
\left\{
\ba{l}
\ds\|\bfx(t)\|\,\leq\, \|\bfx_0\|
+\,\frac{\eps}{B}\,t\,\|\bE\|_{L^\infty}
\,+\,C\,\eps\,e^{C\,\frac{\eps}{B}\,t\,\|\bE\|_{\dot{W}^{1,\infty}}}
\left(1\,+\,\|\bfv_0\|\,+\,\frac{\eps}{B}\,\|\bE\|_{L^\infty}\right),
\\[0.95em]
\ds\|\bfv(t)\|\,\leq\, C\,e^{C\,\frac{\eps}{B}\,t\,\|\bE\|_{\dot{W}^{1,\infty}}}
\left(1\,+\,\|\bfv_0\|\,+\,\frac{\eps}{B}\,\|\bE\|_{L^\infty}\right)\,.
\ea\right.
$$
\el
\begin{proof}
Integrating \eqref{toy-e1} yields for a.e. $t\geq0$
\begin{eqnarray*}
\e(\bfv)&\leq&  \e(\bfv_0) \,+\,\frac{\eps}{B}\|\bE\|_{L^\infty}\|\bfv_0\|
\,+\,\frac{\e(\bfv)}{2}\,+\,\frac{\eps^2}{2B^2}\|\bE\|_{L^\infty}^2\\
&+&\frac{\eps}{B}t\,\|\d_t\bE\|_{L^\infty}
+\frac{\eps}{B}(2\|\dD_\bfx\bE\|_{L^\infty}+\|\d_t\bE\|_{L^\infty})\int_0^t
    \e(s)\,\dD s,
\end{eqnarray*}
which after a few algebraic manipulations and an application of the Gr\"onwall lemma proves the claim on $\|\bfv\|$. The bound on $\bfx$ is obtained simarly by integrating \eqref{toy-x1}.
\end{proof}

We now focus on large-time asymptotics for $(\bfx,\e(\bfv))$.

\bpr
\label{toy-long-xe}
Assume  $\bE\in W^{2,\infty}$. There exists a constant $C>0$, depending polynmially on $\|\bE\|_{W^{2,\infty}}$ and $B^{-1}$, such that any solution to~\eqref{eq:xv-hom} starting from $(\bfx_0,\bfv_0)$ satisfies for a.e. $t\geq0$
$$
\left\{
\ba{l}
\ds\|\bfx^\eps(t)-\bfy^\eps(t)\|
\,\leq\; C\,\eps\,e^{C\,\eps\,t}
\,(1+\eps+\|\bfv_0\|),
\\[0.95em]
\ds\|\e(\bfv^\eps(t))-\e^\eps(t)\| \,\leq\, C\,\eps\,e^{C\,\eps\,t}\,(1+\eps+\|\bfv_0\|)^3,
\ea\right.
$$
where $\e(\bfv)=\tfrac12\|\bfv\|^2$ and
$(\bfy^\eps,\e^\eps)$ solves (\ref{eq:mhv}) with initial data
$\bfy^\eps(0)=\bfx_0$ and $\e^\eps(0)\,=\,\tfrac12\|\bfv_0\|^2$.
\epr
\begin{proof}
The estimate on $\bfx^\eps-\bfy^\eps$ follows  from Lemma~\ref{toy-bnd-long} and the Gr\"onwall lemma after an integration of \eqref{toy-x1}. To proceed, we use \eqref{toy-e2} in the form
\begin{eqnarray*}
&&\frac{\dD}{\dD
  t}\left[\e(\bfv^\eps)\,-\,\eps\,\left\langle\EcB(t,\bfx^\eps),\bfv^\eps\right\rangle+\frac{\eps^2}{B^2}\,\chi^\eps\right]
\\
&&=\,-\eps\,\e(\bfv^\eps)\,\Div_\bfx\EcB(t,\bfy^\eps)\,+\,\frac{\eps^2}{B^2}\,\eta^\eps
\,-\,\eps\,\e(\bfv^\eps)\,\left(\Div_\bfx\EcB(t,\bfx^\eps)-\Div_\bfx\EcB(t,\bfy^\eps)\right),
\end{eqnarray*}
thus, for a.e. $t\geq0$,
\begin{eqnarray*}
|\e(\bfv^\eps(t))-\e^\eps(t)|\
&\leq&\,C\eps\,\int_0^t\,|\e(\bfv^\eps(s))-\e^\eps(s)|\,\dD s\,
\\
&&+\,C\,\eps\,(\|\bfv_0\|+\|\bfv^\eps(t)\|)+C\eps^2\,(|\chi^\eps(0)|+|\chi^\eps(t)|)\,
\\
&& +C\,\eps^2\,\int_0^t\,|\eta^\eps(s)|\,\dD s
+C\eps\,\int_0^t\,\|\bfv^\eps(s)\|^2\,\|\bfx^\eps(s)-\bfy^\eps(s)\|\,\dD s,
\end{eqnarray*}
where $C$ depends polynomially on $\|\bE\|_{W^{2,\infty}}$ and $B^{-1}$. One may conclude again with Lemma~\ref{toy-bnd-long} and the Gr\"onwall lemma.
\end{proof}

\br
\label{rk:toy-scaled-time}
The proof also yields the analysis of dynamics involving fields depending on $\eps$ but satisfying bounds uniform with respect to $\eps$. In particular the result may be extended without change to the case where $\bE^\eps(t,\bfx)=\bE(\eps\,t,\bfx)$, $0<\eps\lesssim 1$. This somehow simpler problem is the one classically considered because then the asymptotic dynamics is essentially independent of $\eps$ at leading order since
$$
(\bfy^\eps,\e^\eps)(t)\,=\,(\bfy,\e)(\eps\,t),
$$
with $(\bfy,\e)$ independent of $\eps$.
\er

Of course we may also use Lemma~\ref{toy-bnd-long} to refine time dependences in Proposition~\ref{toy-xe} so as to fill the gap concerning what happens at leading-order for intermediate times $1\lesssim t\lesssim \eps^{-1}$. For possible external reference let us store without proof the corresponding result.

\bpr
\label{toy-bis-xe}
Assume  $\bE\in W^{2,\infty}$. Then there exists a constant $C>0$, depending
polynomially on $\|\bE\|_{W^{2,\infty}}$ and $B^{-1}$, such that the following holds. Let
$(\bfx^\eps,\bfv^\eps)$ be the solution to~\eqref{eq:xv-hom} starting from
$(\bfx_0,\bfv_0)$. Then, for a.e. $t\geq0$,
$$
\left\{\ba{l}
\|\GC^\eps(t)-\bfy^\eps(t)\| \,\leq\, C\,\eps^3\,e^{C\,\eps\,t}
\,\Big[1+t\,(1+\eps+\|\bfv_0\|)\Big]\,(1+\eps+\|\bfv_0\|),
\\[0.95em]
\|\eGC^\eps(t)-\e^\eps(t)\|\,\leq\, C\,\eps^2\,e^{C\,\eps\,t}\,
\Big[1+t\,(1+\eps+\|\bfv_0\|)\Big]\,(1+\eps+\|\bfv_0\|)^2,
\ea\right.
$$
where $(\GC^\eps,\eGC^\eps)$ is as in \eqref{GC} and \eqref{toy-egc} and $(\bfy^\eps,\e^\eps)$ solves \eqref{eq:mhv} with initial data $\bfy^\eps(0)=\GC^\eps(0)$ and
$\e^\eps(0)=\eGC^\eps(0)$. 
\epr

\subsection{PDE counterparts}

Now let us translate the foregoing results at the PDE level. 

On the reduced phase-space where $\bfZ=(\bfy,\e)$ lives the relevant macroscopic velocity is $\eps\cW_1(t,\bfZ)$ where
\be
\label{eq:gc2d-1}
\cW_1(t,\bfZ)
\,=\, \bp
\ba{l}
\ds\EcB(t,\bfy)
\\[0.75em]
-\e \,\Div_{\bfx}\EcB(t,\bfy)
\ea
\ep\,,
\ee
which corresponds to the velocity field of system \eqref{eq:mhv} defining the characteristic curves of the equation
\be
\label{eq:gc2d-2}
\d_t G^\eps\,+\,\eps\,\Div_\bfZ \left(\cW_1\,G^\eps\right)\,=\,0.
\ee

With this in hands we may deduce from Proposition~\ref{toy-long-xe} and Proposition~\ref{toy-bis-xe} the following statement, where we have made explicit push-forwards that were easy to compute.

\bt
\label{toy-density}
Let $\bE\in W^{2,\infty}$. There exists a constant $C$ depending polynomially on $\|\bE\|_{W^{2,\infty}}$ and $B^{-1}$ such that the following holds for any solution $f^\eps$ to \eqref{eq:vla2d} with initial datum a nonnegative density $f_0$.\\
{\bf (i).} \emph{Long-time first-order asymptotics}. $F^\eps$ defined by
$$
F^\eps(t,\bfx,\e)=\int_0^{2\pi}\,f^\eps(t,\bfx,\sqrt{2\,\e}\,\beD(\theta))\dD \theta
$$
satisfies for a.e. $t\geq0$
$$
\|F^\eps(t,\cdot)-G^\eps(t,\cdot)\|_{\dot{W}^{-1,1}}\leq C\,\eps\,e^{C\,\eps\,t}\,\int_{\R^2\times\R^2}(1+\eps+\|\bfv\|)^3\,f_0(\bfx,\bfv)\,\dD\bfx\,\dD\bfv
$$
where $G^\eps$ solves \eqref{eq:gc2d-2} with initial datum $G_0=F^\eps(0,\cdot)$ given by
$$
G_0(\bfx,\e)=\int_0^{2\pi}\,f_0(\bfx,\sqrt{2\,\e}\,\beD(\theta))\dD \theta\,.
$$
{\bf (ii).} \emph{Short-time second-order asymptotics}. The push-forwards $F_{\rm gc}^\eps(t,\cdot)$ of $f^\eps(t,\cdot)$ by the maps
$$
(\bfx,\bfv)\,\mapsto\,\left(\bfx+\frac{\eps}{B}\,\bJ\bfv\,,\,
\tfrac12\|\bfv\|^2\,-\,\eps\,\left\langle\EcB(t,\bfx),\bfv\right\rangle\right)
$$
satisfy for a.e. $t\geq0$
$$
\|F_{\rm gc}^\eps(t,\cdot)-G_{\rm gc}^\eps(t,\cdot)\|_{\dot{W}^{-1,1}}\leq C\,\eps^2\,e^{C\,\eps\,t}\,(1+t)\,\int_{\R^2\times\R^2}(1+\eps+\|\bfv\|)^3\,f_0(\bfx,\bfv)\,\dD\bfx\,\dD\bfv
$$
where $G_{\rm gc}^\eps$ solves \eqref{eq:gc2d-2} with initial datum $G_0=F_{\rm gc}^\eps(0,\cdot)$, defined as the push-forward of $f_0$ by the map
$$
(\bfx,\bfv)\,\mapsto\,\left(\bfx+\frac{\eps}{B}\,\bJ\bfv\,,\,
\tfrac12\|\bfv\|^2\,-\,\eps\,\left\langle\EcB(0,\bfx),\bfv\right\rangle\right)\,.
$$
\et
\begin{proof}
The first result is a direct consequence of the abstract Proposition~\ref{p:ODEtoPDE} and the estimates provided in Proposition~\ref{toy-long-xe} on the characteristic curves. The second one follows the same lines with the help of Proposition \ref{toy-bis-xe} instead of Proposition~\ref{toy-long-xe}.   
\end{proof}

Due to the special structure of the homogeneous two-dimensional case, with essentially the same proof one may also provide versions focusing only on the spatial variables and its $\eps$-corrections. This involves the asymptotic equation
\be\label{eq:gc2d-2x}
\d_t r^\eps\,+\,\eps\,\Div_\bfy\left(r^\eps\,\EcB\right)
\,=\,0\,.
\ee
 
\bpr
Let $\bE\in W^{2,\infty}$. There exists a constant $C$ depending polynomially on $\|\bE\|_{W^{2,\infty}}$ and $B^{-1}$ such that the following holds for any solution $f^\eps$ to \eqref{eq:vla2d} with initial datum a nonnegative density $f_0$.\\
{\bf (i).} \emph{Long-time first-order asymptotics}. $\rho^\eps$ defined by
$$
\rho^\eps(t,\bfx)=\int_{\R^2}\,f^\eps(t,\bfx,\bfv)\dD \bfv
$$
satisfies for a.e. $t\geq0$
$$
\|\rho^\eps(t,\cdot)-r^\eps(t,\cdot)\|_{\dot{W}^{-1,1}}\leq C\,\eps\,e^{C\,\eps\,t}\,\int_{\R^2\times\R^2}(1+\eps+\|\bfv\|)\,f_0(\bfx,\bfv)\,\dD\bfx\,\dD\bfv
$$
where $r^\eps$ solves \eqref{eq:gc2d-2x} with initial datum $r_0=\rho^\eps(0,\cdot)$ given by
$$
r_0(\bfx)=\int_{\R^2}\,f_0(\bfx,\bfv)\dD \bfv\,.
$$
{\bf (ii).} \emph{Short-time third-order asymptotics}. $\rho_{\rm gc}^\eps$ defined by
$$
\rho_{\rm gc}^\eps(t,\bfy)=\int_{\R^2}\,f^\eps\left(t,\bfy-\frac{\eps}{B}\,\bJ\bfv,\bfv\right)\dD\bfv
$$
satisfies for a.e. $t\geq0$
$$
\|\rho_{\rm gc}^\eps(t,\cdot)-r_{\rm gc}^\eps(t,\cdot)\|_{\dot{W}^{-1,1}}\leq C\,\eps^3\,e^{C\,\eps\,t}\,(1+t)\,\int_{\R^2\times\R^2}(1+\eps+\|\bfv\|)^2\,f_0(\bfx,\bfv)\,\dD\bfx\,\dD\bfv
$$
where $r_{\rm gc}^\eps$ solves \eqref{eq:gc2d-2x} with initial datum $(r_{\rm gc}^\eps)_0=\rho_{\rm gc}^\eps(0,\cdot)$, given by
$$
(r_{\rm gc}^\eps)_0(\bfy)=\int_{\R^2}\,f_0\left(\bfy-\frac{\eps}{B}\,\bJ\bfv,\bfv\right)\dD\bfv\,.
$$
\epr

\section{General three-dimensional case}\label{s:general}

We come back to the three-dimensional system
\be
\label{e:xv}
\left\{
\ba{l}
\ds\frac{\dD\bfx}{\dD t}\,=\,\bfv\,,
\\[0.95em]
\ds\frac{\dD\bfv}{\dD t}\,=\,\frac{\bfv\wedge \bB(t,\bfx)}{\eps} \,+\,\bE(t,\bfx)\,,
\ea
\right.
\ee
and follow the pattern of the short-time analysis of Section~\ref{s:homogeneous}. As there we do not mark $\eps$-dependences as long as no confusion is possible.

\subsection{Slow variables and uniform bounds}

First, Lemma~\ref{toy-bnd0} stands without change in its statement or its proof.

\bl\label{bnd0}
Solutions to \eqref{e:xv} starting from $(\bfx_0,\bfv_0)$ are defined
globally in time and satisfy for any $t\geq0$
$$
\left\{\ba{l}
\|\bfx(t)\|\,\leq\, \|\bfx_0\|
+\,t\,\|\bfv_0\|+\,t^2\,\|\bE\|_{L^\infty},
\\[0.95em]
\|\bfv(t)\| \,\leq\, \|\bfv_0\|+2\,t\,\|\bE\|_{L^\infty}\,.
\ea\right.
$$
\el

Here some geometric preparation is needed to identify some set of slow variables. At leading order the fast motion is locally a rotation of $\bfv$ around $\eDpar(t,\bfx)$ where we recall that $\eDpar$ is defined through
$$
B(t,\bfx)=\|\bB(t,\bfx)\|\,,\qquad \bB(t,\bfx)\,=\,B(t,\bfx)\,\eDpar(t,\bfx)\,.
$$
As aforementioned this naturally suggests first a separation of $\bfv$ between a component aligned on $\eDpar(t,\bfx)$, $\vpar(t,\bfx,\bfv)\,\eDpar(t,\bfx)$, and a perpendicular component $\bvperp(t,\bfx,\bfv)$, and second by mimicking the homogeneous case the introduction of a kinetic energy variable associated with $\bvperp(t,\bfx,\bfv)$, $\eperp(t,\bfx,\bfv)=\tfrac12\|\bvperp(t,\bfx,\bfv)\|^2$\,.

We recall that the above decomposition is explicitly given as 
$$
\left\{
\begin{array}{l}
\vpar(t,\bfx,\bfv) \,=\,\langle\bfv\,,\,\eDpar(t,\bfx)\rangle,
\\[0.95em]
\bvperp(t,\bfx,\bfv)=\bfv-\vpar(t,\bfx,\bfv)\,\eDpar(t,\bfx),
\end{array}\right.
$$
and that correspondingly we introduce the decomposition of the electric field\footnote{That is, $\Epar(t,\bfx)=\vpar(t,\bfx,\bE(t,\bfx))$, $\Eperp(t,\bfx)=\bvperp(t,\bfx,\bE(t,\bfx))$.} $\bE=\Epar\,\eDpar+\Eperp$,
$$
\Epar(t,\bfx,\bfv)=\left\langle\bE(t,\bfx)\,,\,\eDpar(t,\bfx)\right\rangle,\qquad\qquad
\Eperp(t,\bfx,\bfv)=\bE(t,\bfx)-\Epar(t,\bfx,\bfv)\,\eDpar(t,\bfx)\,.
$$
Both to ease computations and to emphasize analogies with the two-dimensional case it is expedient to introduce, for any $\bfx\in\R^3$, the linear operator $\bJ(t,\bfx)$ defined as
$$
\bJ(t,\bfx) \,\bfa\,=\, \bfa\wedge \eDpar(t,\bfx)\,.
$$
Going on with geometric considerations, we note that the following simple relations are of pervasive use in latter computations:
\be\label{e:J}
\left\{
\begin{array}{l}
\bJ(t,\bfx)\eDpar(t,\bfx)\,=\,0\,,\qquad\,\,
\eDpar(t,\bfx)\cdot\bJ(t,\bfx)\bfa=0\,,
\\[0.95em]
\bJ(t,\bfx)^2\,\bfa\,=\,-\bvperp(t,\bfx,\bfa)\,,\qquad\,\,
\bJ(t,\bfx)^*=-\bJ(t,\bfx)\,,
\end{array}
\right.
\ee
and $\eDpar(t,\bfx)\cdot\d_t\eDpar(t,\bfx)\,=\,0$, $\eDpar(t,\bfx)\cdot\dD_\bfx\eDpar(t,\bfx)\,\bfa\,=\,0$.

For the sake of concision, but somewhat inconsistently, from now on we shall use the shorthand $\vpar(t)$ for $\vpar(t,\bfx(t),\bfv(t))$ and similarly for $\bvperp$ and $\eperp$. We shall also identify functions of $(\bfx,\bfv)$ with functions of $(\bfx,\vpar,\bvperp)$. Then, we may split \eqref{e:xv} as
\be
\label{e:1}
\left\{
\ba{l}\ds
\frac{\dD\bfx}{\dD t} \,=\,\bfv,
\\[0.95em]
\ds\frac{\dD\vpar}{\dD t}\,=\,\Epar(t,\bfx)\,+\,\left\langle \bvperp ,
\d_t\eDpar(t,\bfx)+\dD_{\bfx}\eDpar(t,\bfx)\bfv \right\rangle,
\\[0.95em]
\ds\frac{\dD\eperp}{\dD t}\,=\, \left\langle \Eperp(t,\bfx)\,-\,\vpar\,\left(\d_t\eDpar(t,\bfx)\,+\,\dD_{\bfx}\eDpar(t,\bfx)\,\bfv\right),  \bvperp\right\rangle,
\ea \right.
\ee
and
\be
\label{e:2}
\ds\frac{\dD\bvperp}{\dD t}\,=\,\ds
\frac{B(t,\bfx)}{\eps}\, \bJ(t,\bfx)\,\bvperp
\,+\, \bF(t,\bfx,\bfv), 
\ee
where $\bfv=\vpar\eDpar+\bvperp$ and the force field $\bF$ is
$$
\bF(t,\bfx,\bfv)=\,\bF_0(t,\bfx,\vpar)+\bF_1(t,\bfx,\bfv)+\bF_2(t,\bfx,\bfv),
$$
with $\bF_1$ depending linearly on $\bvperp$, $\bF_2$ quadratic in $\bvperp$, explicitly
\be
\label{def:F}
\left\{\ba{lll}
\ds\bF_0(t,\bfx,\vpar)  &=  & \ds \Eperp(t,\bfx)\,-\,
\vpar\left(\d_t\eDpar(t,\bfx)\,+\,\vpar\,\dD_{\bfx}\eDpar(t,\bfx)\,\eDpar(t,\bfx)\right)\,,
\\[0.95em]\ds
\bF_1(t,\bfx,\bfv)  &=  & \ds -\,
\left\langle\d_t\eDpar(t,\bfx)\,+\,\vpar\,\dD_{\bfx}\eDpar(t,\bfx)\,\eDpar(t,\bfx)\,,\,
  \bvperp
\right\rangle\,\eDpar(t,\bfx)
\\[0.95em]\ds
\,& \, & \ds-\,\vpar\,\dD_{\bfx}\eDpar(t,\bfx)\,\bvperp\,,
\\[0.95em]
\ds\bF_2(t,\bfx,\bvperp) & = & \ds -\, \left\langle
\dD_{\bfx}\eDpar(t,\bfx)\,\bvperp\,,\, \bvperp \right\rangle\,\eDpar(t,\bfx)\,.
\ea\right.
\ee

\br
As already pointed out in Remark~\ref{rddnt-hom} along the analysis of the homogeneous case, it is convenient to work with a formulation containing some redundancy such as \eqref{e:1}-\eqref{e:2}. Indeed, here, to suppress the apparent overdetermination one could for instance replace \eqref{e:2} with an equation for an angle of $\bvperp$ but then one loses track of an important property of System~\eqref{e:xv}: at principal order oscillations are linear in $\bvperp$. In contrast, as already apparent in the homogeneous case or in the splitting of $\bF$, all our algebraic manipulations will be organized by the degree of linearity in $\bvperp$.
\er

\subsection{Elimination of linear terms}

A direct consequence of Lemma~\ref{bnd0} and \eqref{e:1} is that $(\bfx,\vpar,\eperp)$ in $W^{1,\infty}_{loc}$ and $\bvperp$ in $L_{loc}^\infty$ are uniformly bounded with respect to $\eps$. This is sufficient to extract converging sequences but not to take limits in the equations because of the nonlinearity in $\bvperp$.

Instead, to proceed, we begin an uncoupling process similar to the one carried out in Section~\ref{s:homogeneous}. Elimination, at leading order, of linear terms in $\bfv$ is summarized as

\bl
\label{1st}
For any $\bL\in W^{1,\infty}(\R^+_t; \,\cL_1(\R^3,\R^p))$, $p\in\N^*$,
solutions $(\bfx,\bfv)$ to \eqref{e:xv} satisfy for a.e. $t\geq 0$,
$$
\bL(t)\bvperp \,=\,-\eps\frac{\dD}{\dD t}\left[\bL(t)\left(\frac{\bJ\bvperp}{B}\right)\right]
\,+\,\eps\, \bL'(t)\left(\frac{\bJ\bvperp}{B}\right)
\,+\,\eps \bL(t)\bU,
$$ 
with the macroscopic velocity $\bU$ given by
\be
\label{U}
\bU(t,\bfx,\bfv)\,=\, \frac{\bJ\,\bF}{B}
+\left[\d_t\left(\frac{\bJ}{B}\right)
\,+\,\dD_\bfx\left(\frac{\bJ}{B}\right)\,\bfv\right]\,\bvperp\,.
\ee
\el

\begin{proof}
Applying $\frac{\eps\bJ}{B}$ to \eqref{e:2} and combining with the first line of \eqref{e:1} yields
\be\label{e:homologic-1st}
\bvperp \,=\,-\eps\frac{\dD}{\dD t}\left[\frac{\bJ\,\bvperp}{B(\bfx)}\right]
\,+\,\eps\, \bU(t,\bfx,\bfv)\,.
\ee
Then the result follows from the chain rule. 
\end{proof}

Note that the macroscopic velocity $\bU$ is split according to degree
in $\bvperp$ as $\bU=\bU_{10}+\bU_{11}+\bU_{12}$ where  $\bU_{10}$
contains terms which do not depend on $\bvperp$, 
\begin{eqnarray*}
\bU_{10}(t,\bfx,\vpar) &=&
                           \frac{\bJ(t,\bfx)}{B(t,\bfx)}\,\bF_{0}(t,\bfx,\vpar)
\\
&= &\EcB(t,\bfx)\,+\,\vpar^2\,\curvB(t,\bfx) \,+\,\vpar\,\dtB(t,\bfx)
\\
&= &\EcB(t,\bfx)\,+\,\vpar\,\Sig(t,\bfx,\vpar),
\end{eqnarray*}
which corresponds to the classical drifts defined in (\ref{drift-0})-(\ref{drift-2}), whereas $\bU_{11}$ is given by 
$$
\bU_{11}(t,\bfx,\vpar,\bvperp)
\,=\, \frac{\bJ\,\bF_1}{B}(t,\bfx,\vpar,\bvperp)
\,+\,\left[\d_t\left(\frac{\bJ}{B}\right) \,+\,\vpar\dD_\bfx\left(\frac{\bJ}{B}\right)\,\eDpar\right]\,\bvperp\,
$$
and observing that $\bJ\,\bF_2=0$, we have for $\bU_{12}$,
\begin{eqnarray*}
\bU_{12}(t,\bfx,\bvperp)
\,=\,\left[\dD_\bfx\left(\frac{\bJ}{B}\right)(t,\bfx)\,\bvperp\right]\,\bvperp\,.
\end{eqnarray*}

Let us anticipate that the partial elimination of $\bU_{12}$ will give a contribution known as the grad-$B$ drift and that encodes the influence of the variations of the intensity $B$ on the spatial trajectory. 

With Lemma~\ref{1st}, at leading-order in $\eps$ one may eliminate from \eqref{e:1} terms that are linear in $\bvperp$. We first treat the first equation in \eqref{e:1} by applying Lemma~\ref{1st} with $\bL(t)\bvperp=\bvperp$, which reduces to \eqref{e:homologic-1st}. This leads to 
\be\label{e:x}
\frac{\dD \GC}{\dD t}\,=\,\vpar\,\eDpar(t,\bfx) \,+\,\eps\,\bU(t,\bfx,\bfv),
\ee
where $\bU$ is as in \eqref{U} and we have introduced the so-called
guiding center already defined in (\ref{GC}). 

Then we consider the
second equation in \eqref{e:1} and apply Lemma~\ref{1st} with
$$
\bL(t)\bvperp=\langle\bvperp,\d_t\eDpar(t,\bfx(t))+\vpar(t)\dD_{\bfx}\eDpar(t,\bfx(t))\eDpar(t,\bfx(t))\rangle,
$$
to remove the linear part with respect to $\bvperp$ in the right hand
side and to derive an equation for a first correction of the parallel velocity,
\begin{align}
\label{e:vp}
\frac{\dD }{\dD t}\left[\vpar
\,+\,\frac{\eps}{B}\,\left\langle\bJ\,\bvperp\,,\,\d_t\eDpar+\vpar\,\dD_{\bfx}\eDpar\,\eDpar\right\rangle\right] \,=\,\Epar\,+\,\left\langle\bvperp\,,\,\dD_{\bfx}\eDpar\,\bvperp\right\rangle
+\eps \,u_1(t,\bfx,\bfv),
\end{align}
where $u_1=u_{10}+u_{11}+u_{12}+u_{13}$ is obtained from
\begin{eqnarray*}
u_{10}&=&\langle\bU_{10},\d_t\eDpar+\vpar\,\dD_{\bfx}\eDpar\,\eDpar\rangle
\\
&=&\frac{1}{B}\,\left\langle {\bJ}\,\Eperp\,,\,\d_t\eDpar+\vpar\,\dD_{\bfx}\eDpar\,\eDpar\right\rangle
\\
&=& \left\langle \bE\,,\, \dtB + \vpar \curvB\right\rangle
\,=\,\left\langle \bE\,,\,\Sig\right\rangle,
\end{eqnarray*}
where  $\curvB$ and $\dtB$ are given in \eqref{drift-0}-\eqref{drift-1} and $u_{11}$ is
\begin{eqnarray*}
\ds u_{11} &=&\left\langle \bU_{11},
\d_t\eDpar+\vpar\,\dD_{\bfx}\eDpar\,\eDpar\right\rangle
\,+\,\frac{\Epar}{B}\,\left\langle {\bJ}\,\bvperp, \dD_{\bfx}\eDpar\eDpar\right\rangle
\\  
&+&\ds\frac{1}{B}\, \left\langle
\bJ\bvperp,\d_t^2\eDpar+\vpar\,\d_t(\dD_{\bfx}\eDpar)\eDpar+\vpar\,\dD_{\bfx}\eDpar\d_t\eDpar\right\rangle
\\
&+&\ds\frac{\vpar}{B}\, \left\langle
\bJ\bvperp,\left[\d_t(\dD_\bfx\eDpar)+\vpar\,\dD_{\bfx}^2\eDpar\eDpar+\vpar\,\left(\dD_{\bfx}\eDpar\right)^2\right]\eDpar\right\rangle\,,
\end{eqnarray*}
whereas the last terms $(u_{12},u_{13})$ are
\begin{eqnarray*}
u_{12}&=&  \ds\left\langle\bU_{12}(t,\bfx,\bfv)\,,\,\d_t\eDpar+\vpar\,\dD_{\bfx}\eDpar\,\eDpar\right\rangle
\\ 
&+&\ds\frac{1}{B}\, \left\langle
\bJ\bvperp,\left[\d_t(\dD_\bfx\eDpar)+\vpar\,\dD_{\bfx}^2\eDpar\eDpar+\vpar\,\left(\dD_{\bfx}\eDpar\right)^2\right]\bvperp \right\rangle
\\ 
&+& \ds\frac{1}{B} \, \left\langle \bvperp,\d_t\eDpar+\vpar\,\dD_{\bfx}\eDpar\,\eDpar\right\rangle\,
\left\langle \bJ\bvperp\,,\, \dD_{\bfx}\eDpar\,\eDpar\right\rangle\,,
\end{eqnarray*}
and
$$
u_{13}\,=\,
\frac{1}{B} \, \left\langle \bvperp, \dD_{\bfx}\eDpar\,\bvperp\right\rangle\,
\left\langle{\bJ}\,\bvperp\,,\, \dD_{\bfx}\eDpar\,\eDpar\right\rangle\,.
$$

Finally we conclude the elimination of linear terms by reformulating the third equation in \eqref{e:1}. To proceed we apply Lemma~\ref{1st} with
\begin{align*}
\bL(t)\bvperp&=\left\langle \bvperp\,,\,\bF_0(t,\bfx(t),\vpar(t))\right\rangle\\[0.5em]
&=\left\langle \bvperp\,,\,\bE(t,\bfx(t))
\,-\,\vpar(t)\,\d_t\eDpar(t,\bfx(t))
\,-\,(\vpar(t))^2\,\dD_{\bfx}\eDpar(t,\bfx(t))\, \eDpar(t,\bfx(t))\right\rangle
\end{align*}
and naturally obtain an equation for a first correction of the  kinetic energy in
the perpendicular plan to the magnetic field,
\be
\label{e:ep}
\frac{\dD}{\dD t}\left[\eperp \,+\, \frac{\eps}{B}\,\left\langle
{\bJ}\,\bvperp,\bF_0\right\rangle
\right]  \,=\,
-\vpar\left\langle\bvperp,
\dD_{\bfx}\eDpar(t,\bfx)\bvperp\right\rangle+\eps \,d_1(t,\bfx,\bfv),
\ee
where $d_1=d_{11}+d_{12}+d_{13}$ is obtained from
\begin{eqnarray*}
d_{11} &=& \left\langle \bU_{11}\,,\,\bF_0\right\rangle
\,-\,\ds  \frac{\Epar}{B}\,\left\langle
{\bJ}\bvperp\,,\,\d_t\eDpar+2\vpar\,\dD_{\bfx}\eDpar\,\eDpar\right\rangle
\\ 
&+&  \left\langle
\frac{\bJ}{B}\bvperp\,,\,\d_t\bE\,+\, \vpar\dD_\bfx\bE\eDpar\right\rangle
\\ 
&-&\ds\frac{\vpar}{B}\, \left\langle
\bJ\bvperp,\d_t^2\eDpar+\vpar\,\d_t(\dD_{\bfx}\eDpar)\eDpar+\vpar\,\dD_{\bfx}\eDpar\d_t\eDpar\right\rangle
\\
&-& 
\ds\frac{\vpar^2}{B}\, \left\langle\bJ\bvperp\,,\,
\left[\d_t(\dD_\bfx\eDpar)+\vpar\,\dD_{\bfx}^2\eDpar\eDpar+\vpar\,\left(\dD_{\bfx}\eDpar\right)^2\right]\eDpar
\right\rangle\,,
\end{eqnarray*}
and $d_{12}$ is given by
\begin{eqnarray*}
d_{12}&=& 
\left\langle \bU_{12}\,,\,\bF_0\right\rangle
\,+\,
\left\langle\frac{\bJ}{B}\bvperp\,,\,\dD_\bfx\bE\bvperp\right\rangle
\\ 
&-&
\frac{1}{B}
\, \langle \bvperp,\d_t\eDpar+\vpar\,\dD_{\bfx}\eDpar\,\eDpar\rangle\,
\left\langle\bJ\bvperp\,,\,\d_t\eDpar+2\vpar\,\dD_{\bfx}\eDpar\,\eDpar\right\rangle
\\ 
&-& \ds\frac{\vpar}{B}\, \left\langle \bJ\bvperp\,,\,
\left[\d_t(\dD_\bfx\eDpar)+\vpar\,\dD_{\bfx}^2\eDpar\eDpar+\vpar\,\left(\dD_{\bfx}\eDpar\right)^2\right]\bvperp \right\rangle\,,
\end{eqnarray*}
and
$$
d_{13}\,=\, -\frac{1}{B} \, \langle \bvperp, \dD_{\bfx}\eDpar\,\bvperp\rangle\,\left\langle
 \bJ\bvperp\,,\,\d_t\eDpar+2\vpar\,\dD_{\bfx}\eDpar\,\eDpar\right\rangle.
$$

To summarize, gathering \eqref{e:x}, \eqref{e:vp}, and  \eqref{e:ep}, we have derived from \eqref{e:1}-\eqref{e:2} the following system of equations, 
\be
\label{e:1new}
\left\{\ba{l}
\ds\frac{\dD \GC}{\dD t}
\,=\, \vpar\eDpar+\eps \,\bU,
\\[1em]
\ds\frac{\dD}{\dD t}\left[\vpar
\,+\,\frac{\eps}{B}\,\left\langle\bJ\,\bvperp\,,\,\d_t\eDpar+\vpar\,\dD_{\bfx}\eDpar\,\eDpar\right\rangle\right] \,=\,\Epar
+\left\langle\bvperp,\dD_{\bfx}\eDpar\bvperp\right\rangle
+\eps \,u_1,\\[1em]
\ds\frac{\dD}{\dD t}\left[\eperp \,+\, \frac{\eps}{B}\,\left\langle
{\bJ}\,\bvperp,\bF_0\right\rangle
\right]\,=\,-\vpar\left\langle\bvperp,\dD_{\bfx}\eDpar\bvperp\right\rangle+\eps \,d_1\,.
\ea\right.
\ee
Now to derive the leading order of an uncoupled slow dynamics it remains to analyze the contribution of the quadratic term $\left\langle\bvperp,\dD_{\bfx}\eDpar\bvperp\right\rangle$ that appears --- at zeroth order with respect to
$\eps$ --- in the last equations of \eqref{e:1new}. 

\subsection{Elimination of quadratic terms}
\label{s:quadratic}

Lemma~\ref{1st} encodes that all terms linear in $\bvperp$ --- the variable whose angle is oscillating at frequency $1/\eps$ --- are $\eps$-small in $W^{-1,\infty}$. This is directly related to the fact that they all have zero mean with respect to the fast angle. In contrast, as in the homogeneous case, quadratic terms do produce slow contributions that are asymptotically relevant. The next result identifies what are those contributions.

To state it we introduce a notion of trace restricted to the plane orthogonal to $\eDpar$. For any $\bA\in\cL_2(\R^3,\R^p)$, $p\in\N^*$, at any point $\bfx$ and time $t$
\be
\label{Trperp}
\Trperp{t,\bfx}\,\bA\,=\,\Tr\,\bA\,-\,\bA(\eDpar(t,\bfx),\eDpar(t,\bfx)).
\ee 
In particular, for any $\bfa\in\R^3$ orthogonal to $\eDpar(t,\bfx)$,
we observe that
$$
\|\bfa\|^2\,\Trperp{t,\bfx}\,\bA
\,=\,\bA(\bfa,\bfa)\,+\,\bA(\bJ(t,\bfx)\,\bfa,\bJ(t,\bfx)\,\bfa)\,.
$$
Since with any linear operator $\bA\in\cL_1(\R^3,\R^3\otimes\R^p)$ one may
associate a quadratic operator in $\cL_2(\R^3,\R^p)$ by
$(\bfa,\bfb)\mapsto \langle \bfa,\bA\bfb\rangle$ the above definitions
may be extended to such operators by identification. We also recall that $\Re$ denotes the symmetric part.


\bl
\label{2nd}
For any $\bA\in W^{1,\infty}(\R^+_t; \cL_2(\R^3,\R^p))$, $p\in\N^*$, solutions to \eqref{e:xv} satisfy at a.e. $t$
$$
\bA(t)(\bvperp(t),\bvperp(t))
\,=\,\eperp(t)\Trperp{t,\bfx(t)}\left(\bA(t)\right)\,-\,\eps\,\frac{\dD\chi_\bA}{\dD t}(t)\,+\,\eps\,\eta_\bA(t),
$$
where 
$$\left\{\ba{l}
\ds\chi_\bA \,=\,\frac{1}{2\,B}\,\Re(\bA)\left(\bvperp,\bJ\bvperp\right)\,,
\\[0.95em]
\ds\eta_\bA\,=\,
\frac12\,\Re(\bA)(\bvperp,\bU)
\,+\,\frac{1}{2\,B}\,\Re(\bA)'\left(\bvperp,\bJ\bvperp\right)
\,+\,\frac{1}{2\,B}\,\Re(\bA)\left(\bF,\bJ\bvperp\right)\,.
\ea\right.
$$
\el

\br\label{trace-free-bis}
Consistently with Remarks~\ref{antiderivative} and~\ref{trace-free}, note that $\chi_\bA$ has itself no slow component at leading order since $\Re(\bA(t))\left(\cdot,\bJ(t,\bfx)(\cdot)\right)$ is trace-free on the plane orthogonal to $\eDpar(t,\bfx)$.
\er

\begin{proof} 
One may assume without loss of generality that $\bA$ is symmetric. Then by combining \eqref{e:2} and \eqref{e:homologic-1st} and using \eqref{e:J} we derive that
\begin{eqnarray*}
\eps\frac{\dD}{\dD t}\left[\bA\left(\bvperp,\frac{\bJ\,\bvperp}{B}\right)\right] 
&=&-\bA(\bvperp,\bvperp)
\,+\,\bA\left(\bJ\bvperp,\bJ\bvperp\right)
\\
&+&\eps\left[\bA(\bvperp,\bU)
\,+\,\frac{1}{B}\,\bA'\left(\bvperp,\bJ\bvperp\right)
\,+\,\frac{1}{B}\,\bA\left(\bF,\bJ\bvperp\right)\right]\,.
\end{eqnarray*}
By multiplying the latter by $1/2$ and adding $\bA(\bvperp,\bvperp)$ on both sides one achieves the proof. 
\end{proof}

\subsection{Proof of Theorem \ref{th:1}}
\label{sec:proofTh1}

We have now sufficient materials to prove Theorem \ref{th:1} on the
asymptotic behavior on solutions to \eqref{e:1}-\eqref{e:2} when $\eps\rightarrow 0$. 

On the one hand, applying Lemma~\ref{2nd} with the quadratic form associated with $\dD_\bfx \eDpar(t,\bfx)$, one may partially eliminate quadratic terms in $\bvperp$ from \eqref{e:1new}. As a result
\be
\label{e:voy}
\left\langle\bvperp, \dD_\bfx \eDpar\bvperp\right\rangle
=\eperp\, \Div_\bfx(\eDpar)
\,-\,\frac{\eps}{2}\frac{\dD}{\dD t}\left[\left\langle \bvperp,\frac{\Re(\dD_\bfx\eDpar)\bJ}{B}\,\bvperp\right\rangle\right] 
+\eps\, u_2,
\ee
where $u_2 =u_{21}+u_{22}+u_{23}$ with $u_{21}$,
\begin{eqnarray*}
 u_{21}
\,=\,\frac{1}{2} \,\left[ \left\langle \bvperp, \Re(\dD_\bfx
  \eDpar)\bU_{10}\right\rangle +  \left\langle
  \bF_0,\frac{\Re(\dD_\bfx \eDpar)\bJ}{B}\,\bvperp\right\rangle
  \right],
\end{eqnarray*}
whereas $u_{22}$ is given by
\begin{eqnarray*}
\ds u_{22}&=& \frac{1}{2} \,\left[ \left\langle \bvperp, \Re(\dD_\bfx \eDpar)\bU_{11}\right\rangle 
+\left\langle \bF_1,
\frac{\Re(\dD_\bfx\eDpar)\bJ}{B}\,\bvperp\right\rangle \right]
\\
&+&\frac{1}{2 \,B} \,\left\langle \bvperp,
{\Re\left(\d_t(\dD_\bfx\eDpar)+\vpar\dD_\bfx^2 \eDpar\eDpar\right)\bJ}\bvperp \right\rangle
\end{eqnarray*}
and $u_{23}$
\begin{eqnarray*}
u_{23}&=& \ds\frac{1}{2} \,\left[ \left\langle \bvperp, \Re(\dD_\bfx \eDpar)\bU_{12}\right\rangle 
+\left\langle \bF_2,
\frac{\Re(\dD_\bfx\eDpar)\bJ}{B}\,\bvperp\right\rangle \right]
\\
&+&\frac{1}{2\,B} \,\left\langle \bvperp,
{\Re\left(\dD_\bfx^2 \eDpar\bvperp\right)\bJ}\bvperp \right\rangle\,.
\end{eqnarray*}
Substituting \eqref{e:voy} in the second equation of \eqref{e:1new},
we get an equation for the corrected parallel velocity $\vGC$, defined in (\ref{vGC}), that is,
\be
\label{vp:2}
\frac{\dD\vGC }{\dD t} \,=\,\Epar\,+\,\eperp\, \Div_\bfx\eDpar+\eps\,(u_1 +u_2)\,.
\ee

On the other hand, we proceed in the same way with the quadratic term associated with
$-\vpar(t)\,\dD_{\bfx}\eDpar(t,\bfx(t))$ to transform the third
equation of \eqref{e:1new} into a new equation for a correction to the part of the kinetic energy in the perpendicular plan to the magnetic field direction $\eGC$, already defined in \eqref{eGC},
\be
\label{ep:2}
\ds\frac{\dD\eGC }{\dD t} \,=\,-\vpar \eperp\, \Div_\bfx\eDpar+\eps \,(d_1+d_2),
\ee
where $d_2=d_{21}+d_{22}+d_{23}+d_{24}$ with 
$$
\left\{
\ba{rcl}
\ds d_{21}&=&\, \ds-\vpar\,u_{21}\,,
\\\,\\
\ds d_{22}&=&\, \ds-\vpar\,u_{22}
\,-\,
\frac{\Epar}{2\,B} \, \left\langle \bvperp,\Re\left(\dD_\bfx \eDpar\right)\bJ\,\bvperp\right\rangle,
\\ \,\\
\ds d_{23}&=&\ds-\vpar\,u_{23}\,-\,  \frac{1}{2B}\,\left\langle \bvperp ,
\d_t\eDpar+\vpar\dD_{\bfx}\eDpar\eDpar \right\rangle\,\left\langle \bvperp,
\Re\left(\dD_\bfx \eDpar\right)\bJ\bvperp\right\rangle\,\,,
\ea\right.
$$
and the quartic term $d_{24}$ is
$$
\ds d_{24}\,=\, -\frac{1}{2\,B}\,  \left\langle \bvperp ,
\dD_{\bfx}\eDpar\bvperp \right\rangle\, \left\langle \bvperp,
{\Re\left(\dD_\bfx \eDpar\right)\bJ}\bvperp\right\rangle\,.
$$
With \eqref{vp:2} and \eqref{ep:2}, System~\eqref{e:1}-\eqref{e:2} yields 
\be
\label{e:1new2}
\left\{\ba{l}
\ds
\frac{\dD\GC}{\dD t}\,=\,\vpar\eDpar(t,\bfx) +\eps \,\bU(t,\bfx,\bfv),
\\[1em]
\ds\frac{\dD \vGC}{\dD t}\,=\,\Epar(t,\bfx)\,+\,\eperp\, \Div_\bfx\eDpar(t,\bfx) 
+\eps\,(u_1 +u_2)(t,\bfx,\bfv),
\\[1em]
\ds \frac{\dD \eGC}{\dD t}\,=\,
-\vpar \eperp\, \Div_\bfx\eDpar(t,\bfx)  +\eps \,(d_1+d_2) (t,\bfx,\bfv)\,.
\ea\right.
\ee
At this juncture, the leading-order part of the slow evolution system \eqref{e:1new2} is already uncoupled from the fast equation \eqref{e:2}. This allows to derive the following asymptotic result by mimicking the analysis of Section~\ref{s:homogeneous}, relying this time on Lemma~\ref{bnd0} to bound remainders.

\bpr
\label{1st-xve}
Under the assumptions of Theorem \ref{th:1}, there exists a constant $C>0$, depending
polynomially on $\|\bE\|_{W^{1,\infty}}$,  $\|B^{-1}\|_{W^{1,\infty}}$ and $\|\eDpar\|_{W^{2,\infty}}$, such that the following holds. Consider $(\bfx^\eps,\bfv^\eps)$ a solution to~\eqref{e:xv} starting from $(\bfx_0,\bfv_0)$. Then,  for a.e. $t\geq0 $
\begin{align*}
\|\bfx^\eps(t)-\bfy(t)\|
&+\|\vpar^\eps(t)-v(t)\|
+\|\eperp^\eps(t)-\e(t)\|\\
&\leq C\,\eps\,e^{C\,t\,(\|\bfv_0\|^3+t^3)}
\,\|\bfv_0\|\,(1+\|\bfv_0\|^2),
\end{align*}
where $(\bfy,v,e)$ solves
$$
\left\{\ba{l}
\ds\frac{\dD\bfy}{\dD t} \,=\,v\,\eDpar(t,\bfy),\\[1em]
\ds\frac{\dD v}{\dD t}\,=\,\Epar(t,\bfy)\,+\,\e\, \Div_\bfx\eDpar(t,\bfy),\\[1em]
\ds\frac{\dD \e}{\dD t}\,=\,
-v\,\e\, \Div_\bfx\eDpar(t,\bfy)\,,
\ea\right.
$$
with initial data $\bfy(0)=\bfx_0$, $v(0)={\vpar}_{0}$ and $\e(0)=\tfrac12\|{\bvperp}_0\|^2$.
\epr

Finally to achieve the proof of Theorem \ref{th:1}, we simply apply
Proposition \ref{p:ODEtoPDE} where the slow map $\cA(t,\cdot)$  is
given by $(\bfx,\bfv)\mapsto (\bfx, \vpar(t,\bfx,\bfv),\eperp(t,\bfx,\bfv))$ and the weights
$\cM$ are given by
$$
\cM(t,\bfx,\bfv)\,=\,  C\,\eps\,e^{C\,t\,(\|\bfv\|^3+t^3)}
\,\|\bfv\|\,(1+\|\bfv\|^2)
$$
with $C$ as in Proposition~\ref{1st-xve}.

\subsection{Elimination of higher-order  terms}
Though the latter result does provide some insights, in general it
fails to capture leading-order dynamics of all slow variables, since
some of them are slower than what can be described with a system of
zeroth order in $\eps$. A simple example is the essentially
two-dimensional case where $\eDpar$ is constant and asymptotically at
zeroth order only $\vpar$ and the parallel component of $\bfy$ are
moving.

To provide a more comprehensive picture, we need a system containing
terms of order $\eps$. The purpose of Theorem
\ref{th:2} is to take into account this correction. Lemmas~\ref{1st} and~\ref{2nd} already contains the basis to clean $\eps$-terms of \eqref{e:1new2} that are of order at most $2$ with respect to $\bvperp$. Yet, $d_1$, $d_2$ and $u_2$ contain cubic terms and $u_2$ also exhibits a quartic term. Therefore we need to investigate how to handle those.

We first show how to eliminate cubic terms.
\bl
\label{3rd}
Let $\bE\in W^{1,\infty}$ and $\bB$ be such that $1/B\in W^{1,\infty}$
and $\eDpar\in W^{1,\infty}$. There exists a constant $C$ depending
polynomially on $\|\bE\|_{W^{1,\infty}}$,  $\|B^{-1}\|_{W^{1,\infty}}$ and
$\|\eDpar\|_{W^{1,\infty}}$ such that for any $\bA\in
W^{1,\infty}\left(\R^+_t;\,\cL_3(\R^3,\R^p)\right)$, $p\in\N^*$,
solutions to \eqref{e:xv} satisfy at a.e. $t\geq 0$,
$$
\bA(t)(\bvperp(t),\bvperp(t),\bvperp(t))
\,=\,-\eps\frac{\dD\chi_\bA}{\dD t}(t)\,+\,\eps\,\eta_\bA(t),
$$
for some $(\chi_\bA,\eta_\bA)$ such that for a.e. $t$
$$
\left\{\ba{l}
\|\chi_\bA(t)\|\,\leq\,C\,\|\bA(t)\|\,\|\bvperp(t)\|^3\,,
\\[0.95em]
\|\eta_\bA(t)\|\,\leq\,C\,\|\bvperp(t)\|^2\,\left[\,\|\bA'(t)\|\,\|\bvperp(t)\|
+\|\bA(t)\|\,(1+\|\bfv(t)\|^2)\,\right]\,.
\ea\right.
$$
\el
\begin{proof}
One may assume without loss of generality that $\bA$ is symmetric-valued. Then from \eqref{e:2} and \eqref{e:homologic-1st} stem
\begin{eqnarray*}
\eps\,\frac{\dD}{\dD  t}\left[\bA\left(\frac{\bJ\,\bvperp}{B},\bvperp,\bvperp\right)\right](t)
&=&-\,\bA(\bvperp,\bvperp,\bvperp)+2\bA(\bJ\bvperp,\bJ\bvperp,\bvperp)
\\ &+&\,\frac{\eps}{B}\,\bA'\left(\bJ\,\bvperp,\bvperp,\bvperp\right)
+\eps\,\bA\left(\bU,\bvperp,\bvperp\right)\\
&+&\,\frac{2\eps}{B}\,\bA\left(\bJ\,\bvperp,\bF,\bvperp\right)
\end{eqnarray*}
and
\begin{eqnarray*}
\eps\,\frac{\dD}{\dD t}\left[ B^2\ds
\bA\left(\frac{\bJ\,\bvperp}{B},\frac{\bJ\,\bvperp}{B},\frac{\bJ\,\bvperp}{B}\right)\right](t)
  &=&\,-3\,\bA(\bJ\bvperp,\bJ\bvperp,\bvperp)
\\ 
&+&\,\frac{\eps}{B}\,\bA'\left(\bJ\,\bvperp,\bJ\,\bvperp,\bJ\,\bvperp\right)
+\,3\,\eps\,\bA\left(\bU,\bJ\,\bvperp,\bJ\,\bvperp\right)\\
&+&\,\frac{2\,\eps}{B}\,\frac{\d_tB+\dD_\bfx B\,\bfv}{B}\,\,\bA\left(\bJ\,\bvperp,\bJ\,\bvperp,\bJ\,\bvperp\right)\,.
\end{eqnarray*}
Then summing the former with $2/3$ of the latter yields the result.
\end{proof}

To complete the uncoupling at order $\eps$ remains the task of analyzing the possible contribution of quartic terms. By using \eqref{e:2}, \eqref{e:homologic-1st} and the fact that $\bJ(t,\bfx)^2\bfa=-\bfa$ for any $\bfa$ orthogonal to $\eDpar(t,\bfx)$, it is possible to achieve this task at the level of generality considered so far. As a result one would prove that in general the elimination of quartic terms may indeed leave relevant slow terms. However, for concision's sake we choose to specialize the discussion to the specific form required by 
$$
d_{24}\,=\,\left\langle \bvperp ,
\Re\left(\dD_{\bfx}\eDpar\right)\bvperp \right\rangle\,\left\langle \bvperp,
\Re\left(\dD_\bfx \eDpar\right)\bJ\bvperp\right\rangle
$$
and that may be eliminated at leading-order.

\bl
\label{4th-specific}
Let $\bE\in W^{1,\infty}$ and $\bB$ be such that $1/B\in W^{1,\infty}$ and $\eDpar\in W^{1,\infty}$. There exists a constant $C$ depending polynomially on $\|\bE\|_{W^{1,\infty}}$,  $\|B^{-1}\|_{W^{1,\infty}}$ and $\|\eDpar\|_{W^{1,\infty}}$ such that for any symmetric-valued $\bA\in W^{1,\infty}\left(\R^+_t;\,\cL_2(\R^3,\R^p)\right)$, $p\in\N^*$, solutions to \eqref{e:xv} satisfy at a.e. $t$
$$
\bA(t)(\bvperp(t),\bvperp(t))
\ \times\ \bA(t)(\bvperp(t),\bJ(t,\bfx(t))\bvperp(t))
\,=\,-\eps\frac{\dD\chi_{\bA,2}}{\dD t}(t)+\eps\eta_{\bA,2}(t),
$$
for some $(\chi_{\bA,2},\eta_{\bA,2})$ such that for a.e. $t$
$$
\left\{\ba{l}
\ds\|\chi_{\bA,2}(t)\| \,\leq\,C\,\|\bA(t)\|^2\,\|\bvperp(t)\|^4\,,
\\[0.95em]
\|\eta_{\bA,2}(t)\| \,\leq\,C\,\|\bA(t)\|\,\|\bvperp(t)\|^3\left(\|\bA'(t)\|\,\|\bvperp(t)\|
+\|\bA(t)\|\,(1+\|\bfv(t)\|^2)\right)\,.
\ea\right.
$$
\el
\begin{proof}
We introduce $\tbA(t):\ (\bfa,\bfb)\mapsto \bA(t)(\bfa,\bJ(t,\bfx(t))\bfb)$ and recall from Remark~\ref{trace-free-bis} that $\Trperp{t,\bfx(t)}\tbA(t)=0$. Thus revisiting the proof of Lemma~\ref{2nd} without symmetrization yields
\begin{align*}
\tbA(\bvperp,\bvperp)
&=\,-\eps\frac{\dD}{\dD t}\left(\frac{1}{2\,B}\,\tbA\left(\bvperp,\bJ\bvperp\right)\right)+\eps\,\teta_{\tbA}
\end{align*}
with 
$$
\teta_{\tbA}\,=\,\frac12\,\tbA(\bvperp,\bU)
\,+\,\frac{1}{2\,B}\,\tbA'\left(\bvperp,\bJ\bvperp\right)
\,+\,\frac{1}{2\,B}\,\tbA\left(\bF,\bJ\bvperp\right)\,.
$$
By multiplying with $\bA(\bvperp,\bvperp)$, one derives the result with
$$
\left\{\ba{l}
\ds\chi_{\bA,2}\,=\,-\frac{1}{4\,B}\,(\bA\left(\bvperp,\bvperp\right))^2\,,
\\[0.95em]
\ds\eta_{\bA,2}\,=\,\bA\left(\bvperp,\bvperp\right)\teta_{\tbA}
-\frac{1}{4\,B}\frac{\d_tB+\dD_\bfx B\,\bfv}{B}\,(\bA\left(\bvperp,\bvperp\right))^2\,.
\ea\right.
$$
\end{proof}

\subsection{Proof of Theorem \ref{th:2}}
\label{sec:proofTh2}

To spare some pieces of notation, in the justification of the foregoing claim we shall use $\lesssim$ to denote $\leq C\times$ with $C$ a local variable depending only and polynomially on $\|\bE\|_{W^{2,\infty}}$,  $\|B^{-1}\|_{W^{2,\infty}}$ and $\|\eDpar\|_{W^{3,\infty}}$. All along we consider $(\bfx(t),\bfv(t))$ a solution to \eqref{e:1}-\eqref{e:2}. We observe that
$$
\left\|\frac{\dD\bfx}{\dD t}\right\|\,\lesssim\,\|\bfv\|\,, \qquad
\left|\frac{\dD\vpar}{\dD t}\right| \,\lesssim\,1+\|\bvperp\|\,\|\bfv\|
$$
and
$$
\|\bF\|\,\lesssim\,1+\|\bfv\|^2\,, \qquad
\|\bU\| \,\lesssim\,1+\|\bfv\|^2\,.
$$

First we apply Lemma~\ref{1st} with the linear application $\bvperp\mapsto\bfU_{11}(t,\bfx(t),\vpar(t), \bvperp)$ and Lemma~\ref{2nd} with the quadratic function $\bvperp\mapsto\bfU_{12}(t,\bfx(t),\vpar(t), \bvperp)$. As a result there exist functions $\bfchi_x^1$, $\bfchi_x^2$, $\bfeta_x^1$, $\bfeta_x^2$ such that
\begin{align*}
\bfU_{11}&\ds\,=\, -\eps\frac{\dD
  \bfchi_\bfx^1}{\dD t} + \eps\,\bfeta_\bfx^1\\
\bfU_{12}&\ds\,=\,
\eperp\,\Trperp{t,\bfx(t)}\left(\dD_\bfx\left(\frac{\bJ}{B}\right)\right) \,-\,\eps\,\frac{\dD
\bfchi_\bfx^2}{\dD t} \,+\, \eps\,\bfeta_\bfx^2
\end{align*}
so that, with $\bfchi_\bfx=\bfchi_\bfx^1+\bfchi_\bfx^2$ and $\bfeta_\bfx =\bfeta_\bfx^1
+\bfeta_\bfx^2$ 
$$
\bU\,=\,\bU_{10}+\eperp\,\Trperp{t,\bfx(t)}\left(\dD_\bfx\left(\frac{\bJ}{B}\right)\right)
-\eps\frac{\dD\bfchi_\bfx}{\dD t}\,+\,\eps\bfeta_\bfx
$$
with 
\be\label{x-remainders}
\|\bfchi_\bfx\|\,\lesssim\, \|\bvperp\|\,\|\bfv\|\,,\qquad
\|\bfeta_\bfx\| \,\lesssim\, \|\bfv\|\,(1+\|\bfv\|^2)\,.
\ee

Before going on we make $\Trperp{t,\bfx(t)}\left(\dD_\bfx\left(\frac{\bJ}{B}\right)\right)$
more explicit. First, by differentiating $\frac{\bJ}{B} \eDpar=0$ we get 
$$
0 = \dD_\bfx\left[\left(\frac{\bJ}{B}\right)\eDpar\right]\,\eDpar \,=\, \left[\dD_\bfx\left(\frac{\bJ}{B}\right)\eDpar\right]\,\eDpar  + \frac{\bJ}{B}\dD_\bfx \eDpar\,\eDpar,
$$
thus, with $\curvB$ defined as in \eqref{drift-0},
$$
\left[\dD_\bfx\left(\frac{\bJ}{B}\right)\eDpar\right]\,\eDpar  \,=\,
-\frac{\bJ}{B}\dD_\bfx \eDpar\,\eDpar \,=\, \curvB\,.
$$
Therefore, recalling definition \eqref{Trperp} and Lemma~\ref{lem:00},
$$
\Trperp{t,\bfx(t)}\left(\dD_\bfx\left(\frac{\bJ}{B}\right) \right)\,=\, \Div_\bfx\left(\frac{\bJ}{B}\right)
-\curvB\,=\, \rotB + \gradB.
$$
In particular the equation on $\GC$ takes the form
\begin{eqnarray*}
\frac{\dD}{\dD t}\left[\GC+\eps^2\bfchi_x\right]
 \,=\, \vpar\eDpar(t,\bfx)\,+\,\eps^2\,\bfeta_x \,+\,\eps\,\Ud(t,\bfx,\vpar,\eperp),
\end{eqnarray*}
where $\Ud$ is defined in \eqref{drift-3} and $(\bfchi_x,\bfeta_x)$ satisfies \eqref{x-remainders}.

Likewise we may clean up the second equation of \eqref{e:1new2}. After some calculations, with arguments identical to those used hereabove, we obtain
\begin{eqnarray*}
u_1 &=& \left\langle\bE\,,\Sig\,\right\rangle-\eps\frac{\dD\chi_\mypar^1}{\dD t}+\eps\eta_\mypar^1\\
&+& \eperp \,\left[ \left\langle\Div_\bfx\left(\frac{\bJ}{B}\right),
\d_t\eDpar+\vpar\dD_\bfx\eDpar\eDpar\right\rangle \,+\,\left\langle\frac{\bJ}{B}\dD_\bfx\eDpar\eDpar,\d_t\eDpar\right\rangle  \right]
\\
&-&\eperp\,\Tr\left(\frac{\bJ}{B}\left(\d_t(\dD_\bfx\eDpar)+\vpar\dD_{\bfx}^2\eDpar\,\eDpar+\vpar\left(\dD_{\bfx}\eDpar\right)^2\right)\right)
\\
&-& \eperp\,\Tr\left[\frac{\bJ}{B}\left(\dD_\bfx\eDpar\,\eDpar\,(\d_t\eDpar)^*\right)\right].
\end{eqnarray*}
and
$$
u_2 \,=\,-\eps\frac{\dD\chi^2_\mypar}{\dD t}\,+\,\eps\,\eta^2_\mypar
$$
with $\chi_\mypar = \chi_\mypar^1+\chi_\mypar^2$ and $\eta_\mypar=\eta_\mypar^1+\eta_\mypar^2$ satisfying
\be\label{par-remainders}
|\chi_\mypar|\lesssim \|\bvperp\|\,(1+\|\bfv\|^2)\,,\qquad\qquad
|\eta_\mypar|\lesssim 1+\|\bfv\|^4\,.
\ee
In the computations aforementioned we stress that we have made extensive use of relations $\bJ\eDpar=0$, $\left(\dD_\bfx\bJ\eDpar\right)\,\eDpar=-\bJ (\dD_\bfx\eDpar)\,
\eDpar$, $\langle \bfa, \bJ\,\bfa\rangle=0$ for
any $\bfa\in\R^3$, and $(\dD_\bfx\eDpar)^*\eDpar=0$. In particular we point out that, by the skew-symmetry of values of $\dD_{\bfx}\!\bJ\ \eDpar$ and $\d_t\bJ$
\begin{eqnarray*}
\Trperp{}\left[\Re(\dD_{\bfx}\eDpar)\dD_{\bfx}\bJ\,\eDpar\right]
&=&\Tr\left[\Re(\dD_{\bfx}\eDpar)\dD_{\bfx}\bJ\eDpar\right]
-\frac12\langle\dD_{\bfx}\eDpar\eDpar,[\dD_{\bfx}\bJ\eDpar]\eDpar\rangle
\\
&=&\frac12\langle\dD_{\bfx}\eDpar\eDpar,\bJ\dD_{\bfx}\eDpar\eDpar\rangle
\,=\,0
\end{eqnarray*}
and 
\begin{eqnarray*}
\Trperp{}\left[\Re(\dD_{\bfx}\eDpar)\d_t\bJ\right]
&=&\Tr\left[\Re(\dD_{\bfx}\eDpar)\d_t\bJ\right]
-\frac12\langle\dD_{\bfx}\eDpar\eDpar,\d_t(\bJ)\eDpar\rangle\\
&=&\frac12\langle\dD_{\bfx}\eDpar\eDpar,\bJ\d_t\eDpar\rangle\,.
\end{eqnarray*}

To simplify further the expression of $u_1$ we observe that 
$$
\Tr\left[\frac{\bJ}{B}\left(\dD_\bfx\eDpar\,\eDpar\,(\d_t\eDpar)^*\right)\right]
\,=\, \left\langle\frac{\bJ}{B}\dD_\bfx\eDpar\eDpar,\d_t\eDpar\right\rangle
$$
and, by the skew-symmetry of values of $\bJ$, that
\begin{align*}
\left\langle\Div_\bfx\left(\frac{\bJ}{B}\right),
\d_t\eDpar\right\rangle
-\Tr\left(\frac{\bJ}{B}\d_t(\dD_\bfx\eDpar)\right)
&=\,-\Div_\bfx\left(\frac{\bJ}{B}\d_t\eDpar\right)\,=\,\Div_\bfx\dtB \\
\left\langle\Div_\bfx\left(\frac{\bJ}{B}\right),
\dD_\bfx\eDpar\eDpar\right\rangle
-\Tr\left(\frac{\bJ}{B}\left(\dD_{\bfx}^2\eDpar\,\eDpar+\left(\dD_{\bfx}\eDpar\right)^2\right)\right)
&=-\Div_\bfx\left(\frac{\bJ}{B}\dD_\bfx\eDpar\,\eDpar\right)\,=\,\Div_\bfx\curvB\,.
\end{align*}
As a result
$$
u_1 \,=\, -\eps\,\frac{\dD\chi_\mypar^1}{\dD t}\,+\,\eps\,\eta_\mypar^1  \,+\,
\left\langle\Sig\,,\,\bE\right\rangle +\eperp\Div_\bfx\Sig
$$
with $(\chi_\mypar^1,\eta_\mypar^1)$ as before.

Therefore gathering the expressions for $u_1$ and $u_2$, we derive
$$
\frac{\dD }{\dD t}\left[\vGC\,+\,\eps^2\chi_\mypar\right]
  \,=\,  \Epar+\eperp\,\Div_\bfx\eDpar \,+\,\eps \,\left[\left\langle\Sig\,,\,\bE\right\rangle 
\, +\,  \eperp\,\Div_\bfx\Sig\right] \,+\, \eps^2\,\eta_\mypar
$$
with $(\chi_\mypar,\eta_\mypar)$ satisfying \eqref{par-remainders}.

Finally we treat the last equation of \eqref{e:1new2} in the same manner. This leads to
\begin{eqnarray*}
d_1 &=& -\eps\frac{\dD\chi_\perp^1}{\dD t}+\eps\eta_\perp^1 -
  \vpar\eperp\Div_\bfx\Sig\\
&+&\eperp\,\left[\left\langle \rotB + \gradB\,,\, \bE
    \right\rangle -
    \Tr\left(\frac{\bJ}{B}\dD_{\bfx}\bE\right)\right],\\
d_2&=&-\eps\frac{\dD\chi_\perp^2}{\dD t}\,+\,\eps\,\eta_\perp^2,
\end{eqnarray*}
with $\chi_\perp=\chi^1_\perp+\chi^2_\perp$ and $\eta_\perp=\eta_\perp^1+\eta_\perp^2$ satisfying
\be\label{perp-remainders}
|\chi_\perp|\lesssim \|\bvperp\|\,(1+\|\bfv\|^3)\,,\qquad\qquad
|\eta_\perp|\lesssim 1+\|\bfv\|^5\,.
\ee
Now to simplify the expression for $d_1$ we observe that from Lemma~\ref{lem:00} stems
$$
\left\langle \rotB+\curvB + \gradB\,,\, \bE\right\rangle
\,-\,\Tr\left(\frac{\bJ}{B}\dD_{\bfx}\bE\right) \,=\,
-\Div_\bfx\left( \frac{\bJ \,\bE}{B}\right) \,=\, -\Div_\bfx \EcB, 
$$
so that
$$
d_1 \,=\, -\eps\frac{\dD\chi_\perp^1}{\dD t}+\eps\eta_\perp^1 \,-\,\eperp\,\left[\vpar\,\Div_\bfx\Sig\,+\,\left\langle \curvB\,,\, \bE\right\rangle  \,+\,
    \Div_\bfx\EcB  \right].
$$
The upshot is
$$
\frac{\dD }{\dD t}\left[\eGC+\eps^2\,\chi_\perp\right]
  \,=\,  -\vpar \eperp\Div_\bfx\eDpar  
-\eps\eperp\left[\Div_\bfx \left(\EcB+ \vpar\Sig\right)+\left\langle \curvB,\bE\right\rangle  \right] +\eps^2\,\eta_\perp
$$
with $(\chi_\perp,\eta_\perp)$ satisfying \eqref{perp-remainders}.

Altogether we have derived
\be
\label{e:1new3}
\left\{\ba{ll}
\ds\frac{\dD}{\dD t}\left[\GC+\eps^2\bfchi_\bfx\right]
& \ds =\, \vpar\eDpar(t,\bfx)\,+\,\eps\,\Ud(t,\bfx,\vpar,\eperp)\,+\,\eps^2\,\bfeta_\bfx,
\\[0.99em]
\ds\frac{\dD }{\dD t}\left[\vGC\,+\,\eps^2\chi_\mypar\right]
  & \ds =\,  \Epar(t,\bfx)\,+\,\eperp\,\Div_\bfx\eDpar(t,\bfx) \,+\,\eps^2\,\eta_\mypar
\\[0.95em]
\, & \ds +\,\eps\, \left(\,\left\langle\Sig,\bE\right\rangle \,+\,\eperp\,\Div_\bfx\Sig\,\right)(t,\bfx,\vpar)\,,
\\[0.99em]
\ds\frac{\dD }{\dD t}\left[\eGC\,+\,\eps^2\chi_\perp\right] & \ds
=\,-\vpar\,\eperp\, \Div_\bfx\eDpar(t,\bfx) \,+\, \eps^2\,\eta_\perp
\\[0.95em]
\, &\ds
-\,\eps\,\eperp\,\left(\, \Div_\bfx\left(\EcB+\vpar\Sig\right)+\left\langle\curvB\,,\,
  \bE\right\rangle \,\right) (t,\bfx,\vpar),
\ea\right.
\ee
with error bounds \eqref{x-remainders}-\eqref{par-remainders}-\eqref{perp-remainders}. Now we want to write \eqref{e:1new3} in terms of $(\GC,\vGC,\eGC)$ plus remainders. As in Section~\ref{s:hom-full} corrections --- of size $\eps^2$ --- stemming from terms of size
$\eps$ may be considered directly as error terms. Yet here some terms
of size $1$ are present and to deal with corrections arising from
those we follow a different path: first linearize them --- a process
that produces errors of size $\eps^2$ that can be handled directly ---
then remove the terms of size $\eps$ introduced in this way by using
Lemmas~\ref{1st} and~\ref{2nd} and the fact that $(\bfx,\vpar,\eperp)$ differs from $(\GC,\vGC,\eGC)$ by terms that are either linear in $\bvperp$ or quadratic in $\bvperp$ but trace-free in the plane orthogonal to $\eDpar(t,\bfx)$, as follows from Remark~\ref{trace-free-bis}. Besides aforementioned estimates this elimination also requires 
$$
\left|\frac{\dD\eperp}{\dD t}\right|\,\lesssim\,\|\bvperp\|\,(1+|\vpar|\,\|\bfv\|)
$$
and results in new functions $(\widehat\bfchi_\bfx,
\widehat\chi_\mypar, \widehat\chi_\perp)$, and  $(\widehat\bfeta_\bfx,
\widehat\eta_\mypar, \widehat\eta_\perp)$ such that,
\be
\label{e:1new4}
\left\{
\ba{ll}
\ds\frac{\dD}{\dD t}\left[\GC+\eps^2 \,\widehat\bfchi_\bfx\right]
& \ds =\,\vGC\,\eDpar(t,\GC)\,+\,\eps\Ud(t,\GC,\vGC,\eGC)
\,+\,\eps^2\, \widehat\bfeta_\bfx\,,
\\[0.99em]
\ds\frac{\dD }{\dD t}\left[\vGC+\eps^2\, \widehat\chi_\mypar\right]
& \ds =\,\Epar(t,\GC)\,+\,\eGC \Div_\bfx\eDpar(t,\GC) \,+\,\eps^2\, \widehat\eta_\mypar
\\[0.95em]
\,&\ds +\,\eps\, \left(\,\left\langle\Sig,\bE\right\rangle \,+\,\eGC\Div_\bfx\Sig\,\right)(t,\GC,\vGC)\,,
\\[0.99em] 
\ds\frac{\dD }{\dD t}\left[\eGC+\eps^2 \,\widehat\chi_\perp\right]
& \ds =\,-\vGC\,\eGC\, \Div_\bfx\eDpar(t,\GC)  \,+\, \eps^2\,
\widehat\eta_\perp
\\[0.95em]
\, &\ds -\,\eps\,\eGC\,\left(\, \Div_\bfx\left(\EcB
    +\vGC\Sig\right) + \left\langle\curvB\,,\,
    \bE\right\rangle\,\right)(t,\GC,\vGC),
\ea\right.
\ee
with 
$$
\|\widehat\bfchi_\bfx\|\,\lesssim\,\|\bvperp\|\,\|\bfv\|\,,\qquad
|\widehat\chi_\mypar|\,\lesssim\,\|\bvperp\|\,(1+\|\bfv\|^2)\,,\qquad
|\widehat\chi_\perp|\,\lesssim\,\|\bvperp\|\,(1+\|\bfv\|^3)
$$
and
$$
\|\widehat\bfeta_\bfx\|\,\lesssim\,\|\bfv\|\,(1+\|\bfv\|^2)\,,
\qquad
|\widehat\eta_\mypar|\,\lesssim\,1+\|\bfv\|^4\,,
\qquad
|\widehat\eta_\perp|\,\lesssim\,1+\|\bfv\|^5\,.
$$

At this stage arguing as in Section~\ref{s:hom-full} we prove the following

\bpr
\label{prop:4.9}
Under the assumptions of Theorem \ref{th:2}, there exists a constant $C>0$, depending
polynomially on $\|\bE\|_{W^{2,\infty}}$,  $\|B^{-1}\|_{W^{2,\infty}}$ and $\|\eDpar\|_{W^{3,\infty}}$ such that the following holds. Let
$(\bfx^\eps,\bfv^\eps)$ be the solution to~\eqref{e:xv} starting from
$(\bfx_0,\bfv_0)$ and $\ZGC^\eps=(\GC^\eps,\vGC^\eps,\eGC^\eps)$ be deduced from it through 
\eqref{GC}-\eqref{eGC}. Then, for a.e. $t\geq0$
\begin{eqnarray*}
\left\|\ZGC^\eps(t)-\bfZ^\eps(t)\right\| \;\leq\, C\,\eps^2\,e^{C\,t\,(\|\bfv_0\|^3+t^3)(1+\eps\,(\|\bfv_0\|+t))}
\,\|\bfv_0\|\,(1+\|\bfv_0\|^3),
\end{eqnarray*}
where $\bfZ^\eps=(\bfy^\eps,v^\eps,\e^\eps)$ solves
$$
\left\{
\ba{l}
\ds\frac{\dD\bfy^\eps}{\dD t}\; =\,v^\eps\,\eDpar(t,\bfy^\eps)\,+\,\eps\,\Ud(t,\bfy^\eps,v^\eps,\e^\eps),
\\[0.99em]
\ds\frac{\dD v^\eps}{\dD t}\; =\,\Epar(t,\bfy^\eps)\,+\,\e^\eps\, \Div_\bfx\eDpar(t,\bfy^\eps)\;+\,\eps\, \left(\,\left\langle\Sig,\bE\right\rangle+\e^\eps\Div_\bfx\Sig\,\right)(t,\bfy^\eps,v^\eps),
\\[0.99em]
\ds\frac{\dD \e^\eps}{\dD t} \, =\,-v^\eps\,\e^\eps\, \Div_\bfx\eDpar(t,\bfy^\eps) -\,\eps\,\e^\eps\,\left(\, \Div_\bfx\left(\EcB+v^\eps\Sig\right)+\left\langle\curvB,\bE\right\rangle \,\right) (t,\bfy^\eps,v^\eps),
\ea\right.
$$
with  $\bfZ^\eps(0)=\ZGC^\eps(0)$.
\epr

We may then use Proposition~\ref{p:ODEtoPDE} to derive Theorem~\ref{th:2} from Proposition~\ref{prop:4.9}.

\section{A toroidal axi-symmetric geometry: proof of Theorem~\ref{th:3}}\label{s:axi}

We now want to provide a three-dimensional analogous to Section~\ref{s:toy-long}, that is, a description of a long-time slow dynamics. Yet the presence of terms of order $1$ in the (short) time asymptotics prevents this from happening unless those terms generates a confined purely oscillatory dynamics in some components and one focuses on the remaining ones.

This requires a special form of geometry of magnetic field lines. We introduce now an example of such a configuration. 

\subsection{Geometric framework}
\label{s:geom-framework}

Let us fix a unitary vector $\eDz\in\bS^2$ to define an axis of symmetry. For vectors $\bfx\notin \R\eDz$, it is expedient to introduce 
$$
z(\bfx)=\langle\bfx,\eDz)\,,\qquad
r(\bfx)=\|\eDz\wedge\bfx\|\,,\qquad
\eDr(\bfx)=\frac{\bfx-z(\bfx)\eDz}{r(\bfx)}\,.
$$
Note that then by construction $\eDr(\bfx)$ is unitary, orthogonal to $\eDz$ and
$$
\bfx\,=\,r(\bfx)\,\eDr(\bfx)+z(\bfx)\,\eDz\,.
$$

We now assume that far from the axis $\R\eDz$ the magnetic field is stationary, toroidal, axi-symmetric and non vanishing, that is (up to a change of $\eDz$ with $-\eDz$), for some $r_0>0$, when $r(\bfx)\geq r_0$ 
$$
\eDpar(\bfx)\,=\,\frac{\eDz\wedge\bfx}{r(\bfx)}\qquad\qquad\textrm{and}\qquad\qquad
B(\bfx)\,=\,b(r(\bfx),z(\bfx))
$$
for some function $b$ with $1/b\in L^\infty([r_0,+\infty[\,\times\R)$. Note that the first equality already ensures $\Div_\bfx(\eDpar)(\bfx)=0$ when $r(\bfx)\geq r_0$ so that the second one is actually equivalent to the natural condition $\Div_\bfx\bB\equiv0$. 

In this context straightforward computations yield when $r(\bfx)\geq r_0$ 
\begin{align*}
\dD_\bfx r(\bfx)&=\langle\eDr(\bfx),\,\cdot\,\rangle\,,&
\dD_\bfx \eDr(\bfx)&=\frac{\eDpar(\bfx)}{r(\bfx)}\,\langle\eDpar(\bfx),\,\cdot\,\rangle\,,&
\dD_\bfx \eDpar(\bfx)&=-\frac{\eDr(\bfx)}{r(\bfx)}\,\langle\eDpar(\bfx),\,\cdot\,\rangle
\end{align*}
and for any $\bfa\in\R^3$,
$$
\dD_\bfx\bJ(\bfx)\,\bfa
\,=\,\frac{\langle\eDpar(\bfx),\bfa\rangle}{r(\bfx)}
\left(\eDz\,\langle\eDpar(\bfx),\,\cdot\,\rangle
-\eDpar(\bfx)\,\langle\eDz,\,\cdot\,\rangle\right)
$$
and 
$$
\dD_\bfx^2\eDpar(\bfx)(\eDpar(\bfx),\,\cdot\,)
\,=\,\frac{\eDr(\bfx)}{(r(\bfx))^2}\,\langle\eDr(\bfx),\,\cdot\,\rangle
-\frac{\eDpar(\bfx)}{(r(\bfx))^2}\,\langle\eDpar(\bfx),\,\cdot\,\rangle\,,
$$
so that in particular
$$
\Div_{\bfx}\left(\frac{\bJ}{B}\right)(\bfx) \,=\,\d_r\left(\frac1b\right)(r(\bfx),z(\bfx))\,\eDz
-\d_z\left(\frac1b\right)(r(\bfx),z(\bfx))\,\eDr(\bfx)
$$
and
$$
\Div_\bfx\Sig(\bfx)
\,=\,\frac{1}{r(\bfx)}\,\d_z\left(\frac1b\right)(r(\bfx),z(\bfx))\,,
$$
whereas the drifts $\bfF$ and $\bfU$ are given by
$$
\bF(t,\bfx,\bfv) \,=\, \Eperp(t,\bfx)
\,+\,\frac{\vpar}{r(\bfx)}\left\langle \eDr(\bfx),\bvperp\right\rangle\,\eDpar(t,\bfx)
\,+\,\frac{\vpar^2}{r(\bfx)}\,\eDr(\bfx)
$$
and
\begin{eqnarray*}
\bU(t,\bfx,\bfv)
&=& \EcB(t,\bfx)
+\frac{\vpar}{r(\bfx)\,b(r(\bfx),z(\bfx))}\,\left(\vpar\eDz
-
\,\langle\eDz,\bvperp\rangle\,\eDpar(\bfx)\right)
\\
&+&\left(\d_r\left(\frac1b\right)(r(\bfx),z(\bfx))\,\langle\eDr(\bfx),\bvperp\rangle
+\d_z\left(\frac1b\right)(r(\bfx),z(\bfx))\,\langle\eDz,\bvperp\rangle\right)\bJ(\bfx)\bvperp\,.
\end{eqnarray*}

We assume moreover that $\Epar(t,\bfx)\equiv 0$ when $r(\bfx)\geq r_0$ and that $\bE$ is axi-symmetric, hence
$$
\Eperp(t,\bfx)\,=\,\Er(t,r(\bfx),z(\bfx))\,\eDr(\bfx)\,+\,\Ez(t,r(\bfx),z(\bfx))\,\eDz,
$$
for some $\Er$ and $\Ez$. With these notational conventions, we may derive from System~\eqref{e:1new3} 
\be
\label{e:rzve}
\left\{\ba{rcl}
\ds\frac{\dD}{\dD t}\Big[r
&\hspace{-1.5em}-&\ds\hspace{-1.5em}\frac{\eps}{b}\,\langle\eDz,\bvperp\rangle
+\eps^2\langle\eDr,\bfchi_x\rangle\Big]\\[0.5em]
&=&\ds\hspace{-0.5em}
\eps\left(-\frac{\Ez}{b}-\eperp\,\d_z\left(\frac1b\right)\right)
+\eps^2\left(\langle\eDr,\bfeta_x\rangle+\frac{\vpar}{r}\langle\eDpar,\bfchi_x\rangle\right),\\[1em]
\ds\frac{\dD}{\dD t}\Big[z
&\hspace{-1.5em}+&\ds\hspace{-1.5em}\frac{\eps}{b}\,\langle\eDr,\bvperp\rangle
+\eps^2\langle\eDz,\bfchi_x\rangle\Big]\\[0.5em]
&=&\ds\hspace{-0.5em}
\eps\left(\frac{\Er}{b}+\frac{\vpar^2}{r\,b}
+\eperp\,\d_r\left(\frac1b\right)\right)
+\eps^2\langle\eDz,\bfeta_x\rangle,\\[1em]
\ds\frac{\dD }{\dD t}\Big[\vpar
&\hspace{-1.5em}+&\ds\hspace{-1.75em}\,\frac{\eps\,\vpar}{r\,b}
\,\left\langle\eDz,\bvperp\right\rangle
+\eps^2\chi_\mypar\Big]\\[1em]
&=&\ds
\frac{\eps\,\vpar}{r}\,\left(\frac{\Ez}{b}
+\eperp\d_z\left(\frac1b\right)\right)+\eps^2\eta_\mypar\,,\\[1em]
\ds\frac{\dD }{\dD t}\Big[\eperp
&\hspace{-1.5em}-&\ds\hspace{-1.5em}\eps\,\,
\left(\left\langle\EcB,\bvperp\right\rangle
+\frac{\vpar^2}{r\,b}
\,\left\langle\eDz,\bvperp\right\rangle
\right)
+\eps^2\chi_\perp\Big]\\[0.5em]
&=&\ds
\eps\,\eperp\,
\left(\d_r\left(\frac{\Ez}{b}\right)-\d_z\left(\frac{\Er}{b}\right)
-\frac{\vpar^2}{r}\d_z\left(\frac1b\right)\right)
+\eps^2\eta_\perp\,.
\ea\right.
\ee

\subsection{Uniform bounds and asymptotics}

As follows from the proof of the following proposition the assumption $\Epar(t,\bfx)\equiv 0$ is sufficient by itself to improve the dependence on time of uniform bounds on velocity. 

\bpr
Under the assumptions of Theorem \ref{th:2}, for any $r_1>r_0$, there exist positive $\eps_0$, $\tau_0$  and $C_0$, $(1/\eps_0,1/\tau_0,C_0)$ depending polynomially on $1/r_0$, $1/(r_1-r_0)$ and $\|(\Er,\Ez,1/b)\|_{W^{2,\infty}([r_0,\infty[\times\R)}$ such that any $(\bfx,\bfv)$ solution to~\eqref{e:xv} starting from $(\bfx_0,\bfv_0)$ with
$$
r(\bfx_0)\geq r_1\,,\qquad\qquad 
0<\eps\,\leq\,\frac{\eps_0}{1+\|\bfv_0\|}
$$
satisfies for a.e. $0\leq t\leq[\tau_0\,(1+\|\bfv_0\|^2)^{-1}]/\eps$ 
$$
r(\bfx(t))\geq r_0\qquad\textrm{and}\qquad
\|\bfv(t)\|\leq C_0\,(\|\bfv_0\|+\eps\,t)\,.
$$
\epr

\begin{proof}
We start with
$$
\frac{\dD}{\dD t}\Big[\frac12\vpar^2+\eperp\Big]
\,=\,\vpar\langle \Eperp(t,\bfx),\bvperp\rangle\,,
$$
then thanks to Lemma~\ref{1st} we derive
\begin{eqnarray*}
&&\frac{\dD}{\dD t}\left[\frac12\vpar^2+\eperp
-\eps\vpar\langle\EcB(t,\bfx),\bvperp\rangle\right]
\\
&& =\,\eps\,\vpar\left(\langle \d_t\Eperp(t,\bfx)+\dD_\bfx\Eperp(t,\bfx)\bfv,\bvperp\rangle
+\langle \Eperp(t,\bfx),\bU(t,\bfx,\bfv)\rangle\right)\,.
\end{eqnarray*}
Therefore as long as $r(\bfx)\geq r_0$
\begin{align*}
\max_{s\in[0,t]}\|\bfv(s)\|^2
&\leq
\|\bfv_0\|^2(1+C\,\eps)+C\,\eps\,\max_{s\in[0,t]}\|\bfv(s)\|^2
+C\,\eps\,t\,\max_{s\in[0,t]}\|\bfv(s)\|\,(1+\max_{s\in[0,t]}\|\bfv(s)\|^2),
\end{align*}
for some $C$ depending on $1/r_0$ and $\|(\Er,\Ez,1/b)\|_{W^{1,\infty}([r_0,\infty[\times\R)}$. Thus as long as $r(\bfx)\geq r_0$ provided 
$$
0<\eps\leq\eps_0\qquad\textrm{and}\qquad
0\leq t\leq \frac{\tau_0}{1+\|\bfv_0\|}\ \frac1\eps\,,
$$
we derive
$$
\|\bfv(t)\|\leq C_0\,(\|\bfv_0\|+\eps\,t),
$$
with $\eps_0$ and $\tau_0$ sufficiently small and $C_0$ sufficiently large
depending on $\|(\Er,\Ez,1/b)\|_{W^{1,\infty}([r_0,\infty[\times\R)}$
and $1/r_0$.  Then using this bound in an integrated version of the first equation of System~\eqref{e:rzve} achieves the proof provided we strengthen the constraint on times to
$$
0\leq t\leq \frac{\tau_0}{1+\|\bfv_0\|^2}\ \frac1\eps\,,
$$
with $\tau_0$ sufficiently small.
\end{proof}

From the foregoing Proposition and System~\eqref{e:rzve} we deduce the following proposition.

\bpr
\label{rzve}
Under the assumption of Theorem \ref{th:3}, there exist positive constants $\eps_0$, $\tau_0$ and $C_0$, $(1/\eps_0,1/\tau_0,C_0)$ depending polynomially on $1/r_0$, $1/(r_1-r_0)$ and $\|(\Er,\Ez,1/b)\|_{W^{2,\infty}([r_0,\infty[\times\R)}$, such that the following holds with
$$
\eps_{max}(\|\bfv_0\|):=\frac{\eps_0}{1+\|\bfv_0\|}\qquad\textrm{and}\qquad
T_{max}(\|\bfv_0\|):=\frac{\tau_0}{1+\|\bfv_0\|^2}\,.
$$
Let $(\bfx^\eps,\bfv^\eps)$ be a solution to~\eqref{e:xv} starting from
$(\bfx_0,\bfv_0)$ satisfying $r(\bfx_0)\geq r_1$. Then provided that 
$$0<\eps\,\leq \eps_{max}(\|\bfv_0\|)\,,$$
$\ZGC^\eps=(r(\bfx^\eps), z(\bfx^\eps), \vpar(\bfx^\eps,\bfv^\eps), \eperp(\bfv^\eps))$ satisfies for a.e. $0\leq t\leq T_{max}(\|\bfv_0\|)/\eps$
$$
\left\|\ZGC^\eps(t)-\bfZ^\eps(t)\right\|
\,\leq\, C\,\eps\,e^{C\,\eps\,t\,\|\bfv_0\|^4}
\,(1+\|\bfv_0\|^3),
$$
where $\bfZ^\eps=(r^\eps,z^\eps,v^\eps,\e^\eps)$ solves
$$
\left\{\ba{l}
\ds\frac{\dD r^\eps}{\dD t}=\eps\left(-\frac{\Ez(t,r^\eps,z^\eps)}{b(r^\eps,z^\eps)}-\e^\eps\,\d_z\left(\frac1b\right)(r^\eps,z^\eps)\right)\,,
\\[1em]
\ds\frac{\dD z^\eps}{\dD t}=
\eps\left(\frac{\Er(t,r^\eps,z^\eps)}{b(r^\eps,z^\eps)}+\frac{(v^\eps)^2}{r\,b(r^\eps,z^\eps)}
+\e^\eps\,\d_r\left(\frac1b\right)(r^\eps,z^\eps)\right)\,,
\\[1em]
\ds\frac{\dD v^\eps}{\dD t}
\,=\;
\eps\,\frac{v^\eps}{r^\eps}\,\left(\frac{\Ez(t,r^\eps,z^\eps)}{b(r^\eps,z^\eps)}
+\e^\eps\d_z\left(\frac1b\right)(r^\eps,z^\eps)\right)\,,
\\[1em]
\ds\frac{\dD \e^\eps}{\dD t}
\;=\,
\eps\,\e^\eps\,
\left(\d_r\left(\frac{\Ez}{b}\right)(t,r^\eps,z^\eps)-\d_z\left(\frac{\Er}{b}\right)(t,r^\eps,z^\eps)
-\frac{(v^\eps)^2}{r^\eps}\d_z\left(\frac1b\right)(r^\eps,z^\eps)\right)\,,
\ea\right.
$$
with $\bfZ^\eps(0)=\ZGC^\eps(0)$.
\epr

Finally from Proposition \ref{rzve}, we derive Theorem~\ref{th:3} through Proposition~\ref{p:ODEtoPDE}.

\br\label{rk:scaled-time}
As in Remark~\ref{rk:toy-scaled-time}, we stress that the proof also yields the analysis of dynamics involving fields depending on $\eps$ but satisfying bounds uniform with respect to $\eps$. In particular the result may be extended without change to the case where $\bE^\eps(t,\bfx)=\bE(\eps\,t,\bfx)$, $0<\eps\lesssim 1$. In this somehow simpler case the asymptotic dynamics is essentially independent of $\eps$ at leading order since
$$
\bfZ^\eps(t)\,=\,\bfZ(\eps\,t),
$$
with $\bfZ$ independent of $\eps$.
\er

\br\label{rk:axi}
Though we have chosen not to delve into this here as it would have lead us too far beyond our scope, one may remove the assumption that $\bE$ is axi-symmetric and still obtain a similar result provided one stays away from $\vpar=0$. It would follow from an analysis similar to the one expounded here but using instead of \eqref{e:2} an equation encoding rotation of $\bfx$ around $\eDz$ at speed $\vpar$.
\er

\section{A self-consistent case}\label{s:nl}

To illustrate that the foregoing analysis may also be carried out in some nonlinear cases we now consider 
\be\label{eq:nl-vlasov}
\left\{
\begin{array}{l}\ds
\d_t f^\eps\,+\,\Div_\bfx(f^\eps\,\bfv)
\,+\,\Div_\bfv\left(f^\eps\,\left(\frac{\bfv\wedge \bB(t,\bfx)}{\eps}\,+\,\bE^\eps(t,\bfx)\right)\right)\,=\,0\,,\\[0.5em]\ds
\bE^\eps(t,\bfx)=(\bK\star_{\bfx}(\rho^\eps(t,\cdot)-\urho(t,\cdot)))(x)\,,\qquad 
\rho^\eps(t,\bfx)=\int_{\R^3}\,f^\eps(t,\bfx,\bfv)\dD \bfv\,,
\end{array}
\right.
\ee
where $\bK$ is a fixed vector-valued kernel and $\urho$ is a fixed background density (representing possible other species). To stay focused on robust ubiquitous mechanisms and reuse as much as possible the estimates of the linear case expounded so far, we assume that $\bK$ is as smooth and localized as required by the analysis.

If we were to allow singular kernels, the foregoing system would include some of the classical Vlasov-Poisson systems. Yet this would lead us to delve into technical details related to the choice of topologies adapted to to the singularity at hand, a case-by-case study. Even for the two-dimensional Vlasov-Poisson case with a uniform magnetic field, the uniform estimates stemming from the divergence-free structure and the conservation of energy are insufficient to remain at the level of smooth solutions and lead the consideration of Di Perna-Lions solutions for \eqref{eq:nl-vlasov} and Delort solutions for the limiting system \cite{Miot-2D-gyrokinetic}. Note that the nature of the singularity depends dramatically on fine details of the modeling: in particular taking it account screening effects already tames the Poisson singularity, and even smoothed kernels play a deep intermediate role in the analysis of mean-field limits and the design and convergence analysis of particle-in-cell methods.

In the following we denote $\cM$ the space of finite Radon measures, $\cM_+$ its subspace of nonnegative finite Radon measures, $BV$ the space of functions of bounded variation, that is, of finite Radon measures with gradient\footnote{Thus Sobolev embeddings imply that those measures are actually absolutely continuous, hence may be identified with densities.} a finite Radon measure, and $BV_+$ the space of finite nonnegative Radon measures with gradient a finite Radon measure. Classical arguments prove the following proposition. 

\bpr
\label{p:WP-nl}
Assume $\bB\in W^{1,\infty}$, $\bK\in W^{1,\infty}$ and $\urho\in \cC^0(\R_+;\cM_+(\R^3)-weak*)$.\\ Then for any $\eps>0$, and any $f_0\in BV_+$, there exists a unique distribution function\\ $f^\eps\in\cC^0(\R_+;\cM_+(\R^6)-weak*)\cap L^\infty_{loc}(\R_+;BV(\R^6))$ solving \eqref{eq:nl-vlasov} starting from $f^\eps(0,\cdot,\cdot)=f_0$. Moreover the above $f^\eps$ is obtained by pushing forward $f_0$ by the characteristic flow of \eqref{eq:nl-vlasov}.
\epr

\subsection{First-order asymptotics}\label{s:1st-nl}

We want to use Propositions~\ref{1st-xve} and~\ref{prop:4.9} in the way already pointed out in Remark~\ref{rk:toy-scaled-time}, that is, with $\eps$-dependent electric fields satisfying bounds uniform with respect to $\eps$.

In this direction, our first observation is that the solution $f^\eps$ from Proposition~\ref{p:WP-nl} satisfies for any $t\geq0$, any $\eps>0$ and any $\ell\geq0$,
\begin{align*}
\int_{\R^3}\rho^\eps(t,\bfx)\,\dD \bfx
&\,=\,\int_{\R^6}f^\eps(t,\bfx,\bfv)\,\dD \bfx\,\dD \bfv
\,=\,\int_{\R^6}f_0(\bfx,\bfv)\,\dD \bfx\,\dD \bfv\,,\\
\|\bE^\eps(t,\cdot)\|_{W^{\ell,\infty}}
&\leq \|\bK\|_{W^{\ell,\infty}}\,(\int_{\R^6}f_0(\bfx,\bfv)\,\dD \bfx\,\dD \bfv+
\int_{\R^3}\urho(t,\dD \bfx))\,,
\end{align*}
with $\bE^\eps$ as in \eqref{eq:nl-vlasov}. This already ensures a uniform use of Lemma~\ref{bnd0}. Our second observation on $f^\eps$ is that since
\[
\d_t \rho^\eps\,+\,\Div_\bfx(\,\bfj^\eps)\,=\,0\,,\qquad\qquad 
\bfj^\eps(t,\bfx)=\int_{\R^3}\,\bfv\,f^\eps(t,\bfx,\bfv)\,\dD \bfv\,,
\]
we have for any $t\geq0$, any $\eps>0$ and any $\ell\geq0$,
\begin{align*}
\int_{\R^3}\|\bfj^\eps\|(t,\bfx)\,\dD \bfx
&\,\leq\,\int_{\R^6}\|\bfv\|\,f_0(\bfx,\bfv)\,\dD \bfx\,\dD \bfv
+2\,t\,\|\bE^\eps\|_{L^\infty([0,t]\times\R^3)}\,\int_{\R^6}f_0(\bfx,\bfv)\,\dD \bfx\,\dD \bfv\,,\\
\|\d_t\bE^\eps(t,\cdot)\|_{W^{\ell,\infty}}
&\leq 
\|\bK\|_{W^{\ell+1,\infty}}\,\int_{\R^3}\|\bfj^\eps\|(t,\bfx)\,\dD \bfx+
\|\bK\|_{W^{\ell,\infty}}\,\int_{\R^3}|\d_t\urho|(t,\dD \bfx)\,.
\end{align*}

This leads to a nonlinear version of Theorem~\ref{th:1}. To state it we modify notation $\cV_0=\cV_0^\bE$ introduced in \eqref{V:0} to mark the dependence of the vector-field on the electric field $\bE$.

\bt\label{th:1-nl}
Assume $\bB\in W^{1,\infty}$ is such that $1/B\in W^{1,\infty}$ and $\eDpar\in W^{2,\infty}$, $\bK\in W^{2,\infty}$ and $\urho\in W^{1,\infty}(\R_+;\cM_+(\R^3))$. For any $M>0$, there exists a constant $C$ depending polynomially on $\|\bK\|_{W^{2,\infty}}$,  $\|B^{-1}\|_{W^{1,\infty}}$, $\|\eDpar\|_{W^{2,\infty}}$, $\|\urho\|_{W^{1,\infty}(\R_+;\cM_+(\R^3))}$ and $M$ such that if $f^\eps$ solves \eqref{eq:nl-vlasov} with initial data $f_0\in BV_+$ such that
\[
\int_{\R^6}(1+\|\bfv\|)\,f_0(\bfx,\bfv)\,\dD \bfx\,\dD \bfv\,\leq\, M, 
\]
then $F^\eps$ defined by
$$
F^\eps(t,\bfx,v_\mypar,\eperp)=\int_{\bS_{t,\bfx}}\,f^\eps(t,\bfx,v_\mypar\,\eDpar(t,\bfx)\,+\,\sqrt{2\,\eperp}\ \widehat{\beD})\ \dD \sigma_{t,\bfx}(\widehat{\beD}),
$$
with $\bS_{t,\bfx}=\{\eDpar(t,\bfx)\}^\perp\cap\bS^2$ and $\sigma_{t,\bfx}$ its canonical line-measure, satisfies for any $t\geq0$
$$
\|F^\eps(t,\cdot)-G(t,\cdot)\|_{\dot{W}^{-1,1}}\leq \,\eps\,\delta_{f_0}(t)\,,
$$
where 
\begin{align*}
\delta_{f_0}(t)&=
C\,\exp\left(C\,e^{C\,t^4}\,\int_{\R^3\times\R^3} e^{C\,t\,\|\bfv\|^3}\,f_0(\bfx,\bfv)\,\dD\bfx\,\dD\bfv\right)\\
&\quad\times e^{C\,t^4}\,\int_{\R^3\times\R^3} e^{C\,t\,\|\bfv\|^3}\,\|\bfv\|\,(1+\|\bfv\|^2)\,f_0(\bfx,\bfv)\,\dD\bfx\,\dD\bfv
\end{align*}
and $G$ solves 
\be\label{eq:1st-nl}
\left\{
\begin{array}{l}\ds
\d_t G\,+\,\Div_\bfZ \left(\cV_0^\bE\,G \right) \,=\,0\,,\\[0.5em]\ds
\bE(t,\bfx)=(\bK\star_{\bfx}(\rho(t,\cdot)-\urho(t,\cdot)))(x)\,,\qquad 
\rho(t,\bfx)=\int_{\R\times\R_+}\,G(t,\bfx,v,w)\dD v\,\dD w\,,
\end{array}
\right.
\ee
with initial datum $G_0$
\be
\label{eq:1st-nlG0}
G_0(\bfZ)=\int_{\bS_{0,\bfy}}\,f_0(\bfy,v\,\eDpar(0,\bfy)\,+\,\sqrt{2\,\e}\ \widehat{\beD})\ \dD \sigma_{0,\bfy}(\widehat{\beD})\,,
\ee
where $\cV_0^\bE$ is still given by formula \eqref{V:0} but $\bE$ is now a self-consistent electric field.
\et

\begin{proof}
Applying\footnote{Actually we rather inspect the proof to check that bounds on $\d_t\bE^\eps$ are required only in consistency errors and not in Lipschitz constants so as to track the effect of the growth in time of bounds on $\d_t\bE^\eps$.} Proposition~\ref{1st-xve} already gives 
\[
\|F^\eps(t,\cdot)-G^\eps(t,\cdot)\|_{\dot{W}^{-1,1}}\leq C\,\eps\,e^{C\,t^4}\,\int_{\R^3\times\R^3} e^{C\,t\,\|\bfv\|^3}\,\|\bfv\|\,(1+\|\bfv\|^2)\,f_0(\bfx,\bfv)\,\dD\bfx\,\dD\bfv,
\]
with $G^\eps$ solving
\[
\left\{
\begin{array}{l}\ds
\d_t G^\eps\,+\,\Div_\bfZ \left(\cV_0^{\bE^\eps}\,G^\eps \right) \,=\,0\,,\\[0.5em]\ds
\bE^\eps(t,\bfx)=(\bK\star_{\bfx}(\rho^\eps(t,\cdot)-\urho(t,\cdot)))(x)\,,\qquad 
\rho^\eps(t,\bfx)=\int_{\R^3}\,f^\eps(t,\bfx,\bfv)\dD \bfv\,,
\end{array}
\right.
\]
with initial datum $G_0$ given by \eqref{eq:1st-nlG0}. It is thus sufficient to compare $G^\eps$
with $G$ the unique solution to \eqref{eq:1st-nl}.

In this direction, we first observe that System~\eqref{eq:1st-nl} support direct counterparts to Proposition~\ref{p:WP-nl} and Lemma~\ref{bnd0} so that $\bE$ satisfies exactly the same bounds as the ones derived above for $\bE^\eps$. A direct comparison of the respective characteristics for $\cV_0^{\bE^\eps}$ and $\cV_0^{\bE}$ show that with a constant $C$ as in Theorem~\ref{th:1-nl}, for any $t\geq0$
\[
\|G(t,\cdot)-G^\eps(t,\cdot)\|_{\dot{W}^{-1,1}}\leq C\,e^{C\,t^4}\,\int_{\R^3\times\R^3} e^{C\,t\,\|\bfv\|^3}\,f_0(\bfx,\bfv)\,\dD\bfx\,\dD\bfv
\ \times\ \int_0^t\|\bE^\eps(s,\cdot)-\bE(s,\cdot)\|_{L^\infty}\dD s
\]
so that for any $t\geq0$
\begin{align*}
&\|G(t,\cdot)-F^\eps(t,\cdot)\|_{\dot{W}^{-1,1}}\\
&\leq \|F^\eps(t,\cdot)-G^\eps(t,\cdot)\|_{\dot{W}^{-1,1}}\\
&+C\,\|\bK\|_{W^{1,\infty}}\,e^{C\,t^4}\,\int_{\R^3\times\R^3} e^{C\,t\,\|\bfv\|^3}\,f_0(\bfx,\bfv)\,\dD\bfx\,\dD\bfv
\ \times\ \int_0^t\|G(s,\cdot)-F^\eps(s,\cdot)\|_{\dot{W}^{-1,1}}\dD s
\end{align*}
and the Gr\"onwall lemma achieves the proof.
\end{proof}

\subsection{Second-order asymptotics}\label{s:2nd-nl}

The analysis of the second-order asymptotics is significantly more involved. Consequently in the present subsection we will enforce a few assumptions beyond those of the external field case. 

To begin with, we stress that in order to build on Proposition~\ref{prop:4.9} we need to bound $\d_t^2\bE^\eps=-\bK\star_\bfx \Div_\bfx(\,\d_t\bfj^\eps)$ uniformly with respect to $\eps$. Unlike bounds on $\d_t\bE^\eps$, this leads to restrictions on initial data. To discuss the corresponding consequences let us from now on consider $\eps$-dependent initial data $f_0^\eps$ and set
\[
\bfj_0^\eps(\bfx)\,=\,
\int_{\R^3}\,\bfv\,f_0^\eps(\bfx,\bfv)\,\dD \bfv\,,\qquad\qquad
\bfj_{\perp,0}^\eps(\bfx)\,=\,
\int_{\R^3}\,\bvperp(0,\bfx,\bfv)\,f_0^\eps(\bfx,\bfv)\,\dD \bfv\,.
\]

To see how the aforementioned restrictions arise, note that from a direct integration stems
\[
\d_t\bfj^\eps(t,\bfx)
\,=\,\frac{(B\bJ)(t,\bfx)}{\eps}\bfj^\eps(t,\bfx)
+\rho^\eps(t,\bfx)\,\bE^\eps(t,\bfx)
-\Div_\bfx\left(\int_{\R^3}f^\eps(t,\bfx,\bfv)\,\bfv\otimes\bfv\dD \bfv\right)\,.
\]
Specializing the latter inequality to initial time shows that there is little hope to bound $\d_t^2\bE^\eps$ uniformly with respect to $\eps$ if $\bfj_{\perp,0}^\eps/\eps$ is not uniformly bounded. Yet we also need to be able to propagate this condition on a time interval independent of $\eps$. To study this particular point, we rewrite the equation on $\bfj^\eps$,  following the strategy applied so far on characterstics, as
\begin{align*}\ds
\d_t\widetilde{\bfj}^\eps
&\ds
\,=\,\frac{(B\bJ)}{\eps}\widetilde{\bfj}^\eps
+\eps\,\frac{\bJ}{B}\left(
-\d_t(\rho^\eps\bE^\eps)
+\Div_\bfx\left(\int_{\R^3}\d_tf^\eps(\cdot,\cdot,\bfv)\,\bfv\otimes\bfv\dD \bfv\right)\right)
\\
&\ds\quad
+\eps\,\d_t\left(\frac{\bJ}{B}\right)\left(
-\rho^\eps\bE^\eps
+\Div_\bfx\left(\int_{\R^3}f^\eps(\cdot,\cdot,\bfv)\,\bfv\otimes\bfv\dD \bfv\right)\right)
\end{align*}
for 
\[
\widetilde{\bfj}^\eps
\,=\,
\bfj^\eps-\eps\,\frac{\rho^\eps}{B}\bJ\bE^\eps
+\eps\,\frac{\bJ}{B}\left(\Div_\bfx\left(\int_{\R^3}f^\eps(\cdot,\cdot,\bfv)\,\bfv\otimes\bfv\dD \bfv\right)\right)\,.
\]
This suggests that to carry on the argument one should assume some initial control on $\d_tf^\eps$ and propagate it over time. 

At this stage it should also be clear to the reader that the latter strategy is essentially equivalent to assuming initially and propagating a uniform control on 
\[
\Div_\bfv\left(f^\eps\,\frac{\bfv\wedge \eDpar(t,\bfx)}{\eps}\right)
\,=\,\frac{\bfv\wedge \eDpar(t,\bfx)}{\eps}\cdot\nabla_\bfv f^\eps\,.
\]
The link between the latter and the bound on $\bfj_\perp$ is even easier to derive from
\[
\bfj_\perp(t,\bfx)\,=\,
-\int_{\R^3}\,\bfv\quad(\bfv\wedge \eDpar(t,\bfx))\cdot\nabla_\bfv\,f^\eps(t,\bfx,\bfv)\dD \bfv\,.
\]

\subsubsection{Uniform bounds on derivatives}

Since this is a key part of the argument, before going on with the derivation of second-order asymptotics, we focus on the propagation of the well-prepared character. Our scheme here is to interpret preparation of data as a condition on the smallness of the derivative with respect to an angle encoding fast rotation and benefit from separation of fast and slow dynamics to check that the slow part cannot destroy this condition on fast-angle dependency. 

The resulting precise statement is as follows.

\bpr\label{p:bd-nl}
Assume $\bB\in W^{2,\infty}$ is such that $1/B\in W^{2,\infty}$ and $\eDpar\in W^{3,\infty}$, $\bK\in W^{3,\infty}$ and $\urho\in W^{1,\infty}(\R_+;\cM_+(\R^3))$. For any $p_0\in(1,\infty)$, $R_0>0$, $M>0$ and $T>0$, there exist positive constants $C$ and $\eps_0$ such that if $0<\eps<\eps_0$ and $f^\eps$ solves \eqref{eq:nl-vlasov} with initial data $f_0^\eps\in BV_+\cap W^{1,p_0}$ such that
\begin{align*}
&\ds\int_{\R^6}\,f_0^\eps(\bfx,\bfv)\,\dD \bfx\,\dD \bfv\,\leq\, M,&\qquad\qquad
&\ds\operatorname{supp}f_0^\eps\quad\subset\quad
\{\ (\bfx,\bfv)\ ;\ \|\bfv\|\leq R_0\ \}\,,
\end{align*}
then for any $0\leq t\leq T$
\begin{align*}\ds
\|\nabla_{\bfx,\bfv}f^\eps(t,\cdot,\cdot)&\ds\|_{L^{p_0}(\R^6)}
+\frac{1}{\eps}\|(\bfx,\bfv)\mapsto(\bfv\wedge\eDpar(t,\bfx))\cdot\nabla_{\bfv}f^\eps(t,\bfx,\bfv)\|_{L^{p_0}(\R^6)}\\
&\ds\leq
C\,\left(
\|\nabla_{\bfx,\bfv}f_0^\eps\|_{L^{p_0}(\R^6)}
+\frac{1}{\eps}\|(\bfx,\bfv)\mapsto(\bfv\wedge\eDpar(t,\bfx))\cdot\nabla_{\bfv}f_0^\eps(\bfx,\bfv)\|_{L^{p_0}(\R^6)}\right)\,.
\end{align*}
\epr

\begin{proof}
We first recall that from Lemma~\ref{bnd0} stems a uniform bound $R$ for $\|\bfv\|$ on the support of $f^\eps(t,\cdot,\cdot)$, $0\leq t\leq T$, $\eps>0$, and that we have already derived bounds on $\bE^\eps$, $\d_t\bE^\eps$ and their spatial derivatives. 

To carry out the proof, we shall introduce plane coordinates for $\bvperp(t,\bfx,\bfv)$ and correct them according to the gyrocenter dynamics. This requires a (smooth) consistent choice of frames on the planes $\{\eDpar(t,\bfx)\}^\perp$. Thus we pick\footnote{See comments in Remark~\ref{rk:frame}.} $\eDa\in W^{2,\infty}$ and $\eDb\in W^{2,\infty}$ such that $(\eDa,\eDb,\eDpar)$ form a field of direct orthonormal frames. This being done, we define $\bfu(t,\bfx,\bfv)\in \R^2$ through
\be
\label{def:u}
\bfu\,=\,\bfsigma(t,\bfx)^*\,\bfv\,,\qquad
\bfsigma\,=\,\begin{pmatrix}\eDa&\eDb
\end{pmatrix}\,,
\ee
where ${}^*$ denotes the adjoint operator. Note that if $(\bfx,\bfv)$ solves \eqref{e:xv} then the corresponding $\bfu$ solves
\be
\label{e:u}
\ds\frac{\dD\bfu}{\dD t}\,=\,\ds
\frac{b^\eps(t,\bfx,\bfu)}{\eps}\, \bJ_0\,\bfu
\,+\,\bfsigma^*(t,\bfx)\,\bF_0(t,\bfx,\vpar)\,+\,\bA_0(t,\bfx,\vpar)\,\bfu
\ee
where 
\be\label{def:J0}
\bJ_0:=\begin{pmatrix}0&1\\-1&0\end{pmatrix}\,,
\ee
$\bF_0$ is as in \eqref{def:F}, and
\begin{align}
\label{def:b0}
b^\eps(t,\bfx,\bfu)&\ds:=
B(t,\bfx)-\eps\,\langle\eDa(t,\bfx),
\dD_\bfx\eDb(t,\bfx)\,\bfsigma(t,\bfx)\bfu\rangle\,,
\\
\label{def:bA0}
\bA_0(t,\bfx,\vpar)&:=\ds
\left(\d_t\bfsigma(t,\bfx)+\vpar\,\dD_\bfx\bfsigma(t,\bfx)\eDpar(t,\bfx)\right)^*\,\bfsigma(t,\bfx)\\\nn
&\ds\qquad-\vpar\,(\bfsigma(t,\bfx))^*\dD_\bfx\eDpar(t,\bfx)\,\bfsigma(t,\bfx)\,.
\end{align}
Now to replace the change from $\eperp$ to $\eGC^\eps$, we would like to identify a gyrokinetic correction to $\bfu$ so as to ensure that at leading order the norm of its correction satisfies an equation uncoupled from any angle defining $\bfu$.

Since we do not need the full algebraic details of the involved computations, instead of writing explicitly the underlying abstract lemma, providing counterparts to Lemmas~\ref{1st} and~\ref{2nd}, we simply point out that from \eqref{e:u} stem
\begin{align*}
\bfsigma^*\bF_0
&=-\frac{B}{\eps}\, \bJ_0\left(\frac{\eps}{B}\,\bJ_0\,\bfsigma^*\bF_0\right)
+\frac{\dD}{\dD t}\left(\frac{\eps}{B}\,\bJ_0\,\bfsigma^*\bF_0\right)
-\eps\frac{\dD}{\dD t}\left(\frac{1}{B}\,\bJ_0\,\bfsigma^*\bF_0\right)
\end{align*}
and
\begin{align*}
\bA_0\bfu
&=\frac12(\bA_0-\bJ_0\bA_0\bJ_0)\bfu
-\frac{B}{\eps}\, \bJ_0\left(\frac{\eps}{2B}\,\bJ_0\,\bA_0\bfu\right)
+\frac{\dD}{\dD t}\left(\frac{\eps}{2B}\,\bJ_0\,\bA_0\bfu\right)\\
&\quad
-\eps\left(\frac{b^\eps-B}{2\eps B}\bJ_0\bA_0\bJ_0\bfu
+\frac{1}{2B}\,\bJ_0\,\bA_0\left(\bfsigma^*\bF_0+\bA_0\bfu\right)
+\frac{\dD}{\dD t}\left(\frac{1}{2B}\,\bJ_0\,\bA_0\right)\bfu\right)
\end{align*}
and that $\bA_0-\bJ_0\bA_0\bJ_0$ commutes with $\bJ_0$. Incidentally, for the sake of consistency with Lemma~\ref{2nd}, we observe that
\[
\langle\bfu,(\bA_0-\bJ_0\bA_0\bJ_0)\bfu\rangle
\,=\,\Tr(\bA_0)\,\frac{\|\bfu\|^2}{2}\,.
\]
The upshot of the previous considerations is the introduction of $\uGC^\eps(t,\bfx,\bfv)$ defined as
\be\label{def:uGC}
\uGC^\eps
:=\bfu-\,\frac{\eps}{B}\,\bJ_0\,\bfsigma^*\bF_0-\frac{\eps}{2B}\,\bJ_0\,\bA_0\,\bfu\,.
\ee

Note that, for some $\eps_0>0$ independent of $f_0^\eps$ (satisfying the conditions of the proposition), the function 
\[
(t,\bfx,\bfv)\mapsto (t,\GC^\eps(t,\bfx,\bfv),\vGC^\eps(t,\bfx,\bfv),\uGC^\eps(t,\bfx,\bfv))
\] 
defined on $\{\ (t,\bfx,\bfv)\in[0,T]\times\R^6\ ;\ \|\bfv\|\leq R\ \}$ is a bi-Lipschitz map uniformly with respect to $\eps\in(0,\eps_0)$. Thus, for $\eps\in(0,\eps_0)$, we may define $g^\eps$ through
\[
f^\eps(t,\bfx,\bfv)\,=\,g^\eps(t,\GC^\eps(t,\bfx,\bfv),\vGC^\eps(t,\bfx,\bfv),\eGC^\eps(t,\bfx,\bfv),\theta(t,\bfx,\bfv))\,.
\] 
and we observe that for some uniform constant $C_0$, 
\begin{align*}\ds
\|\nabla_{\bfx,\bfv}f^\eps(t,\cdot,\cdot)&\ds\|_{L^{p_0}}
+\frac{1}{\eps}\|(\bfx,\bfv)\mapsto(\bfv\wedge\eDpar(t,\bfx))\cdot\nabla_{\bfv}f^\eps(t,\bfx,\bfv)\|_{L^{p_0}}\\
&\ds\leq
C_0\,\left(
\|\nabla_{\bY}g^\eps(t,\cdot)\|_{L^{p_0}}
+\frac{1}{\eps}\|
\bY\mapsto(\bJ_0\bfu)\cdot\nabla_{\bfu}g^\eps(t,\bY)\|_{L^{p_0}}\right)\,,\\
\|\nabla_{\bY}g^\eps(0,\cdot)&\|_{L^{p_0}}
+\frac{1}{\eps}\|\bY\mapsto(\bJ_0\bfu)\cdot\nabla_{\bfu}g^\eps(0,\bY)\|_{L^{p_0}}\\
&\ds\leq
C_0\,\left(
\|\nabla_{\bfx,\bfv}f_0^\eps\|_{L^{p_0}}
+\frac{1}{\eps}\|(\bfx,\bfv)\mapsto(\bfv\wedge\eDpar(t,\bfx))\cdot\nabla_{\bfv}f_0^\eps(\bfx,\bfv)\|_{L^{p_0}}\right)\,,
\end{align*}
(under the assumptions of the proposition, including $\eps\in(0,\eps_0)$)
where $\bY=(\bfy,v,\bfu)$. 

The gain from the gyrokinetic corrections is that $g^\eps$ solves an equation of the form
\[
\d_tg^\eps+(\cU_0^\eps(t,\bY)
+\eps\,\cU_\bY^\eps(t,\bY))\cdot\nabla_\bY g^\eps
+\left(\frac1\eps B(t,\bfy)+b_\theta^\eps(t,\bY)\right)\,(\bJ_0\bfu)\cdot\nabla_\bfu g^\eps\,=\,0\,,
\]
where again $\bY=(\bfy,v,\bfu)$, with
\begin{itemize}
\item $\cU_\bY^\eps(t,\bY))\cdot\nabla_\bY$ and $(\bJ_0\bfu)\cdot\nabla_\bfu$ commuting;
\item $\nabla_\bY\cU_0^\eps$, $\nabla_{\bY}\cU_\bY^\eps$ and $\nabla_{\bY}b_\theta^\eps$ uniformly bounded in $L^\infty$ (on the support of $g^\eps$). 
\end{itemize}
At this stage, since $\nabla_{\bfy,v}$ obviously commutes with $(\bJ_0\bfu)\cdot\nabla_\bfu$, we only need to pick a version of $\nabla_\bfu$ commuting with $(\bJ_0\bfu)\cdot\nabla_\bfu$ so as to complete the proof by direct estimates. With this aim in mind, we point out that $(-\Delta_{\bfu})^{1/2}$ commutes with $(\bJ_0\bfu)\cdot\nabla_\bfu$ and that $\|\nabla_\bfu(\cdot)\|_{L^{p_0}}$ and $\|(-\Delta_{\bfu})^{1/2}(\cdot)\|_{L^{p_0}}$ are equivalent semi-norms (by standard Calder\'on-Zygmund elliptic regularity theory since $p_0\in(1,\infty)$). Moreover, to estimate harmless corresponding commutators we use the following Kato-Ponce type commutator estimate\footnote{See for instance the case $s=1$ in \cite[Theorem~5.1]{Li_Kato-Ponce} and recall that $\|(-\Delta_{\bfu})^{1/2}(\cdot)\|_{BMO}\leq C_0 \|\nabla_{\bfu}(\cdot)\|_{L^\infty}$ (since Riesz transforms map $L^\infty$ to $BMO$ continuously).}
\[
\|(-\Delta_{\bfu})^{1/2}(f\,g)-f\,(-\Delta_{\bfu})^{1/2}(g)\|_{L^{p_0}}
\leq C_0\,\|\nabla_\bfu f\|_{L^\infty}\,\|g\|_{L^{p_0}}\,,
\]
for some $C_0$ independent of $f$ and $g$. Therefore differentiating the equation for $g^\eps$ with $\nabla_{\bfy,v}$, $(-\Delta_{\bfu})^{1/2}$ and $\eps^{-1}(\bJ_0\bfu)\cdot\nabla_\bfu$, and applying a Gr\"onwall argument achieves the proof.
\end{proof}

\br
Note that the foregoing proof is the only place where we perform a change of variables instead of pushing forward. The main reason is that we are aiming here at preserving throughout the transport nature of \eqref{eq:nl-vlasov} --- that comes with cheap tracking of derivatives --- instead of its conservative character. This gain is however inessential and we could have kept our usual point of view up to a few minor changes. 
\er

\br
The foregoing proof is also the only place where instead of performing manipulations on $\|\bvperp\|^2$ we work directly with $\bvperp$. Alternatively, to reduce as much as possible \emph{new} technical considerations (and spare the use of tools from harmonic analysis) one could opt for a framework in which polar coordinates are non-singular and restrict to initial data such that
\[
\operatorname{supp}f_0^\eps\quad\subset\quad
\{\ (\bfx,\bfv)\ ;\ \|\bfv\|\leq R_0\,,\quad \|\bvperp(0,\bfx,\bfv)\|\geq r_0\ \}
\]
with $0<r_0<R_0$ fixed. If one is willing to pay this price, then one may define the angle $\theta(t,\bfx,\bfv)\in \R/(2\pi\Z)$ wherever $\bfv$ is not colinear with $\eDpar(t,\bfx)$  through
\be
\label{def:theta}
\bfv\,=\,\vpar(t,\bfx,\bfv)\,\eDpar(t,\bfx)
\,+\,\sqrt{2\wperp(t,\bfx,\bfv)}\,\left(
\cos(\theta)\,\eDa(t,\bfx)
+\sin(\theta)\,\eDb(t,\bfx)\right)
\ee
and define $h^\eps$ through
\[
f^\eps(t,\bfx,\bfv)\,=\,h^\eps(t,\GC^\eps(t,\bfx,\bfv),\vGC^\eps(t,\bfx,\bfv),\eGC^\eps(t,\bfx,\bfv),\theta(t,\bfx,\bfv))\,.
\] 
The resulting equation for $h^\eps$ takes the form
\[
\d_th^\eps+(\cV_0^{\bE^\eps}(t,\bZ)
+\eps\,\cV_\bZ^\eps(t,\bZ,\theta))\cdot\nabla_\bZ h^\eps
+\left(\frac1\eps B(t,\bfy)+\cV_\theta^\eps(t,\bZ,\theta)\right)\,\d_\theta h^\eps\,=\,0\,,
\]
where $\bZ=(\bfy,v,w)$, with $\nabla_{\bZ,\theta}\cV_\bZ^\eps$ and $\nabla_{\bZ,\theta}\cV_\theta^\eps$ uniformly bounded in $L^\infty$ (on the support of $h^\eps$). Therefore, under this more stringent assumption, differentiating the equation for $h^\eps$ with $\nabla_\bZ$ and $\eps^{-1}\d_\theta$ and applying a Gr\"onwall argument yield the sought estimates. A small gain is that one may derive in this way $L^{p_0}$-estimates of derivatives for any $p_0\in[1,\infty]$. Yet we feel that the gain in simplicity is not worth the extra unnatural restriction.
\er

\br\label{rk:frame}
The foregoing proof is the only one where we crucially use that $\eDpar(t,\cdot)$
is defined on the \emph{whole} $\R^3$, or put in other words, that $\bB(t,\cdot)$ is defined and non vanishing on $\R^3$. In other places we could have assumed that such a $\bB$ were given on a domain sufficient to contain the support of $f^\eps$. Indeed, in general, the possibility to extend $\eDpar$ into a frame field $(\eDa,\eDb,\eDpar)$ may be constrained by topological obstructions. No such obstruction arise on contractible domains, such as $\R^3$, nor on domains that are homotopic to the circle $\bS^1$, such as $\R^3$ minus a line, or $\R^3$ minus a cylinder (as considered in Section~\ref{s:axi}), since the fundamental group of the sphere $\bS^2$ is trivial. On more general domains $\Omega$, the existence of such a frame choice may be seen as an extra constraint on $\bB$ satisfied when maps $\eDpar (t,\cdot):\Omega\to\bS^2$ are topologically trivial. A typical example of obstruction arises from the case where $\eDpar$ would be nowhere tangent to some (topological) sphere but in this case $\eDpar$ cannot be defined in all the interior of the sphere and thus the domain must contain a hole. Note however that it is sufficient to have a (smooth) consistent choice of a vector-field $\bfz$ nowhere colinear with $\eDpar$ to derive the direct orthonormal frame $\left(\eDpar,\bJ^2\bfz/\|\bJ\bfz\|,\bJ\bfz/\|\bJ\bfz\|\right)$. In particular, $\dD\eDpar\,\eDpar$ provides such a vector field wherever it does not vanish. Moreover a large class of confining geometries are precisely designed to ensure the existence of a smooth level set function\footnote{Typically with compact level sets.} $\bfx\mapsto\psi(\bfx)$ such that $\nabla_\bfx\psi$ provides on the domain of interest a vector-field nowhere vanishing and everywhere orthogonal to $\eDpar$. Thus, for the practical cases that we have in mind this does not appear as a strong constraint.
\er

\subsubsection{Asymptotics}

From the conclusions of Proposition~\ref{p:bd-nl} we derive a nonlinear counterpart to Theorem~\ref{th:2} for well-prepared data.

\bt
\label{th:2-nl}
Let $p_0\in(1,\infty)$ and $p_0'$ be its Lebesgue conjugate, $1/p_0+1/p_0'=1$.\\
Assume $\bB\in W^{2,\infty}$ is such that $1/B\in W^{2,\infty}$ and $\eDpar\in W^{3,\infty}$, $\bK\in W^{3,\infty}\cap W^{1,p_0'}$ and $\urho\in W^{2,\infty}(\R_+;\cM_+(\R^3))$. For any $R_0>0$, $M>0$ and $T>0$, there exist positive constants $C$ and $\eps_0$ such that if $0<\eps<\eps_0$ and $f^\eps$ solves \eqref{eq:nl-vlasov} with initial data $f_0^\eps\in BV_+\cap W^{1,\infty}$ such that
\begin{align*}
&\ds\,\|f_0^\eps\|_{L^1(\R^6)}\,+\,\|\nabla_{\bfx,\bfv}f_0^\eps\|_{L^{p_0}(\R^6)}\,\leq\, M\,,\\[0.5em]
&\ds\,
\|(\bfx,\bfv)\mapsto(\bfv\wedge\eDpar(t,\bfx))\cdot\nabla_{\bfv}f_0^\eps(\bfx,\bfv)\|_{L^{p_0}(\R^6)}\,\leq\,M\,\eps\,,\\[0.5em]
&\ds\operatorname{supp}f_0^\eps\quad\subset\quad
\{\ (\bfx,\bfv)\ ;\ \|\bfv\|\leq R_0\ \}\,,
\end{align*}
then the density $F^\eps$ defined by
$$
F^\eps(t,\cdot)\,=\, \ZGC^{\bE^\eps,\,\eps}(t,\cdot)_* \,(f^\eps(t,\cdot))\,,
$$
where $\ZGC^{\bE^\eps,\,\eps}$ is defined through \eqref{GC}-\eqref{eGC} (but with $\bE^\eps$ the electric field generated by $f^\eps$), satisfies for any $0\leq t\leq T$
$$
\|F^\eps(t,\cdot)-G^\eps(t,\cdot)\|_{\dot{W}^{-1,1}}\leq C\,\eps^2\,,
$$
where $G^\eps$ solves 
\be
\label{eq:2nd-nl}
\left\{
\begin{array}{l}\ds
\d_t G^\eps\,+\,\Div_\bfZ \left(\cV^{\bE_\eps,\,\eps}\,G^\eps \right) \,=\,0\,,\\[0.5em]\ds
\bE_\eps(t,\bfx)=(\bK\star_{\bfx}(\rho_\eps(t,\cdot)-\urho(t,\cdot)))(x)\,,\quad 
\rho_\eps(t,\bfx)=\int_{\R\times\R_+}\,G^\eps(t,\bfx,v,w)\dD v\,\dD w\,,
\end{array}
\right.
\ee
with initial datum $G_{0}^\eps$
\be
\label{eq:2nd-nlG0}
G_{0}^\eps\,=\, \bZ_{\rm gc}^{\bE^\eps,\,\eps}(0,\cdot)_* \,(f_0)\,,
\ee
where $\cV^{\bE_\eps,\,\eps}=\cV^{\bE_\eps}_0+\,\eps\,\cV^{\bE_\eps}_1$ is still given by formula \eqref{V:eps} but $\bE_\eps$ is now a self-consistent electric field.
\et

\br
The condition encoding the well-prepared character of the initial data mixes spatial and kinetic variables in an intricate way. Yet an easy way to enforce it is to take initial data that are radial in velocity (up to a term of order $\eps$).
\er

\begin{proof}
The scheme of the proof is identical to the one of Theorem~\ref{th:1-nl} so we only stress important departures from the latter. We recall that the strategy relies on two intermediate comparisons with the solution of an equation similar to \eqref{eq:2nd-nl} where $\cV^{\bE_\eps,\,\eps}$ is replaced with $\cV^{\bE^\eps,\,\eps}$ (with as above $\bE^\eps$ associated with $f^\eps$ and $\bE_\eps$ associated with $G^\eps$).

In the comparison with $F^\eps$, the main new ingredient is Proposition~\ref{p:bd-nl} that provides uniform bounds on $\bfj_\perp^\eps(t,\cdot)$ in $L^{p_0}$ hence on $\d_t\bfj^\eps(t,\cdot)$ in $L^{p_0}\cap W^{-1,1}$ and thus on $\d_t^2\bE^\eps(t,\cdot)$ in $L^\infty$. This allows to derive the first intermediate comparison from Proposition~\ref{prop:4.9}.

In the comparison with $G^\eps$, the only significantly new constraint is that we need a comparison of $(\bE^\eps)_\mypar$ with $(\bE_\eps)_\mypar$ at order $\eps^2$. This follows from the identity
\begin{align}\label{e:cancel}\ds
&\ds\int \varphi(\bfx)\,\rho^\eps(t,\bfx)\,\dD\bfx
\,=\, 
\int \varphi(\bfy)\, \bF^\eps(t,\bZ)\,\dD\bZ
-\eps\,\int \nabla_\bfx\varphi(\bfx)\cdot\frac{\bJ(t,\bfx)}{B(t,\bfx)}\bfj_\perp^\eps(t,\bfx)\,\dD\bfx\\\nn
&\ds
-\eps^2\int_0^1\left(
\int \dD_\bfx^2\varphi\left(\bfx+\eps\,s\,\frac{\bJ(t,\bfx)\bfv}{B(t,\bfx)}\right)\left(\frac{\bJ(t,\bfx)\bfv}{B(t,\bfx)},\frac{\bJ(t,\bfx)\bfv}{B(t,\bfx)}\right)
\,f^\eps(t,\bfx,\bfv)\,\dD \bfx\dD\bfv\right)(1-s)\dD s
\end{align}
which, after an integration by parts in the last term, gives
\[
\|\rho^\eps(t,\cdot)-\rho_\eps(t,\cdot)\|_{\dot{W}^{-1,1}+\dot{W}^{-1,p_0'}}
\leq 
\|F^\eps(t,\cdot)-G^\eps(t,\cdot)\|_{\dot{W}^{-1,1}}
+C_0\,\eps^2
\]
for some harmless $C_0$, since both $\bfj_\perp^\eps/\eps$ and $\nabla_{\bfx,\bfv}f^\eps$ are uniformly bounded in $L^{p_0}$. This is sufficient to conclude the proof.
\end{proof}

As hinted at by the computation \eqref{e:cancel}, in the well-prepared case considered in Theorem~\ref{th:2-nl} it is also possible to obtain a second-order description for the density of original first-order slow variables, that is, without the guiding-center correction.

\bc
\label{cor:2-nl}
Let $p_0\in(1,\infty)$ and $m_0>3(1-1/p_0)$.\\
Assume $\bB\in W^{2,\infty}$ is such that $1/B\in W^{2,\infty}$ and $\eDpar\in W^{3,\infty}$, $\bK\in W^{3,\infty}$ and $\urho\in W^{2,\infty}(\R_+;\cM_+(\R^3))$. For any $R_0>0$, $M>0$ and $T>0$, there exist positive constants $C$ and $\eps_0$ such that if $0<\eps<\eps_0$ and $f^\eps$ solves \eqref{eq:nl-vlasov} with nonnegative initial data $f_0^\eps\in W^{1,1}\cap W^{1,\infty}$ such that
\begin{align*}
&\ds\,
\|f_0^\eps\|_{L^1(\R^6)}
\,+\,\|(\bfx,\bfv)\mapsto (1+\|\bfx\|^{m_0})\nabla_{\bfx,\bfv}f_0^\eps(\bfx,\bfv)\|_{L^{p_0}(\R^6)}\,\leq\, M\,,\\[0.5em]
&\ds\,
\|(\bfx,\bfv)\mapsto (1+\|\bfx\|^{m_0})\,
(\bfv\wedge\eDpar(t,\bfx))\cdot\nabla_{\bfv}f_0^\eps(\bfx,\bfv)\|_{L^{p_0}(\R^6)}\,\leq\,M\,\eps\,,\\[0.5em]
&\ds\operatorname{supp}f_0^\eps\quad\subset\quad
\{\ (\bfx,\bfv)\ ;\ \|\bfv\|\leq R_0\ \}\,,
\end{align*}
then $F^\eps$ defined by
$$
F^\eps(t,\bfx,v_\mypar,\eperp)=\int_{\bS_{t,\bfx}}\,f^\eps(t,\bfx,v_\mypar\,\eDpar(t,\bfx)\,+\,\sqrt{2\,\eperp}\ \widehat{\beD})\ \dD \sigma_{t,\bfx}(\widehat{\beD}),
$$
with $\bS_{t,\bfx}=\{\eDpar(t,\bfx)\}^\perp\cap\bS^2$ and $\sigma_{t,\bfx}$ its canonical line-measure, satisfies for any $0\leq t\leq T$
$$
\|F^\eps(t,\cdot)-G^\eps(t,\cdot)\|_{\dot{W}^{-1,1}}\leq C\,\eps^2\,,
$$
where $G^\eps$ solves 
\be
\left\{
\begin{array}{l}\ds
\d_t G^\eps\,+\,\Div_\bfZ \left(\cV^{\bE_\eps,\,\eps}\,G^\eps \right) \,=\,0\,,\\[0.5em]\ds
\bE_\eps(t,\bfx)=(\bK\star_{\bfx}(\rho_\eps(t,\cdot)-\urho(t,\cdot)))(x)\,,\quad 
\rho_\eps(t,\bfx)=\int_{\R\times\R_+}\,G^\eps(t,\bfx,v,w)\dD v\,\dD w\,,
\end{array}
\right.
\ee 
with initial data $G_{0}^\eps$ 
$$
G_0^\eps(t,\bfx,v_\mypar,\eperp)=\int_{\bS_{0,\bfx}}\,f_0^\eps(\bfx,v_\mypar\,\eDpar(0,\bfx)\,+\,\sqrt{2\,\eperp}\ \widehat{\beD})\ \dD \sigma_{0,\bfx}(\widehat{\beD}),
$$
where $\cV^{\bE_\eps,\,\eps}$ is as in \eqref{V:eps} but with the self-consistent $\bE_\eps$.
\ec

\br If one completes $\eDpar$ into a frame field $(\eDa,\eDb,\eDpar)$ (as in the proof of  Proposition~\ref{p:bd-nl}), then the definitions of $F^\eps$ and $G_0^\eps$ may be equivalently written as
\begin{align*}
F^\eps(t,\bfx,v_\mypar,\eperp)&=\frac{1}{2\pi}\int_0^{2\pi}\,f^\eps(t,\bfx,v_\mypar\,\eDpar(t,\bfx)\,+\,\sqrt{2\,\eperp}\ \eD^\theta(t,\bfx))\ \dD \theta,\\[0.5em]
G_0^\eps(\bfx,v_\mypar,\eperp)&=\frac{1}{2\pi}\int_0^{2\pi}\,f_0^\eps(\bfx,v_\mypar\,\eDpar(0,\bfx)\,+\,\sqrt{2\,\eperp}\ \eD^\theta(0,\bfx))\ \dD \theta\,,
\end{align*}
with $\eD^\theta(t,\bfx)=\cos(\theta)\eDa(t,\bfx)+\sin(\theta)\eDb(t,\bfx)$.
\er

\br
As is readily derived from a comparison of the proofs of Theorem~\ref{th:2-nl} and Corollary~\ref{cor:2-nl}, one may remove spatial weights from the assumptions of the latter corollary provided the conclusion is weakened into an estimate of $\|F^\eps(t,\cdot)-G^\eps(t,\cdot)\|_{\dot{W}^{-1,1}+\dot{W}^{-1,p_0}}$.
\er

\begin{proof}
To resolve the notational mismatch between Theorem~\ref{th:2-nl} and Corollary~\ref{cor:2-nl}, we use here the subscript ${}_{\rm gc}$ to denote densities introduced in the former theorem. Thus, to derive Corollary~\ref{cor:2-nl} from Theorem~\ref{th:2-nl} we only need to compare on one hand $F^\eps$ with $F_{\rm gc}^\eps$ and on the other hand $G^\eps$ with $G_{\rm gc}^\eps$. The comparison of the latter is readily derived from a standard stability estimate on \eqref{eq:2nd-nl} and a comparison of their initial data. Now, the comparisons of the $G_0^\eps$ and $(G_{\rm gc})_0^\eps$ on one side and of $F^\eps$ and $F_{\rm gc}^\eps$ on the other side follow essentially from the same argument, which is a variation on \eqref{e:cancel} and its use in the proof of Theorem~\ref{th:2-nl}.

As in \eqref{e:cancel}, we may perform second-order expansions to compare quantities of interest. The only significant change in the expansion is that in the terms of order $\eps$ appear not only $\bfj_\perp^\eps$ but also the scalar
\begin{align*}\ds
\sigma_\perp^\eps(t,\bfx)&\ds
\,=\,
\int_{\R^3}\,
\left\langle\bJ(t,\bfx)\,\bfv,\,\Re(\dD_\bfx\eDpar(t,\bfx))\,\bJ(t,\bfx)^2\,\bfv\right\rangle
f^\eps(t,\bfx,\bfv)\dD \bfv\\
&\ds\,=\,
\frac12\int_{\R^3}\,
\left\langle\bJ(t,\bfx)^2\,\bfv,\,\Re(\dD_\bfx\eDpar(t,\bfx))\,\bJ(t,\bfx)^2\,\bfv\right\rangle
\quad(\bJ(t,\bfx)\bfv)\cdot\nabla_\bfv\,f^\eps(t,\bfx,\bfv)\dD \bfv
\end{align*}
that is controlled exactly as $\bfj_\perp^\eps$. The only departure in uses of the algebraic identities stems from our will to get estimates in $W^{-1,1}$ and not in $W^{-1,1}+W^{-1,\infty}$. 

To achieve this goal, we simply observe one one hand that the proof of Proposition~\ref{p:bd-nl} also shows propagation of $\dot{W}^{1,p_0}$ regularity with spatial weights $(1+\|\bfx\|^2)^{m_0/2}$ and on the other that from $m_0>3(1-1/p_0)$ stems for any $f$, 
\[
\|\,f\,\|_{L^1(\R^3)}\leq C_0 \|x\mapsto (1+\|\bfx\|^2)^{m_0/2}\,f\,\|_{L^{p_0}(\R^3)}
\]
for some constant $C_0$ independent of $f$. Incidentally this also relaxes the condition $\bK\in W^{1,p_0'}$ from assumptions of Theorem~\ref{th:2-nl}.
\end{proof}

Note that though the statement of Corollary~\ref{cor:2-nl} does not involve guiding-center coordinates, trying to prove it without essentially following the proof of Theorem~\ref{th:2-nl} would be rather cumbersome since the cancellation used to prove Corollary~\ref{cor:2-nl} and arising from the well-prepared nature of initial data is only present at the level of densities.

\appendix

\section{Comparison with the classical adiabatic invariant formulation}\label{s:mu}

In the present Section, for the sake of comparison with part of the physical literature, we derive counterparts to our three-dimensional results expressed in terms of slow variables 
$$
\left(\bfx,\vpar(t,\bfx,\bfv),\frac{\eperp(t,\bfx,\bfv)}{B(t,\bfx)}\right)
$$ 
and corrections thereof. To do so we explicitly introduce the function $\muperp$ defined by
$$
\muperp(t,\bfx,\bfv)
\,=\,\frac{\eperp(t,\bfx,\bfv)}{B(t,\bfx)}\\
\,=\,\frac{\|\bvperp(t,\bfx,\bfv)\|^2}{2\,B(t,\bfx)}\,.
$$

\subsection{Long-time asymptotics in the toroidal axi-symmetric case}

Since this is slightly less computationally demanding, we provide first a counterpart to Theorem~\ref{th:3}. Our starting point is System~\eqref{e:rzve}. From it we derive
\begin{align*}\ds
\frac{\dD}{\dD t}\Big[\muperp&\ds
-\eps\,\left(
\muperp\,\left\langle \d_r\left(\frac1b\right)\eDz-\d_z\left(\frac1b\right)\eDr,\bvperp\right\rangle
+\frac1b\left\langle\EcB,\bvperp\right\rangle
+\frac{\vpar^2}{r\,b^2}
\,\left\langle\eDz,\bvperp\right\rangle
\right)
+\eps^2\chi_\mu\Big]\\[0.5em]
&\ds=\,
\eps\,\frac{\muperp}{b}\left(\d_r\Ez-\d_z\Er\right)
+\eps^2\,\eta_\mu\\[0.5em]
&\ds-\eps\,\,
\left(\left\langle\EcB,\bvperp\right\rangle
+\frac{\vpar^2}{r\,b}
\,\left\langle\eDz,\bvperp\right\rangle
\right)
\,\left\langle
\d_r\left(\frac1b\right)\eDr
+\d_z\left(\frac1b\right)\eDz
,\bvperp\right\rangle\\[0.5em]
&\ds
-\eps\muperp\,\langle\eDz,\bvperp\rangle\,\left(
\d_r^2\left(\frac1b\right)\,\langle\eDr,\bvperp\rangle
+\d_{r\,z}^2\left(\frac1b\right)\,\langle\eDz,\bvperp\rangle
\right)
\\[0.5em]
&\ds
+\eps\muperp\,\langle\eDr,\bvperp\rangle\,
\left(\d_{r\,z}^2\left(\frac1b\right)\,\langle\eDr,\bvperp\rangle
+\d_z^2\left(\frac1b\right)\,\langle\eDz,\bvperp\rangle
\right)
\\[0.5em]
&\ds
+\frac{\eps}{b}\,
\left\langle\d_z\left(\frac1b\right)\eDr
-\d_r\left(\frac1b\right)\eDz,
\bvperp\right\rangle\left\langle \Eperp(t,\bfx)\,-\,\vpar\,
\dD_{\bfx}\eDpar(t,\bfx)\,\bfv,
\bvperp\right\rangle
\end{align*}
with
$$
\chi_\mu
\,=\,
\left\langle\d_r\left(\frac1b\right)\eDr
+\d_z\left(\frac1b\right)\eDz,\bfchi_x\right\rangle
+\frac{\chi_\perp}{b}
$$
and
\begin{align*}
\eta_\mu&\ds=\,
\eperp\,\left(\left\langle\d_r\left(\frac1b\right)\eDr+\d_z\left(\frac1b\right)\eDz,\bfeta_x\right\rangle+\frac{\vpar}{r}\d_r\left(\frac1b\right)\langle\eDpar,\bfchi_x\rangle\right)
+\frac{\eta_\perp}{b}\\[0.5em]
&\ds
+\left\langle
\d_r\left(\frac1b\right)\eDr
+\d_z\left(\frac1b\right)\eDz,\bfchi_x\right\rangle
\left\langle \Eperp(t,\bfx)\,-\,\vpar\,
\dD_{\bfx}\eDpar(t,\bfx)\,\bfv,
\bvperp\right\rangle
\\[0.5em]
&\ds
+\eperp\,\langle\eDr,\bfchi_x\rangle\,\left(
\d_r^2\left(\frac1b\right)\,\langle\eDr,\bvperp\rangle
+\d_{r\,z}^2\left(\frac1b\right)\,\langle\eDz,\bvperp\rangle
\right)
\\[0.5em]
&\ds
+\eperp\,\langle\eDz,\bfchi_x\rangle\,\left(
\d_{r\,z}^2\left(\frac1b\right)\,\langle\eDr,\bvperp\rangle
+\d_z^2\left(\frac1b\right)\,\langle\eDz,\bvperp\rangle
\right)
\\[0.5em]
&\ds
+\chi_\perp\left\langle
\d_r\left(\frac1b\right)\eDr
+\d_z\left(\frac1b\right)\eDz
,\bvperp\right\rangle\,.
\end{align*}

Now we note that the extra $\cO(\eps)$-terms in the right-hand side of the foregoing system are either third-order with respect to $\bvperp$, or second-order but trace-free in the plan orthogonal to $\eDpar$ (when suitably paired). As a result they may be eliminated, leaving
\begin{align*}\ds
\frac{\dD}{\dD t}\Big[\muperp&\ds
-\eps\,\left(
\muperp\,\left\langle \d_r\left(\frac1b\right)\eDz-\d_z\left(\frac1b\right)\eDr,\bvperp\right\rangle
+\frac1b\left\langle\EcB,\bvperp\right\rangle
+\frac{\vpar^2}{r\,b^2}
\,\left\langle\eDz,\bvperp\right\rangle
\right)
+\eps^2\widehat\chi_\mu\Big]\\[0.5em]
&\ds=\,
\eps\,\frac{\muperp}{b}\left(\d_r\Ez-\d_z\Er\right)
+\eps^2\,\widehat\eta_\mu
\end{align*}
with
\be\label{mu-axi-remainders}
|\widehat\chi_\mu|\lesssim \|\bvperp\|\,(1+\|\bfv\|^3)\,,\qquad\qquad
|\widehat\eta_\mu|\lesssim 1+\|\bfv\|^5\,.
\ee

Therefore, the involved asymptotic vector field is now $\eps\cZ_1$ with $\cZ_1$ defined as
$$
\cZ_1(t,\bfZ)
\,=\,\bp\ds
-\frac{\Ez(t,r,z)}{b(r,z)}+\frac{\mu}{b(r,z)}\,\d_zb(r,z)\\[1em]\ds
\frac{\Er(t,r,z)}{b(r,z)}+\frac{v^2}{r\,b(r,z)}
-\frac{\mu}{b(r,z)}\,\d_rb(r,z)
\\[1em]
\ds\frac{v}{r}\,\left(\frac{\Ez(t,r,z)}{b(r,z)}
-\frac{\mu}{b(r,z)}\,\d_zb(r,z)\right)\\[1em]\ds
\frac{\mu}{b(r,z)}\left(\d_r\Ez-\d_z\Er\right)(t,r,z)
\ep,
$$
where the slow variable is now $\bfZ=(r,z,v,\mu)$. Note that $r\,b\,\cZ_1$ is divergence-free.

\bt
\label{th:3-mu}
Let $\bB$ be a stationary, axi-symmetric and toroidal magnetic field and $\bE$ be
an axi-symmetric electric field orthogonal to $\bB$, with $(\Er,\Ez,1/b)\in
W^{2,\infty}$ in the region where $r(\bfx)\geq r_0$ for some
$r_0$. For any $r_1>r_0$, there exist positive constants $\eps_0$, $\tau_0$ and $C_0$, $(1/\eps_0,1/\tau_0,C_0)$ depending polynomially on $1/r_0$, $1/(r_1-r_0)$ and $\|(\Er,\Ez,1/b)\|_{W^{2,\infty}([r_0,\infty[\times\R)}$, such that the following holds with
$$
\eps_{max}(R_0):=\frac{\eps_0}{1+R_0}\qquad\textrm{and}\qquad
T_{max}(R_0):=\frac{\tau_0}{1+R_0^2}\,.
$$
Consider $f^\eps$ a solution to \eqref{eq:vlasov} with initial datum a nonnegative density $f_0$ supported where 
$$
r(\bfx)\geq r_1\qquad\textrm{and}\qquad \|\bfv\|\leq R_0 
$$
for some $R_0>0$ and define $F^\eps$ as
$$
F^\eps(t,r,z,v,\mu)=\int_0^{2\pi}\!\!\int_0^{2\pi}
\,f^\eps(t,r\,\eDr^\theta+z\eDz,
v\,\eDpar(r\,\eDr^\theta+z\eDz)
\,+\,\sqrt{2\,\mu\,b(r,z)}\ \eDperp^{\theta,\,\varphi})\ r\,b(r,z)\,\dD \varphi\,\dD\theta
$$
with  
$$
\eDr^\theta\,=\,\cos(\theta)\eDx+\sin(\theta)\eDy\,,\qquad\qquad
\eDperp^{\theta,\,\varphi}\,=\,\cos(\varphi)\eDr^\theta+\sin(\varphi)\eDz\,,
$$
where $(\eDx,\eDy,\eDz)$ is a fixed orthonormal basis.  Then provided that
$$
0<\eps\leq \eps_{max}(R_0)\,,
$$
we have for a.e. $0\leq t \leq T_{max}(R_0)/\eps$
$$
\|F^\eps(t,\cdot)-G^\eps(t,\cdot)\|_{\dot{W}^{-1,1}}\leq C_0\,\eps\,
\int_{\R^3\times\R^3} e^{C\,\eps\,t\,\|\bfv\|^4}\,(1+\|\bfv\|^3)\,f_0(\bfx,\bfv)\,\dD\bfx\,\dD\bfv,
$$
where $G^\eps$ solves 
\be\label{eq:axi-mu}
\d_t G^\eps\,+\,\eps\Div_{\bfZ} \left(G^\eps\,\cZ_1\right)
\,=\,0,
\ee
with initial datum $G_0$ given by
$$
G_0(\bfZ)\,=\,
\int_0^{2\pi}\!\!\int_0^{2\pi}
\,f_0(\,r\,\eDr^\theta+z\eDz,
v\,\eDpar(r\,\eDr^\theta+z\eDz)
\,+\,\sqrt{2\,\mu\,b(r,z)}\ \eDperp^{\theta,\,\varphi})\ r\,b(r,z)\,\dD \varphi\,\dD\theta\,.
$$
\et

Incidentally we observe that of course it is easier to derive the balance law
$$
\d_t(\mu G)\,+\,\eps\Div_{\bfZ} \left(\mu\,G^\eps\,\cZ_1\right)
\,=\,\frac{\mu}{b}\,G\,\left(\d_r\Ez-\d_z\Er\right)
$$
for the asymptotic equation \eqref{eq:axi-mu} than the corresponding result for the original formulation. Yet the energy balance law is in turn less straightforward to derive.

At the level of description considered in this section since going from variables $(r,z,\vpar,\eperp)$ to $(r,z,\vpar,\muperp)$ is quite simple we could have deduced Theorem~\ref{th:3-mu} from Theorem~\ref{th:3}. To give two hints in this direction, note that if we denote $F^\eps_e$ and $F^\eps_\mu$ the averaged densities respectively from Theorem~\ref{th:3} and Theorem~\ref{th:3-mu}, then $F^\eps_\mu(t,r,z,v,\mu)\,\dD\!r\,\dD\!z\,\dD\!v\,\dD\!\mu$ is the push-forward of $F^\eps_e(t,r,z,v,\e)\,\dD\!r\,\dD\!z\,\dD\!v\,\dD\!\e$ by the map $(r,z,v,\e)\mapsto (r,z,v,\e/b(r,z))$. Likewise, with the same convention, it may be checked that $G^\eps_\mu(t,r,z,v,\mu)\,\dD\!r\,\dD\!z\,\dD\!v\,\dD\!\mu$ is obtained from $G^\eps_e(t,r,z,v,\e)\,\dD\!r\,\dD\!z\,\dD\!v\,\dD\!\e$ in the same way, by using that
$$
\cZ_1^{\mu}(t,r,z,v,\mu)
\,=\,
\left[\frac{\cW_1^e}{b}
-\frac{\e}{b^2}\left(\d_rb\,\cW_1^r
+\d_zb\,\cW_1^z\right)\right](t,r,z,v,\mu\,b(r,z))
$$
where superscripts denote components. However we have chosen not to follow this path and instead to come back to the normal form \eqref{e:rzve} since the strategy under discussion is more cumbersome when higher-order corrections are taken into account.

\subsection{General second-order asymptotics}

We now provide a counterpart to Theorem~\ref{th:2}. Our starting point is System~\eqref{e:1new3}. From it stems
\begin{align*}
\ds\frac{\dD }{\dD t}\Big[\muperp
&\ds
-\eps\frac{\muperp}{B^2}\dD_\bfx\!B\ \bJ\,\bvperp
-\,\frac{\eps}{B} \,\left\langle\bvperp,
\,\EcB+\vpar\Sig\right\rangle
-\; \frac{\eps\,\vpar}{2B^2} \,\left\langle\bJ\,\bvperp,\,\Re(\dD_\bfx\eDpar)\,\bvperp\right\rangle
+\eps^2\chi_\mu
\Big]\\[1em]
&\ds
=\,-\frac{\muperp}{B}\left(\d_tB+\vpar\Div_\bfx\bB\right)
-\,\eps\,\frac{\muperp}{B}\,\Div_\bfx\left(B\,\EcB+\vpar\,B\,\Sig\right)
\\[1em]
&\ds
-\,\eps\,\muperp\,\left\langle\curvB\,,\,\bE\right\rangle
-\eps\,\muperp^2\dD_\bfx\!B\ \rotB
+\eps^2\eta_\mu\\[1em]
&\ds
+\eps\muperp\dD_\bfx\d_t\left(\frac1B\right)\ \bJ\,\bvperp
+\eps\muperp\dD_\bfx^2\left(\frac1B\right)\left(\bfv,\,\bJ\,\bvperp\right)\\[1em]
&\ds-\frac{\eps}{B^3}\dD_\bfx\!B\ \bJ\bvperp\ \left\langle \Eperp(t,\bfx) \,-\,\vpar\,\left(
\d_t\eDpar(t,\bfx)+\dD_{\bfx}\eDpar(t,\bfx)\,\bfv\right),\bvperp\right\rangle\\[1em]
&\ds-\,
\eps\,\left(\d_t\left(\frac1B\right)+\dD_\bfx\left(\frac1B\right)\bfv\right)
\Big[\left\langle\bvperp,
\,\EcB+\vpar\,\Sig\,\right\rangle
+\; \frac{\vpar}{2B} \,\left\langle\bJ\,\bvperp,\,\Re(\dD_\bfx\eDpar)\,\bvperp\right\rangle\Big]
\end{align*}
with
\begin{align*}
\chi_\mu&
=-\frac{\muperp}{B}\dD_\bfx\!B\ \bfchi_\bfx
+\frac{\chi_\perp}{B}\\[1em]
\eta_\mu&
=-\frac{\muperp}{B}\dD_\bfx\!B\ \bfeta_\bfx+\frac{\eta_\perp}{B}
+\eperp\dD_\bfx\left[\d_t\left(\frac1B\right)\right]\bfchi_\bfx
+\eperp\dD_\bfx^2\left(\frac1B\right)\left(\bfv,\,\bfchi_\bfx\right)
-\frac{\chi_\perp}{B^2}\left(\d_tB+\dD_\bfx\!B\ \bfv\right)\\[1em]
&-\frac{1}{B^2}\dD_\bfx\!B\ \bfchi_\bfx\ \left\langle \Eperp(t,\bfx) \,-\,\vpar\,\left(
\d_t\eDpar(t,\bfx)+\dD_{\bfx}\eDpar(t,\bfx)\,\bfv\right),\bvperp\right\rangle\,.
\end{align*}
The last three lines of the foregoing system may be discarded as being linear or cubic in $\bvperp$, or quadratic in $\bvperp$ but with zero trace in the plane orthogonal to $\bvperp$. This leads to 
\begin{align*}
\ds\frac{\dD }{\dD t}\Big[\muperp
&\ds
-\eps\frac{\muperp}{B^2}\dD_\bfx\!B\ \bJ\,\bvperp
-\,\frac{\eps}{B} \,\left\langle\bvperp,
\,\EcB+\vpar\Sig\right\rangle
-\; \frac{\eps\,\vpar}{2B^2} \,\left\langle\bJ\,\bvperp,\,\Re(\dD_\bfx\eDpar)\,\bvperp\right\rangle
+\eps^2\widehat\chi_\mu
\Big]\\[1em]
&\ds
=\,-\frac{\muperp}{B}\left(\d_tB+\vpar\Div_\bfx\bB\right)
-\,\eps\,\frac{\muperp}{B}\,\Div_\bfx\left(B\,\EcB+\vpar\,B\,\Sig\right)
\\[1em]
&\ds
-\,\eps\,\muperp\,\left\langle\curvB\,,\,\bE\right\rangle
-\eps\,\muperp^2\dD_\bfx\!B\ \rotB
+\eps^2\widehat\eta_\mu
\end{align*}
with
$$
|\widehat\chi_\mu|\lesssim \|\bvperp\|\,(1+\|\bfv\|^3)\,,\qquad\qquad
|\widehat\eta_\mu|\lesssim 1+\|\bfv\|^5\,.
$$

At this stage one may follow the final lines of the proof of Theorem~\ref{th:2}, by replacing in the zeroth and first-order terms of the right-hand side $(\bfx,\vpar,\muperp)$ with $\eps$-corrections $(\GC,\vGC,\muGC)$ up to $\cO(\eps^2)$-terms that may be added to $\eps^2\widehat\eta_\mu$ and $\eps$-terms that may be removed by the by-now familiar elimination process, resulting in another harmless modification of $\widehat\chi_\mu$ and $\widehat\eta_\mu$.

To state the resulting theorem, we introduce 
\be\label{muGC}
\muGC\,=\,
\frac{\eperp}{B}
-\eps\frac{\eperp}{B^3}\dD_\bfx\!B\ \bJ\,\bvperp
-\,\frac{\eps}{B} \,\left\langle\bvperp,
\,\EcB+\vpar\Sig\right\rangle
-\; \frac{\eps\,\vpar}{2B^2}
\,\left\langle\bJ\,\bvperp,\,\Re(\dD_\bfx\eDpar)\,\bvperp\right\rangle
\ee
and
$$
\cY^\eps\,=\,\cY_0+\eps\cY_1
$$
with
$$
\cY_0(t,\bfZ)
\,=\,\bp\ds
v\,\eDpar(t,\bfy)\\[1em]\ds
\Epar(t,\bfy)+\mu\,B(t,\bfy)\,\Div_\bfx\eDpar(t,\bfy)\\[1em]\ds
-\frac{\mu}{B(t,\bfy)}\left(\d_tB+v\Div_\bfx\bB\right)(t,\bfy)
\ep
$$
and
$$
\cY_1(t,\bfZ)
\,=\, \bp\ds
\Ud(t,\bfy,v,\mu\,B(t,\bfy)) 
\\[1em]
\ds\left\langle\Sig(t,\bfy,v),\bE(t,\bfy)\right\rangle
+\mu\,B(t,\bfy)\,\Div_\bfx\Sig(t,\bfy,v)
\\[1em]
\ds 
\left[-\,\frac{\mu}{B}\,\Div_\bfx\left(B\,\EcB+v\,B\,\Sig\right)
-\,\mu\,\left\langle\curvB\,,\,\bE\right\rangle
-\,\mu^2\dD_\bfx\!B\ \rotB\right](t,\bfy,v)
\ep
$$
defining vector fields on the reduced phase-space where $\bfZ=(\bfy,v,\mu)$ lives.

\bt
\label{th:2-mu}
Let $\bE\in W^{2,\infty}$ and $\bB$ be such that $1/B\in W^{2,\infty}$ and $\eDpar\in W^{3,\infty}$. There exists a constant $C$ depending polynomially on $\|\bE\|_{W^{2,\infty}}$,  $\|B^{-1}\|_{W^{2,\infty}}$ and $\|\eDpar\|_{W^{3,\infty}}$ such that if $f^\eps$ solves \eqref{eq:vlasov} with initial data a nonnegative density $f_0$, then $F^\eps$ defined so that $F^\eps(t,\cdot)$ is the push-forward of $f^\eps(t,\cdot)$ by the map $(\bfx,\bfv)\mapsto (\GC,\vGC,\muGC)(t,\bfx,\bfv)$ satisfies for a.e. $t\geq0$
$$
\|F^\eps(t,\cdot)-G^\eps(t,\cdot)\|_{\dot{W}^{-1,1}}\leq C\,\eps^2\,e^{C\,t^4\,(1+\eps\,t)}\,\int_{\R^6} e^{C\,t\,\|\bfv\|^3\,(1+\eps\,\|\bfv\|)}\,(1+\|\bfv\|^4)\,f_0(\bfx,\bfv)\,\dD\bfx\,\dD\bfv,
$$
where $G^\eps$ solves 
\be
\label{eq:2nd-mu}
\d_t G^\eps\,+\,\Div_{\bfZ} \left(\cY^\eps\,G^\eps\right)
\,=\,0,
\ee
with initial data $G_{0}^\eps=F^\eps(0,\cdot)$.
\et


\bibliographystyle{apalike}
\bibliography{3D-asymptotics}

\begin{thebibliography}{}

\bibitem[Bellan, 2008]{bellan_2006_fundamentals}
Bellan, P.~M. (2008).
\newblock {\em Fundamentals of plasma physics}.
\newblock Cambridge University Press.

\bibitem[Benettin and Sempio, 1994]{Benettin-Sempio}
Benettin, G. and Sempio, P. (1994).
\newblock Adiabatic invariants and trapping of a point charge in a strong
  nonuniform magnetic field.
\newblock {\em Nonlinearity}, 7(1):281--303.

\bibitem[Bogoliubov and Mitropolsky, 1961]{Bogoliubov-Mitropolsky_oscillations}
Bogoliubov, N.~N. and Mitropolsky, Y.~A. (1961).
\newblock {\em Asymptotic methods in the theory of non-linear oscillations}.
\newblock Translated from the second revised Russian edition. International
  Monographs on Advanced Mathematics and Physics. Hindustan Publishing Corp.,
  Delhi, Gordon and Breach Science Publishers, New York.

\bibitem[Bostan, 2010a]{bostan_10}
Bostan, M. (2010a).
\newblock Gyrokinetic {V}lasov equation in three dimensional setting. {S}econd
  order approximation.
\newblock {\em Multiscale Model. Simul.}, 8(5):1923--1957.

\bibitem[Bostan, 2010b]{Bostan_transport}
Bostan, M. (2010b).
\newblock Transport equations with disparate advection fields. {A}pplication to
  the gyrokinetic models in plasma physics.
\newblock {\em J. Differential Equations}, 249(7):1620--1663.

\bibitem[Bostan, 2019]{Bostan_2D-VP}
Bostan, M. (2019).
\newblock Asymptotic behavior for the {V}lasov-{P}oisson equations with strong
  external magnetic field. {S}traight magnetic field lines.
\newblock {\em SIAM J. Math. Anal.}, 51(3):2713--2747.

\bibitem[Brizard and Hahm, 2007]{bri_hahm_07}
Brizard, A.~J. and Hahm, T.~S. (2007).
\newblock Foundations of nonlinear gyrokinetic theory.
\newblock {\em Rev. Modern Phys.}, 79(2):421--468.

\bibitem[Chen, 2016]{chen_introduction}
Chen, F. (2016).
\newblock {\em Introduction to Plasma Physics and Controlled Fusion}.
\newblock Springer, third edition.

\bibitem[Cheverry, 2017]{cheve2}
Cheverry, C. (2017).
\newblock Anomalous transport.
\newblock {\em J. Differential Equations}, 262(3):2987--3033.

\bibitem[Degond and Filbet, 2016]{PDFF}
Degond, P. and Filbet, F. (2016).
\newblock On the {A}symptotic {L}imit of the {T}hree {D}imensional
  {V}lasov--{P}oisson {S}ystem for {L}arge {M}agnetic {F}ield: {F}ormal
  {D}erivation.
\newblock {\em J. Stat. Phys.}, 165(4):765--784.

\bibitem[Filbet and Rodrigues, 2016]{FR1}
Filbet, F. and Rodrigues, L.~M. (2016).
\newblock Asymptotically stable particle-in-cell methods for the
  {V}lasov-{P}oisson system with a strong external magnetic field.
\newblock {\em SIAM J. Numer. Anal.}, 54(2):1120--1146.

\bibitem[Filbet and Rodrigues, 2017]{FR2}
Filbet, F. and Rodrigues, L.~M. (2017).
\newblock Asymptotically preserving particle-in-cell methods for inhomogeneous
  strongly magnetized plasmas.
\newblock {\em SIAM J. Numer. Anal.}, 55(5):2416--2443.

\bibitem[Freidberg, 2008]{freidberg2008plasma}
Freidberg, J. (2008).
\newblock {\em Plasma Physics and Fusion Energy}.
\newblock Cambridge University Press.

\bibitem[Fr\'enod and Lutz, 2014]{FrenodLutz_geometrical_gyro-kinetic}
Fr\'enod, E. and Lutz, M. (2014).
\newblock On the geometrical gyro-kinetic theory.
\newblock {\em Kinet. Relat. Models}, 7(4):621--659.

\bibitem[Fr\'enod and Sonnendr\"ucker, 1998]{fre_son_97}
Fr\'enod, E. and Sonnendr\"ucker, {\'E}. (1998).
\newblock Homogenization of the {V}lasov equation and of the {V}lasov-{P}oisson
  system with a strong external magnetic field.
\newblock {\em Asymptot. Anal.}, 18(3-4):193--213.

\bibitem[Fr\'enod and Sonnendr\"ucker, 2000]{fre_son_98}
Fr\'enod, E. and Sonnendr\"ucker, {\'E}. (2000).
\newblock Long time behavior of the two-dimensional {V}lasov equation with a
  strong external magnetic field.
\newblock {\em Math. Models Methods Appl. Sci.}, 10(4):539--553.

\bibitem[Garbet et~al., 2010]{Garbet-et-al_2010}
Garbet, X., Idomura, Y., Villard, L., and Watanabe, T.~H. (2010).
\newblock Gyrokinetic simulations of turbulent transport.
\newblock {\em Nuclear Fusion}, 50:043002.

\bibitem[Golse and Saint-Raymond, 1999]{gol_lsr_99}
Golse, F. and Saint-Raymond, L. (1999).
\newblock The {V}lasov-{P}oisson system with strong magnetic field.
\newblock {\em J. Math. Pures Appl. (9)}, 78(8):791--817.

\bibitem[Han-Kwan, 2011]{HanKwan_PhD}
Han-Kwan, D. (2011).
\newblock {\em Contribution {\`a} l'{\'e}tude math{\'e}matique des plasmas
  fortement magn{\'e}tis{\'e}s}.
\newblock PhD thesis, Universit{\'e} Pierre et Marie Curie-Paris VI.

\bibitem[Hazeltine and Meiss, 2003]{haz_mei_03}
Hazeltine, R. and Meiss, J. (2003).
\newblock {\em Plasma Confinement}.
\newblock Dover Publications.

\bibitem[Herda, 2016]{herda_2016_massless}
Herda, M. (2016).
\newblock On massless electron limit for a multispecies kinetic system with
  external magnetic field.
\newblock {\em J. Differential Equations}, 260(11):7861--7891.

\bibitem[Herda, 2017]{Herda_PhD}
Herda, M. (2017).
\newblock {\em Analyse asymptotique et num{\'e}rique de quelques mod{\`e}les
  pour le transport de particules charg{\'e}es}.
\newblock PhD thesis, Universit{\'e} Claude Bernard Lyon 1.

\bibitem[Herda and Rodrigues, 2016]{herda_2016_anisotropic}
Herda, M. and Rodrigues, L.~M. (2016).
\newblock Anisotropic {B}oltzmann-{G}ibbs dynamics of strongly magnetized
  {V}lasov-{F}okker-{P}lanck equations.
\newblock {\em arXiv preprint arXiv:1610.05138}.

\bibitem[Krommes, 2012]{Krommes}
Krommes, J.~A. (2012).
\newblock The gyrokinetic description of microturbulence in magnetized plasmas.
\newblock In {\em Annual review of fluid mechanics. {V}olume 44, 2012},
  volume~44 of {\em Annu. Rev. Fluid Mech.}, pages 175--201. Annual Reviews,
  Palo Alto, CA.

\bibitem[Lee, 1983]{Lee}
Lee, W. (1983).
\newblock Gyrokinetic approach in particle simulation.
\newblock {\em Phys. Fluids}, 26(2):556--562.

\bibitem[Li, 2019]{Li_Kato-Ponce}
Li, D. (2019).
\newblock On {K}ato-{P}once and fractional {L}eibniz.
\newblock {\em Rev. Mat. Iberoam.}, 35(1):23--100.

\bibitem[Littlejohn, 1979]{littleJ1}
Littlejohn, R.~G. (1979).
\newblock A guiding center hamiltonian : A new approach.
\newblock {\em J. Math. Phys.}, 20:2445--2458.

\bibitem[Littlejohn, 1981]{littleJ2}
Littlejohn, R.~G. (1981).
\newblock Hamiltonian formulation of guiding center motion.
\newblock {\em Phys. Fluids}, 24:1730--1749.

\bibitem[Littlejohn, 1983]{littleJ3}
Littlejohn, R.~G. (1983).
\newblock Variational principles of guiding center motion.
\newblock {\em J. Plasma Physics}, 29:111--124.

\bibitem[Lutz, 2013]{Lutz_PhD}
Lutz, M. (2013).
\newblock {\em {\'E}tude math{\'e}matique et num{\'e}rique d'un mod{\`e}le
  gyrocin{\'e}tique incluant des effets {\'e}lectromagn{\'e}tiques pour la
  simulation d'un plasma de Tokamak}.
\newblock PhD thesis, Universit{\'e} de Strasbourg.

\bibitem[Matteo, 2017]{Matteo-PhD}
Matteo, V.~F. (2017).
\newblock {\em Gyrokinetic theory for particle transport in fusion plasmas}.
\newblock PhD thesis, Universit\`a di Roma Tre.

\bibitem[{Miot}, 2016]{Miot-2D-gyrokinetic}
{Miot}, {\'E}. (2016).
\newblock {On the gyrokinetic limit for the two-dimensional Vlasov-Poisson
  system}.
\newblock {\em arXiv preprint}, arXiv-1603.04502.

\bibitem[Miyamoto, 2006]{miyamoto_2006_plasma}
Miyamoto, K. (2006).
\newblock {\em Plasma physics and controlled nuclear fusion}, volume~38 of {\em
  Springer Series on Atomic, Optical, and Plasma Physics}.
\newblock Springer-Verlag Berlin-Heidelberg.

\bibitem[Piel, 2010]{piel2010plasma}
Piel, A. (2010).
\newblock {\em Plasma Physics: An Introduction to Laboratory, Space, and Fusion
  Plasmas}.
\newblock Springer Berlin Heidelberg.

\bibitem[Possanner, 2018]{Possanner}
Possanner, S. (2018).
\newblock Gyrokinetics from variational averaging: {E}xistence and error
  bounds.
\newblock {\em J. Math. Phys.}, 59(8):082702, 34.

\bibitem[Saint-Raymond, 2002]{laure0}
Saint-Raymond, L. (2002).
\newblock Control of large velocities in the two-dimensional gyrokinetic
  approximation.
\newblock {\em J. Math. Pures Appl. (9)}, 81(4):379--399.

\bibitem[Sanders et~al., 2007]{Sanders-Verhulst-Murdock_averaging}
Sanders, J.~A., Verhulst, F., and Murdock, J. (2007).
\newblock {\em Averaging methods in nonlinear dynamical systems}, volume~59 of
  {\em Applied Mathematical Sciences}.
\newblock Springer, New York, second edition.

\bibitem[{Scott}, 2017]{Scott_gyrokinetic}
{Scott}, B.~D. (2017).
\newblock {Gyrokinetic Field Theory as a Gauge Transform or: gyrokinetic theory
  without Lie transforms}.
\newblock {\em arXiv preprint}, arXiv-1708.06265.

\end{thebibliography}

\end{document}